\RequirePackage[english]{babel}
\documentclass[11pt]{amsart}
\usepackage[latin1]{inputenc}
\usepackage[T1]{fontenc}
\usepackage{ae,aecompl}
\usepackage{epsfig}
\usepackage{color}
\psfigdriver{dvips}

\title[Convex Foliated Projective Structures]{Convex Foliated
  Projective Structures and the Hitchin component for $\PSL_4(\RR)$}
\author{Olivier Guichard}
\address{CNRS, Laboratoire de Mathématiques d'Orsay, Orsay cedex, F-91405\\
Université Paris-Sud, Orsay cedex, F-91405}
\author{Anna Wienhard}
\address{Department of Mathematics, University of Chicago\\
5734 University Avenue, Chicago, IL 60637-1514 USA}
\thanks{A.W. was partially supported by the National Science
  Foundation under agreement No. DMS-0111298 and
  No. DMS-0604665. }
\date{\today}

\newtheorem{thm}{Theorem}[section]
\newtheorem{prop}[thm]{Proposition}
\newtheorem{lem}[thm]{Lemma}
\newtheorem{cor}[thm]{Corollary}
\newtheorem{remark}[thm]{Remark}
\newtheorem*{thm*}{Theorem}
\newtheorem{defi}[thm]{Definition}
\newtheorem{fact}[thm]{Fact}
\newtheorem{nota}[thm]{Notation}

\newtheorem{prop_intro}{Proposition}

\newtheorem{thm_intro}[prop_intro]{Theorem}

\newtheorem{cor_intro}[prop_intro]{Corollary}

\newcommand{\NN}{\mathbf{N}}
\newcommand{\ZZ}{\mathbf{Z}}
\newcommand{\RR}{\mathbf{R}}
\newcommand{\CC}{\mathbf{C}}

\newcommand{\PP}{\mathbb{P}}

\newcommand{\HH}{\mathbb{H}}

\newcommand{\G}{\Gamma}
\newcommand{\g}{\gamma}

\newcommand{\Flag}{\mathcal{F}lag}

\newcommand{\Ll}{\mathcal{L}}
\newcommand{\Gr}{\textup{Gr}}
\newcommand{\dev}{\textup{dev}}
\newcommand{\hol}{\textup{hol}}
\newcommand{\PSL}{\textup{PSL}}
\newcommand{\SL}{\textup{SL}}

\newcommand{\SO}{\textup{SO}}
\newcommand{\PGL}{\textup{PGL}}

\newcommand{\PSp}{\textup{PSp}}

\newcommand{\GL}{\textup{GL}}

\newcommand{\rep}{\textup{Rep}}

\newcommand{\id}{\textup{Id}}

\newcommand{\PkR}[1]{{\PP}^{#1}(\RR)}
\newcommand{\PkC}[1]{{\PP}^{#1}(\CC)}
\newcommand{\PnR}{\PkR{n}}
\newcommand{\PTR}{\PkR{3}}
\newcommand{\POC}{\PkC{1}}
\newcommand{\POR}{\PkR{1}}
\newcommand{\PTRd}{\PP^{3}( \RR)^{*}}
\newcommand{\PkRd}[1]{\PP^{#1}( \RR)^{*}}

\newcommand{\vfi}{\varphi}
\renewcommand{\hom}{\mathrm{Hom}}

\newcommand{\bqn}{\begin{eqnarray*}}
\newcommand{\eqn}{\end{eqnarray*}}
\newcommand{\bq}{\begin{eqnarray}}
\newcommand{\eq}{\end{eqnarray}}

\begin{document}
\begin{abstract}
  In this article we give a geometric interpretation of the Hitchin
  component $\mathcal{T}^4(\Sigma) \subset \rep(\pi_1(\Sigma),
  \PSL_4(\RR))$ of a closed oriented surface of genus $g\geq 2$. We
  show that representations in $\mathcal{T}^4(\Sigma)$ are precisely
  the holonomy representations of properly convex foliated projective
  structures on the unit tangent bundle of $\Sigma$.  From this we
  also deduce a geometric description of the Hitchin component
  $\mathcal{T}(\Sigma, \textup{Sp}_4(\RR))$ of representations into
  the symplectic group.
\end{abstract}

\maketitle

\section*{Introduction}
\label{sec_intro}
In his article [11] N.~Hitchin discovered a special connected
component, the ``Teichmüller component'', of the representation variety
of the fundamental group of a closed Riemann surface $\Sigma$ into a
simple adjoint $\RR$-split Lie group $G$. He showed that the
Teichm\"uller component, now usually called ``Hitchin component'', is
diffeomorphic to a ball of dimension $|\chi(\Sigma )| \dim(G)$.  In
this article we interpret the Hitchin component for $G =\PSL_4(\RR)$
as moduli space of certain locally homogeneous geometric structures.
Our main result is the following
\begin{thm_intro}
  \label{intro:main_thm}
  The Hitchin component for $\PSL_4(\RR)$ is naturally homeomorphic to
  the moduli space of (marked) properly convex foliated projective
  structures on the unit tangent bundle of $\Sigma$.
\end{thm_intro}

Let us briefly describe these geometric structures (see
Section~\ref{sec_defPCFP} for a precise definition). Properly convex
foliated projective structures are locally homogeneous $( \PGL_4( \RR),
\PTR)$-structures on the unit tangent bundle $M=S \Sigma$ of the
surface $\Sigma$ satisfying the following additional conditions:
\begin{itemize}
\item every orbit of the geodesic flow on $M$ is locally a projective
  line,
\item every (weakly) stable leaf of the geodesic flow is locally a projective
  plane and the projective structure on the leaf obtained by restriction
  is convex.
\end{itemize}

There is a natural map from the moduli space of projective structures
to the variety of representation $\pi_1(M) \to \PGL_4(\RR)$. We show
that the restriction of this map to the moduli space of properly
convex foliated projective structures is a  omeomorphism onto the
Hitchin component; in particular, the holonomy representation of a
properly convex foliated projective structure factors through the
projection $\pi_1(M) \rightarrow \pi_1( \Sigma)$ and takes values in
$\PSL_4(\RR)$.

Appealing to N.~Hitchin's result, we conclude

\begin{cor_intro}
  The moduli space of (marked) properly convex foliated projective
  structures on the unit tangent bundle of $\Sigma$ is a ball of
  dimension $-15\chi(\Sigma)$.
\end{cor_intro}

We describe (see Section~\ref{sec_examples}) several examples of
projective structures on $M$, including families of projective
structures with ``quasi-Fuchsian'' holonomy $\pi_1(M) \to \pi_1(\Sigma)
\to \PSL_2(\CC) \to \PSL_4(\RR)$. Those examples show that for the holonomy
representation to lie in the Hitchin component the
above additional conditions cannot be weakened.

Geometric interpretations of the Hitchin component were previously
known when $G$ is $\PSL_2(\RR)$ or $\PSL_3(\RR)$. For $\PSL_2(\RR)$
the Hitchin component is the image of the embedding of the
Fricke-Teichmüller space of $\Sigma$ into $\rep(\pi_1(\Sigma ),
\PSL_2(\RR))$ obtained by associating to a (marked) hyperbolic
structure its holonomy homomorphism. For $\PSL_3(\RR)$, S.~Choi and
W.~Goldman showed in \cite{GoldmanChoi_1993} that the Hitchin
component is homeomorphic to the moduli space of (marked) convex real
projective structures on $\Sigma$ . In both cases it can be proved
directly that these moduli spaces of geometric structures are balls of
the expected dimension (see \cite{Goldman_1990} for the
$\PSL_3$-case).

As a consequence of Theorem~\ref{intro:main_thm} we obtain a similar
description for the Hitchin component of the symplectic group. The
symplectic form on $\RR^4$ induces a natural contact structure on
$\PTR$ and $\PSp_4(\RR)$ is the maximal subgroup of $\PSL_4(\RR)$
preserving this contact structure. We call a locally homogeneous
structure modeled on $\PTR$ with this contact structure a projective
contact structure.

\begin{thm_intro}
  \label{intro:symp_thm}
  The Hitchin component for $\PSp_4(\RR)$ is naturally homeomorphic to
  the moduli space of (marked) properly convex foliated projective
  contact structures on the unit tangent bundle $M = S\Sigma$.
\end{thm_intro}

The symplectic form on $\RR^4$ gives rise to an involution on the
representation variety $\rep(\pi_1(\Sigma ), \PSL_4(\RR))$ and on the
Hitchin component for $\PSL_4(\RR)$. The set of fixed points of this
involution is the Hitchin component for $\PSp_4(\RR)$.

As above we conclude
\begin{cor_intro}
  The moduli space of (marked) properly convex foliated projective
  contact structures on the unit tangent bundle $M = S\Sigma$ is a
  ball of dimension $-10\chi(\Sigma)$.
\end{cor_intro}

\smallskip

Our results rely on recent progress in understanding representations
in the Hitchin component for $\PSL_n(\RR)$ made by F.~Labourie
\cite{Labourie_2003}. In particular, he proved that all
representations in the Hitchin components are faithful, discrete and
semisimple. V.~Fock and A.~Goncharov \cite{Fock_Goncharov} showed that
representations in the Hitchin component for a general simple
$\RR$-split Lie group have the same properties, using among other
things the positivity theory for Lie groups developed by G.~Lustzig.
More important for our work is that F.~Labourie (supplemented by the
first author \cite{Guichard_2006}) gives the following geometric
characterization of representations inside the Hitchin component of
$\PSL_n(\RR)$.
\begin{thm_intro}[Labourie \cite{Labourie_2003}, Guichard
  \cite{Guichard_2006}]
  \label{intro:thm_hiteqconvex}
  A representation $\rho: \pi_1(\Sigma) \to \PSL_n(\RR)$ lies in the
  Hitchin component if, and only if, there exists a continuous
  $\rho$-equivariant convex curve $\xi^1: \partial \pi_1(\Sigma) \to
  \PkR{n-1}$.
\end{thm_intro}

A curve $\xi^1: \partial \pi_1(\Sigma) \to \PkR{n-1}$ is said to be
\emph{convex} (see Definition~\ref{def_curvconv}) if for every
$n$-tuple of pairwise distinct points in $\partial \pi_1(\Sigma )$ the
corresponding lines are in direct sum. Convex curves into $\PkR{2}$
are exactly injective maps whose image bounds a strictly convex domain
in $\PkR{2}$.

It is easy to prove that the existence of such a curve for
$\PSL_2(\RR)$ implies that the representation is in the Teichmüller
space (see Lemma~\ref{lem_curveteich}).

Let us quickly indicate that Theorem~\ref{intro:thm_hiteqconvex} is
equivalent to the result of S.~Choi and W.~Goldman
\cite{GoldmanChoi_1993} that the representations in the Hitchin
component for $\PSL_3(\RR)$ are precisely the holonomy representations
of (marked) convex real projective structure on $\Sigma$.

A (marked) convex real projective structure on $\Sigma$ is a pair $(N,
f )$, where $N$ is a convex real projective manifold, that is $N$ is the
quotient $\Omega /\Gamma$ of a strictly convex domain $\Omega$ in
$\PkR{2}$ by a discrete subgroup $\Gamma$ of $\PSL_3(\RR)$, and $f :
\Sigma \rightarrow N$ is a diffeomorphism.

Given a representation $\rho: \pi_1(\Sigma_g) \to \PSL_3(\RR)$ in the
Hitchin component for $\PSL_3(\RR)$, let $\Omega_\xi \subset \PkR{2}$
be the strictly convex domain bounded by the convex curve
$\xi^1(\partial \pi_1(\Sigma)) \subset \PkR{2}$. Then
$\rho(\pi_1(\Sigma))$ is a discrete subgroup of the group of Hilbert
isometries of $\Omega_\xi$ and hence acts freely and properly
discontinuously on $\Omega_\xi$. The quotient
$\Omega_\xi/\rho(\pi_1(\Sigma))$ is a real projective convex manifold,
diffeomorphic to $\Sigma$. Conversely given a real projective
structure on $\Sigma$, we can $\rho$-equivariantly identify $\partial
\pi_1(\Sigma)$ with the boundary of $\Omega$ and get a convex curve
$\xi^1: \partial \pi_1(\Sigma) \to \partial \Omega \subset \PkR{2}$.

\smallskip

Our starting point to associate a geometric structure to a
representation $\rho: \pi_1(\Sigma) \to \PSL_4(\RR)$ in the Hitchin
component was to look for domains of discontinuity for the action of
$\rho(\pi_1(\Sigma))$ in $\PkR{n-1}$, similar to $\Omega_\xi$ for
$\PSL_3(\RR)$. It turns out that for this it is useful to consider the
convex curve $\xi^{n-1}: \partial \pi_1(\Sigma) \to \Gr_{n-1}^n(\RR)
\simeq \PkRd{n-1}$ associated to the contragredient representation of
$\rho$.  For example, the set $\Lambda \subset \PkR{n-1}$ of lines
which are contained in $n-1$ distinct hyperplanes $\xi^{n-1}(t)$ is
invariant and the action of $\Gamma$ on it is free and properly
discontinuous when $n \geq 4$, but cocompact only if $n = 4$.  The
open set $\Omega= \PkR{n-1} -\overline{\Lambda}$, which coincides with
$\Omega_\xi$ for $n=3$, turns out to have the right topology and
cohomological dimension in order to admit a cocompact action of
$\rho(\pi_1(\Sigma))$. Unfortunately, the action on $\Omega$ is proper
only for $n=3$ and $n=4$ ($\Omega$ is empty for $n = 2$).

When $n = 4$, the case of our concern, the action of
$\rho(\pi_1(\Sigma ))$ on $\PTR$ is proper (see
Paragraph~\ref{sec_domdisc}) precisely on the complement of the ruled
surface given by the union of projective lines tangents to the curve
$\xi^1(\partial \pi_1(\Sigma ))$. This complement is an open set which
has two connected components, namely $\Omega$ and $\Lambda$.  The
quotient $\Omega / \rho(\pi_1(\Sigma ))$ is a projective manifold
homeomorphic to the unit tangent bundle $M$ of $\Sigma$ and induces a
properly convex foliated projective structure on $M$. The quotient
$\Lambda / \rho( \pi_1(\Sigma ))$ is a projective manifold which is
homeomorphic to a quotient of $M$ by $\ZZ/3\ZZ$ and induces a
projective structure on $M$ which is foliated, but not properly
convex.

The construction of these domains of discontinuity gives a map from
the Hitchin component to the moduli space of properly convex foliated
structures on $M$. Conversely, starting with a properly convex
foliated projective structure on $M$ we construct an equivariant
convex curve and show that the projective structure is obtained by the 
above construction. This gives a more precise description of the
content of Theorem~\ref{intro:main_thm}.

\smallskip

Let us conclude with some open questions.  One might try to describe
how properly convex foliated projective structures on the unit tangent
bundle of $\Sigma$ can be glued from or decomposed into simpler pieces
similar to pair of pants decompositions of a hyperbolic surface. This
would probably lead to some generalized Fenchel-Nielsen coordinates
for the Hitchin component for $\PSL_4(\RR)$ and hopefully to a
geometric proof that the Hitchin component is a ball. Gluing convex
projective manifolds (with boundary) is one of the tools in
W.~Goldman's work \cite{Goldman_1990} and was also used in a recent
work of M.~Kapovich \cite{Kapovich_2006} to produce convex projective
structures on the manifolds, constructed by M.~Gromov and W.~Thurston
in \cite{Gromov_Thurston}, that admit Riemannian metrics of pinched
negative sectional curvature but no metrics of constant sectional
curvature. This approach might be extendable to the construction of
convex foliated projective manifolds, where the first and subtle point
is to find natural submanifolds at which one should cut and glue.

The interpretation of the Hitchin component for $\PSL_3(\RR)$ as
holonomy representations of convex real projective structures plays an
important role in independent work of F.~Labourie \cite{Labourie_2006}
and J.~Loftin \cite{Loftin_2001} who associate to a convex real
projective structure on $\Sigma$ a complex structure and a cubic
differential on $\Sigma$.  One might hope that our interpretation of
the Hitchin component for $\PSL_4(\RR)$ as holonomy representations of
properly convex foliated projective structures on $M$ could help to
associate a complex structure, a cubic and a quartic differential on
$\Sigma$ to every such representation.  This would answer, for
$\PSL_4( \RR)$, a conjecture of F.~Labourie describing the quotient of
the Hitchin component by the modular group.

Turning to higher dimensions, one might suspect that there is a moduli
space of suitable geometric structures associated to the Hitchin
component for $\PSL_n(\RR)$. As we alluded to above, using the convex
curve, natural domains of discontinuity for a convex representation
can be described, but in general none will admit a cocompact action of
$\pi_1( \Sigma)$.  To find the right geometric structures on a
suitable object associated to $\Sigma$ for general $n$ seems to be a
delicate and challenging problem.

In this paper, we concentrate on the projective manifold
$\Omega/\rho(\pi_1(\Sigma))$ obtained from one of the connected
components of the domain of discontinuity for a representation in the
Hitchin component for $\PSL_4( \RR)$. The projective manifolds
$\Lambda/\rho(\pi_1(\Sigma))$ obtained from the other connected
component also satisfy some additional properties, and it should be
possible to obtain a similar description for the Hitchin component
from them.

In a subsequent paper, we will concentrate more specifically on the
question of describing the Hitchin component for $\PSp_4( \RR)$ by
geometric structures modeled on the space of Lagrangians rather than
on projective space. Using the isomorphism $\PSp_4( \RR) =
\SO^\circ(2,3)$ Hitchin representations give rise to flat conformal
Lorentzian structures on $M$. Since the Hitchin component of
$\SO^\circ(2,3)$ embeds into the Hitchin component of $\PSL_5(\RR)$
this might also help to understand the geometric picture for
$\PSL_5(\RR)$.

\smallskip

{\bf Structure of the paper:} In Section~\ref{sec_conventions} we
review some classical facts about the unit tangent bundle $M= S\Sigma$
and fix some conventions and notation which are used throughout the
paper. Properly convex foliated projective structures are introduced
in Section~\ref{sec_geomstruct}, where we also recall some basic facts
about locally homogeneous geometric structures. In
Section~\ref{sec_examples} we describe several examples of projective
structures on $M$ whose properties justify the definitions made in
Section~\ref{sec_geomstruct}. After reviewing the example of a
properly convex foliated projective structure given in
Section~\ref{sec_examples} we construct in
Section~\ref{sec_reptostruct} a properly convex foliated projective
structure on $M$ starting from the convex curve associated to a
representation in the Hitchin component. Conversely, in
Section~\ref{sec_proptorep} we construct an equivariant convex curve
starting from a properly convex foliated projective structure on $M$,
and show that every properly convex foliated projective structure on
$M$ arises by the construction given in Section~\ref{sec_reptostruct}.
The consequences for representations into the symplectic group are
discussed in Section~\ref{sec_selfdual}. In the
Appendix~\ref{sec_appendix} we collect some useful facts.

\smallskip

{\bf Acknowledgments.}  We thank Bill Goldman for many useful
discussions and for explaining to us the construction of
``quasi-Fuchsian'' real projective structures on the unit tangent
bundle of $\Sigma$.  Both authors thank the Institute for Advanced
Study, Princeton for its hospitality, the second author also enjoyed
the hospitality of the Institut des Hautes \'Etudes Scientifiques,
Bures-sur-Yvette while working on this project.

\section{Conventions}
\label{sec_conventions}

\begin{quote}
  In this section we review basic facts about the geometry of surfaces and
  their unit tangent bundle and introduce some notation.
  Everything is
  classical except maybe Section~\ref{sec_actiongamma}.
\end{quote}

\subsection{The Unit Tangent Bundle}
\label{sec_tangentbundle}
Let $\Sigma$ be a connected oriented closed surface of genus $g\geq 2$ and
$\Gamma = \pi_1( \Sigma, x)$ its fundamental group. We denote by
$\widetilde{\Sigma}$ the universal covering of $\Sigma$.
\begin{nota}
  We denote by $M$ the circle bundle associated to the tangent bundle of
  $\Sigma$. $M$ is homeomorphic to the unit tangent bundle of $\Sigma$ with
  respect to any Riemannian metric on $\Sigma$.
\end{nota}
The fundamental group $\overline\Gamma = \pi_1(M,m)$ is a central extension of
$\Gamma$

\begin{equation*}
  1 \longrightarrow \ZZ \longrightarrow \overline{ \Gamma}
  \overset{p}{ \longrightarrow} \Gamma \longrightarrow 1.
\end{equation*}

This central extension and the group $\Gamma$ admit the classical
presentations:
\begin{eqnarray*}
  \Gamma & = & \big \langle a_1, \dots, a_g, b_1, \dots, b_g \mid [ a_1,
  b_1] \cdots [a_g, b_g] \big\rangle \\
  \overline{ \Gamma} & = & \big \langle a_1, \dots, a_g, b_1, \dots, b_g, \tau \mid [ a_1,
  b_1] \cdots [a_g, b_g] \tau^{2g}, [a_i, \tau], [b_i, \tau] \big \rangle,
\end{eqnarray*}
where $[a,b]=a b a^{-1} b^{-1}$ is the commutator of two elements $a, b$.

There is an important covering $\overline{M} = S \widetilde{\Sigma}$, the unit
tangent bundle of $\widetilde{\Sigma}$ which is a Galois covering of $M$ with
covering group $\Gamma$. Topologically $S \widetilde{ \Sigma}$ is the product
$S^1 \times \RR^2$ so the universal cover $\widetilde{M} = \widetilde{
  \overline{M}}$ is a $3$-dimensional ball and the covering $\widetilde{M}
\to \overline{M}$ is an abelian covering with covering group $\ZZ$.

\subsubsection{Canonical Foliations}
\label{sec_canonfol}
We fix for a moment a hyperbolic metric on the surface $\Sigma$, that is a
Riemannian metric of constant sectional curvature $-1$.  The induced geodesic
flow $g_t: M \to M$ on the unit tangent bundle of $\Sigma$ is Anosov
(see \cite[\S~17.4--6]{Katok_1995}).
In particular, $M$ is endowed with two codimension one foliations,
namely the (weakly) stable foliation and the (weakly) unstable foliation of the geodesic flow,
and a codimension two foliation given by the flow lines, which we call the
geodesic foliation.
\begin{nota}
  We denote by $\mathcal{F}$ the set of leaves of the (weakly) stable foliation and by
  $\mathcal{G}$ the set of leaves of the geodesic foliation.
\end{nota}
The set $\mathcal{G}$ coincides with the set of (unparametrized) oriented
geodesics in $\Sigma$.  Correspondingly the lifts of the geodesic flow to
$\overline{M}$ and $\widetilde{M}$ induce foliations denoted
$\overline{\mathcal{F}}, \overline{\mathcal{G}}$ and respectively
$\widetilde{\mathcal{F}}, \widetilde{\mathcal{G}}$.
 \begin{nota}
   A typical element of $\mathcal{G}$ will be denoted by $g$ and a typical
   element of $\mathcal{F}$ by $f$, similarly for $\overline{\mathcal{F}}, \overline{\mathcal{G}}$ and
$\widetilde{\mathcal{F}}, \widetilde{\mathcal{G}}$. 
\end{nota}
The sets $\overline{\mathcal{F}}, \overline{\mathcal{G}}$ and
$\widetilde{\mathcal{F}}, \widetilde{\mathcal{G}}$ carry a natural
topology 
coming from the Hausdorff distance on subsets, and admit a natural action of
the corresponding covering groups $\Gamma$ and $\overline{\Gamma}$
respectively.

\subsubsection{A Topological Description of the Foliations }
\label{sec_topodesc}
The geodesic flow on $M$ depends on the choice of hyperbolic metric on
$\Sigma$, but it is well known that the foliations $\mathcal{F}$ and
$\mathcal{G}$ admit a description which shows that topologically they are indeed independent of the metric. We recall this description
briefly.

The group $\Gamma$ is a hyperbolic group and hence there is a canonical
boundary at infinity of $\Gamma$ (see \cite{GhysHarpe_1990} for definition and
properties).

\begin{nota}
  The boundary at infinity of the group $\Gamma$ is denoted by $\partial
  \Gamma$. It is a topological circle with a natural $\Gamma$-action.
\end{nota}

The dynamics of the action of any element $\gamma\in \Gamma - \{ 1\}$ is well
understood; the element $\gamma$ has exactly two fixed points $t_{+,\gamma},
t_{-,\gamma} \in \partial \Gamma$. For any $t\in \partial \Gamma$ distinct
from $t_{\mp,\gamma}$ we have $\lim_{n \to \pm \infty} \gamma^{n} \cdot t =
t_{\pm,\gamma}$.

Let us define $\partial \Gamma^{(2)} := \partial \Gamma^2 - \Delta$, where
$\Delta = \{ (t,t) \mid t \in \partial \Gamma \}$ is the diagonal in $\partial
\Gamma^2$. Then we have the following classical facts.

\begin{lem}
\label{lem_minactforGamma}
\label{lem_closedsetofgeodesics}
  The action of $\Gamma$ on $\partial \Gamma$ 
  is minimal.

  The subset of pairs of fixed points $\{ (t_{+, \gamma}, t_{-,\gamma}) 
  \mid \gamma \neq 1\}$ is dense in $\partial \Gamma^{(2)}$.
\end{lem}

The hyperbolic metric on $\Sigma$ isometrically identifies $\widetilde{
  \Sigma}$ with the hyperbolic plane $\HH^2$. Let $\iota: \Gamma \to
\PSL_2(\RR)= \textup{Isom}^{+}( \HH^2)$ be the homomorphism which makes this
identification equivariant. Such a homomorphism is faithful and
  discrete, and will be called a \emph{Fuchsian} representation or a
  uniformization. 
Then any orbital application $\Gamma \rightarrow
\HH^2$, $\gamma \mapsto \iota( \gamma) \cdot x_0$ is a quasiisometry, and
hence induces a $\Gamma$-equivariant homeomorphism $\partial \Gamma \overset{
  \sim}{ \rightarrow} \partial \HH^2$, which is easily seen to be independent
from the base point $x_0$. Using the upper half space model for $\HH^2$ the
boundary $\partial \HH^2$ is identified with the projective line $\POR$ with
the natural $\PSL_2( \RR)$ action.

The orientation of $\Sigma$ induces an orientation on $\partial \Gamma
\simeq \partial \widetilde{\Sigma} \simeq S_x\widetilde{\Sigma}$, for any $x$
in $\widetilde{ \Sigma}$.
This enables us say when a triple of pairwise distinct points of the boundary is positively
oriented.
\begin{nota}
  We denote by $\partial \Gamma^{3+}$ the subset of $\partial \Gamma^3$
  consisting of pairwise distinct positively oriented triples.
\end{nota}

\begin{figure}[htbp]
  \centering
  \input{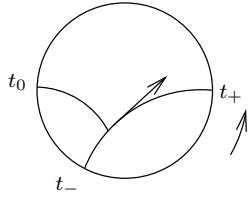}
  \caption{A positively oriented triple}
  \label{fig_orientriple}
\end{figure}

The unit tangent bundle $\overline{M}= S \widetilde{ \Sigma}$ can be
$\Gamma$-equivariantly identified with $\partial \Gamma^{3+}$
\begin{eqnarray*}
  \overline{M} & \longrightarrow & \partial \Gamma^{3+}\\
  v &\longmapsto& (t_+, t_0, t_-),
\end{eqnarray*}
where $t_+$ is the endpoint at
$+\infty$, $t_-$ is the endpoint at $-\infty$ of the geodesic $g_v$ defined by
$v$, and $t_0$ is the unique endpoint in $\partial \Gamma\cong S^1$ of the
geodesic perpendicular to $g_v$ at the foot point of $v$ such that $(t_+, t_0,
t_-)$ is positively oriented (see Figure~\ref{fig_orientriple}).

In this model the leaves of the geodesic and (weakly) stable foliations
$\overline{\mathcal G}$, $\overline{\mathcal F}$ of $\overline{M}$ through the
point $v = (t_+, t_0, t_-)$ are explicitly given by
\begin{eqnarray*}
  \overline{g}_v & = & \big\{ (s_+, s_0, s_-) \in \partial\Gamma^{3+} \mid s_+
  = t_+, s_- = t_- \big \}\\
  \text{and } \overline{f}_v & = & \big\{ (s_+, s_0, s_-) \in
  \partial\Gamma^{3+} \mid s_+ = t_+ \big \}.
\end{eqnarray*}
In particular, the set of
leaves of the geodesic foliation $\overline{ \mathcal{G}}$ on $\overline{M}$
is $\Gamma$-equiv\-ari\-antly identified with $\partial \Gamma^{(2)}$, where
the oriented geodesic $\overline{g}_v$ through $v=(t_+, t_0, t_-)$ is
identified with $(t_+, t_-) \in \partial \Gamma^{(2)} $ The set of (weakly) stable
leaves $\overline{ \mathcal{F}}$ is $\Gamma$-equi\-var\-iant\-ly identified
with $\partial \Gamma$ by mapping the (weakly) stable leaf $\overline{f}_v$ to $t_+
\in \partial \Gamma$.

From now on, we will not distinguish anymore between a (weakly) stable leaf seen
\emph{as a subset of $\overline{ M}$} or as \emph{an element of $\overline{
    \mathcal{F}}$} or as \emph{an element of $\partial \Gamma$}. So for
example $m \in t \in \partial \Gamma$ will denote a point $m$ of $\overline{
  M}$ in the leaf of $\overline{ \mathcal{ F}} \simeq \partial \Gamma$
corresponding to $t$. We will use a similar language for geodesics leaves.

\subsubsection{Identifying $\overline{M}$ with $\PSL_2(\RR)$}
\label{sec_Masquotient}
We will sometimes consider $M$ as the quotient of $\PSL_2( \RR)$ by the
subgroup $\iota( \Gamma) < \PSL_2( \RR)$ and identify leaves of the foliations
${\mathcal G}$ and ${\mathcal F}$ with orbits of subgroups of $\PSL_2( \RR)$.

The isometry $ \widetilde{ \Sigma} \simeq \HH^2$, which is $\iota$-equivariant
for some $\iota: \Gamma \rightarrow \PSL_2( \RR)$, induces a diffeomorphism
$\overline{M}=S \widetilde{ \Sigma} \simeq S \HH^2$ of the unit tangent
bundles.  Since the action of $\PSL_2( \RR)$ on $S \HH^2$ is simply
transitive, we can identify $S \HH^2$ with $\PSL_2( \RR)$ and obtain a
diffeomorphism:
\begin{equation*}
  M = \Gamma \backslash \overline{M} \overset{ \sim}{ \longrightarrow} \iota(
  \Gamma) \backslash \PSL_2( \RR).
\end{equation*}
Note that $\PSL_2(\RR)$ acts by right multiplication on itself and hence also
on $\iota( \Gamma) \backslash \PSL_2( \RR)$.

\begin{lem}
  \label{lem_leavesorbits}
  Under these identifications $M \simeq \iota( \Gamma) \backslash \PSL_2(
  \RR)$ and $\overline{M} \simeq \PSL_2( \RR)$
  \begin{enumerate}
  \item the leaves of the geodesic foliation are the (right) orbits of the
    Cartan subgroup:
    \begin{equation*}
      A = \Big \{ \left( 
        \begin{array}{cc}
          e^{t/2} & 0 \\ 0 & e^{-t/2}
        \end{array} \right) \mid t \in \RR \Big \}, 
    \end{equation*}
  \item the leaves of the (weakly) stable foliation are the (right) orbits of the
    parabolic subgroup:
    \begin{equation*}
     P = \Big \{ \left( 
        \begin{array}{cc}
          a & b \\ 0 & c
        \end{array} \right) \in \PSL_2(\RR) \Big \}.
    \end{equation*}
  \end{enumerate}
\end{lem}

\begin{remark}
  \begin{itemize}
  \item The leaves of the (weakly) unstable foliation are the right orbits of $P^{opp}
    = \Big \{ \left(
        \begin{array}{cc}
          a & 0 \\ b & c
        \end{array} \right) \in \PSL_2(\RR) \Big \}$.
    \item The (left) action of $\PSL_2( \RR)$ on $\overline{M}$ is transitive
      on the set of leaves $\overline{ \mathcal{F}}$ and $\overline{
        \mathcal{G}}$.
    \item The above identification $M \simeq \iota( \Gamma) \backslash \PSL_2(
      \RR)$ endows $M$ with a locally homogeneous $( \PSL_2( \RR), \PSL_2(
      \RR))$-structure with $\PSL_2( \RR)$ acting by left multiplication on
      itself.
  \end{itemize}
\end{remark}

\subsection{The Action of $\overline{ \Gamma}$ on the Leaf Spaces}
\label{sec_actiongamma}
In this paragraph we establish some facts about the action of $\overline{
  \Gamma} = \pi_1( M)$ on the space of geodesics $\widetilde{ \mathcal{G}}$ of
$\widetilde{ M}$ and on the space of (weakly) stable leaves $\widetilde{ \mathcal{F}}$,
which we will use frequently.

Note first that there is an identification of $\widetilde{ \mathcal{F}}$ with
$\widetilde{ \partial \Gamma}$, the universal cover of $\partial \Gamma \simeq
S^1$, lifting the natural isomorphism $\overline{ \mathcal{F}} \simeq \partial
\Gamma$. Two such identifications $\widetilde{ \mathcal{F}} \simeq \widetilde{
  \partial \Gamma}$ differ only by a translation by an element of the central
subgroup $\langle \tau \rangle < \overline{ \Gamma}$. In particular, there is
a well defined action of $\overline{ \Gamma}$ on $\widetilde{ \partial
  \Gamma}$ making any of these isomorphisms $\widetilde{ \mathcal{F}} \simeq
\widetilde{ \partial \Gamma}$ equivariant.

The chosen orientation on $\partial \Gamma$ induces an orientation on
$\widetilde{ \partial \Gamma}$. We choose the element $\tau$ generating the
center of $\overline{ \Gamma}$ so that $( \tau^n \tilde{f}, \tau^m \tilde{f},
\tilde{f})$ is positively oriented precisely when $ n > m > 0$.

\subsubsection{Minimality}
\begin{lem}
  \label{lem_minstable}
  The action of $\overline{ \Gamma}$ on $\widetilde{ \mathcal{F}} \simeq
  \widetilde{ \partial \Gamma}$ is minimal.
\end{lem}

\begin{proof}
  Recall that the action of $\overline{\G}$ on $\widetilde{ \partial \Gamma}$
  is minimal if any closed $\overline{\G}$-invariant subset of $\widetilde{
    \partial \Gamma}$ is either empty or equal to $\widetilde{ \partial
    \Gamma}$. Any $\overline{\G}$-invariant subset $\widetilde{ A}$ of
  $\widetilde{ \partial \Gamma}$ is in particular $\langle \tau
  \rangle$-invariant. Hence it is of the form $\widetilde{ A}=\pi^{-1} (A)$
  where $\pi : \widetilde{ \partial \Gamma} \rightarrow \partial \Gamma$ is
  the natural projection and $A$ is a $\Gamma$-invariant subset of $\partial
  \Gamma$.  The set $\widetilde{ A}$ is closed if and only if $A$ is closed.
  Since $\G$ acts minimally on $\partial \Gamma$, any closed
  $\overline{\G}$-invariant subset $\widetilde{ A} \subset \widetilde{
    \partial \Gamma}$ is either $\pi^{-1}(\emptyset) = \emptyset$ or
  $\pi^{-1}(\partial\G)= \widetilde{ \partial \Gamma}$.
\end{proof}

The space of leaves of the (weakly) unstable foliation of the geodesic flow on
$\widetilde{ M}$ can also be identified with $\widetilde{ \partial \Gamma}$.
This enables us to identify $\widetilde{ \mathcal{G}}$ with a set of pairs $(
\widetilde{t}_+ , \widetilde{t}_- )$ in $\widetilde{ \partial \Gamma} \times
\widetilde{ \partial \Gamma}$. Since the identification is a lift of the
natural identification
\begin{center}
  $\displaystyle \overline{ \mathcal{G}} \simeq \partial \Gamma^{(2)} = \{
  (t_+, t_-) \in \partial \Gamma^2 \mid t_+ \neq t_- \},$
\end{center}
the set $\widetilde{ \mathcal{G}}$ will be identified with a subset of:
\begin{equation*}
  \{ ( \widetilde{t}_+ , \widetilde{t}_- ) \in \widetilde{ \partial \Gamma}^2
  \mid \pi( \widetilde{t}_+ ) \neq \pi( \widetilde{t}_- ) \} = \bigcup_{ n \in
  \ZZ} \widetilde{ \partial \Gamma}^{(2)}_{[n]},
\end{equation*}
where
\begin{equation*}
  \widetilde{ \partial \Gamma}^{(2)}_{[n]} := \{ ( \widetilde{t}_+ ,
  \widetilde{t}_- ) \in \widetilde{ \partial \Gamma}^2 \mid ( \tau^{n+1}
  \widetilde{t}_- , \widetilde{t}_+ , \tau^n \widetilde{t}_- ) \text{ is
    oriented}\}. 
\end{equation*}
This is in fact the decomposition of $\{ ( \widetilde{t}_+ , \widetilde{t}_- )
\in \widetilde{ \partial \Gamma}^2 \mid \pi( \widetilde{t}_+ ) \neq \pi(
\widetilde{t}_- ) \}$ into connected components.  As $\widetilde{
  \mathcal{G}}$ is connected, there exist an $n \in \ZZ$ such that
$\widetilde{ \mathcal{G}}$ is equal to $\widetilde{ \partial
  \Gamma}^{(2)}_{[n]}$ and by changing the identification of the set of
(weakly) unstable leaves with $\widetilde{ \partial \Gamma}$ we can suppose that $n=0$.

\begin{remark}
  The manifold $\widetilde{ M}$ can be $\overline{ \Gamma}$-equivariantly
  identified with the set of triples $( \tilde{t}_+, \tilde{t}_0,
  \tilde{t}_-)$ of $(\widetilde{ \partial \Gamma})^3$ where $(\tau
  \tilde{t}_-, \tilde{t}_+, \tilde{t}_0, \tilde{t}_-)$ is positively oriented.
\end{remark}

\subsubsection{Elements of Zero Translation}
\label{sec_noncent}
Any element $\overline{ \gamma}\in \overline{ \Gamma} - \langle \tau \rangle$
projects onto an element $\gamma = p(\overline{\g})\in \Gamma-\{1\}$. Every
element $\gamma\in \Gamma-\{1\}$ has a unique attractive and repulsive fixed
point $t_{+, \gamma}$, respectively $t_{-, \gamma}$ in $\partial \Gamma$,
which lift to two $\langle \tau \rangle$-orbits $ ( \tau^n \widetilde{ t}_{+,
  \gamma} )_{ n \in \ZZ}$ and $ ( \tau^n \widetilde{ t}_{-, \gamma} )_{ n \in
  \ZZ}$ in $\widetilde{ \partial \Gamma}$.  We choose a pair $( \widetilde{
  t}_{+, \gamma}, \widetilde{ t}_{-, \gamma})$ such that the triple $( \tau
\widetilde{ t}_{-, \gamma} \widetilde{ t}_{+, \gamma}, \widetilde{ t}_{-,
  \gamma})$ is oriented. Two different choices of such a pair differ only by
the action of a power of $\tau$.

Since the element $\overline{ \gamma}$ commutes with $\tau$, it leaves the two
orbits $ ( \tau^n \widetilde{ t}_{ \pm, \gamma} )_{ n \in \ZZ}$ invariant and
acts orientation preserving on $\widetilde{ \partial \Gamma}$. Therefore there
is a unique integer $l$ such that
\begin{equation*}
  \overline{ \gamma} \cdot \widetilde{ t}_{ +, \gamma} =
  \tau^l \widetilde{ t}_{ +, \gamma} \text{ and }
  \overline{ \gamma} \cdot \widetilde{ t}_{ -, \gamma} =
  \tau^l \widetilde{ t}_{ -, \gamma}.
\end{equation*}
We will call $l =: \mathbf{t} ( \overline{ \gamma})$ the \emph{translation} of
$\overline{ \gamma}$. Obviously $\mathbf{t}( \overline{ \gamma} \tau^m) =
\mathbf{t} (\overline{ \gamma}) + m$, hence in every orbit $\{ \overline{
  \gamma} \tau^m \mid m \in \ZZ \}$ there is a unique element of zero
translation. 

Elements in $\overline{\G}$ of zero translation can be characterized by
considering the action on the space of geodesics. Given a non-trivial element
$\gamma \in \Gamma$, it fixes exactly two geodesic leaves in $\overline{
  \mathcal{G}} \simeq \partial \Gamma^{(2)}$, namely $\overline{g_\g} = (
t_{+, \gamma}, t_{-, \gamma})$ and $\overline{g_{\g^{-1}}} = ( t_{-, \gamma},
t_{+, \gamma})$, and its action on $\overline{g_\g} = ( t_{+, \gamma}, t_{-,
  \gamma})$ corresponds to a positive time map of the geodesic flow. The
geodesic leave $\overline{g_\g}$ lifts to one $\langle \tau \rangle$-orbit of
geodesic leaves $( g_n)_{ n \in \ZZ}$ in $\widetilde{ \mathcal{G}}$.  Among
the lifts of $\gamma$ there is a unique element $\overline{ \gamma} \in
\overline{\G}$ fixing each of these geodesics This element $\overline{
  \gamma}$ is the unique element of translation zero in $\{ \overline{ \gamma}
\tau^m \mid m \in \ZZ \}$.

\begin{lem}
  \label{lem_closedgeo}
  The set of pairs of fixed points of zero translation elements
  \begin{equation*}
    \left \{ (\tilde{t}_{+, \gamma}, \tilde{t}_{-,\gamma}) \mid \gamma\in 
    \overline{\Gamma} - \{ 1\} \right \}
  \end{equation*}
  is dense in $ \widetilde{\partial\Gamma}_{[0]}^{(2)}$.
\end{lem}

\begin{proof}
  The closure of this set is $\tau$-invariant (since we take all possible 
  fixed pairs for a given element) so it is of the form $\pi^{-1}(A)$
  where $A \subset \partial \Gamma^{(2)}$ is closed and contains all the 
  pairs $(t_{+, \gamma}, t_{-, \gamma})$ with $\gamma \in \Gamma -\{ 1\}$. 
  We conclude by Lemma~\ref{lem_closedsetofgeodesics}.
\end{proof}

\subsection{Projective Geometry}
\label{sec_projgeom}
Since we will consider manifolds which are locally modeled on the
three-dimensional real projective space, we recall some basic notions
from projective geometry.

Let $E$ be a vector space, then $\mathbb{P}( E)$ denotes the space of
lines in $E$. We write $\PkR{n-1}$ when $E = \RR^n$. 
The Grassmanian of $m$-dimensional subspaces of $E$ is
denoted by $\textup{Gr}_m(E)$ or $\Gr_{m}^{n}(\RR)$ when $E =
\RR^n$. We denote by $\PkRd{n-1}$ the Grassmanian of hyperplanes in
$\RR^n$. The variety of full flags in $E$ is denoted by $\mathcal{F}lag(E)$.

If $F \subset E$ is a subvector space, $\mathbb{P}(F)$ is naturally a
subspace of $\mathbb{P}(E)$, which is called a \emph{projective line} when $\dim
F = 2$ and a \emph{projective plane} when $\dim F = 3$. We will
regularly consider an element of $\textup{Gr}_{2}^{4}( \RR)$, or of
$\PTRd$ as a projective line respectively a projective plane in
$\PTR$ without introducing any
additional
notation. What is meant should always be clear from the context.

\subsubsection{Convexity}
\label{sec_convexity}
A subset $C$ of the projective space is said to be \emph{convex} if its
intersection with any projective line is connected. It is said to be 
\emph{properly convex} if its closure does not contain any projective
line.

A convex subset of the projective plane which is a connected
component of the complement of two projective lines through some point $x$ is called a
\emph{sector}. The point $x$ will be called the \emph{tip} of the sector.

\section{Geometric Structures}
\label{sec_geomstruct}

\begin{quote}
  In this section we introduce the notion of properly convex foliated
  projective structures and define the Hitchin component. 
  For more background on geometric structures we refer the reader to
  \cite{Goldman_88}.
\end{quote}

\subsection{Projective Structures}
\label{sec_projstruc}

\begin{defi}
  A \emph{projective structure} on an $n$-dimensional manifold $M$ is a
  maximal atlas $\{(U, \vfi_U)\}$ on $M$ such that
  \begin{enumerate}
  \item $\{U\}$ is an open cover of $M$ and, for any $U$, $\vfi_U: U \to
    \vfi_U (U)$ is a homeomorphism onto an open subset of $\PnR$.
  \item For each $U, V$ the change of coordinates $\vfi_V \circ \vfi_{U}^{-1}:
    \vfi_U (U \cap V) \to \vfi_V(U \cap V)$ is locally projective, i.e. 
    it is (locally) the restriction of an element of $\PGL_{n+1} (\RR)$.
  \end{enumerate}
  A manifold endowed with a projective structure is called a $\PnR$-manifold.
\end{defi}
A projective structure is a locally homogeneous $(\PGL_{n+1}(\RR),
\PnR)$-structure.

Let $M,N$ be two $\PnR$-manifolds. A continuous map $h: M \to N$ is
\emph{projective} if for every coordinate chart $(U, \vfi_U)$ of $M$ and every
coordinate chart $(V, \vfi_V)$ of $N$ the composition
\begin{equation*}
  \vfi_V \circ h\circ \vfi_{U}^{-1}: \vfi_U (h^{-1}(V) \cap U)
  \longrightarrow \vfi_V (h(U) \cap V)
\end{equation*}
is locally projective. A projective map is always a local homeomorphism.

Two projective structures on $M$ are equivalent if there is a homeomorphism
$h: M\to M$ isotopic to the identity which is a projective map with respect to
the two projective structures on $M$.  The set of equivalence classes of
projective structures on $M$ is denoted by $\mathcal{P}(M)$. We draw the
reader's attention to the above restriction on the differentiability class for
the coordinate changes; the projective structures constructed in
Section~\ref{sec_reptostruct} are only continuous.

Let $\widetilde{M}$ be the universal covering of $M$. A projective structure
on $M$ defines a projective structure on $\widetilde{M}$. Since
$\widetilde{M}$ is simply connected, there exists a global projective map
$\dev: \widetilde{M} \to \PnR$. The action of $\pi_1( M)$ on $\widetilde{M}$
respects the projective structure on $\widetilde{M}$. More precisely, for
every $\gamma \in \pi_1( M)$, there is a unique element $\hol( \gamma)$ in
$\PGL_{n+1}( \RR)$ such that $\dev( \gamma \cdot m) = \hol( \gamma) \cdot
\dev( m)$ for any $m$ in $\widetilde{M}$. This defines a homomorphism $\hol:
\pi_1(M) \to \PGL_{n+1}(\RR)$ such that the map $\dev$ is $\hol$-equivariant.
The map $\dev$ is called the developing map and the homomorphism $\hol$ is
called the holonomy homomorphism.  Conversely the data of a developing pair
$(\dev,\hol)$ defines a $\pi_1(M)$-invariant projective structure on
$\widetilde{M}$, hence a projective structure on $M$:
\begin{prop}[\cite{Goldman_88}]
A projective structure on $M$ is equivalent to the data $(\dev, \hol)$
of a holonomy homomorphism $\hol: \pi_1(M) \to
\PGL_{n+1}(\RR)$ and a $\hol$-equivariant local homeomorphism  $\dev:
\widetilde{M} \to \PnR$.
Two pairs $(\dev_1, \hol_1)$ and $(\dev_2, \hol_2)$ define
equivalent projective structures on $M$ if and only if 
there exists a homeomorphism $h: M\to M$ isotopic to the identity and
an element $g \in \PGL_{n+1}(\RR)$ such that 
\begin{itemize}
\item $\hol_2(\gamma) = g \hol_1(\gamma) g^{-1}$ for all $\gamma \in
  \pi_1(M)$, and 
\item $\dev_1\circ \tilde{h} = g^{-1}\circ \dev_2$, where $\tilde{h}:
  \widetilde{M}\to \widetilde{M}$ is the homeomorphism induced by
  $h$. 
\end{itemize}
\end{prop}

This defines an equivalence relation on the pairs $(\dev, \hol)$ such
that $\mathcal{P}(M)$ is identified with the set of equivalence
classes of pairs $(\dev, \hol)$. 
We endow the set of pairs $(\dev, \hol)$ with the topology coming from
the compact-open topology on the spaces of maps $\widetilde{M} \to
\PnR$ and
$\pi_1(M) \to
\PGL_{n+1}(\RR)$, and consider $\mathcal{P}(M)$ with the induced
quotient 
topology. 

The holonomy map
\begin{equation*}
  \hol: \mathcal{P}(M) \longrightarrow \hom(\pi_1(M),\PGL_{n+1}(\RR))/\PGL_{n+1}(\RR)
\end{equation*}
associates to a pair $(\dev, \hol)$ just the holonomy
homomorphism. 

\subsection{Foliated Projective Structures on $S\Sigma$} 
\label{sec_folprostr}
We consider now the unit tangent bundle $M=S\Sigma$ endowed with the
(weakly) stable foliation $\mathcal{F}$ and the geodesic foliation
$\mathcal{G}$. 

\begin{defi}
  A \emph{foliated projective structure} on $(M,\mathcal{F}, \mathcal{G})$ is 
  a projective structure $\{(U, \vfi_U)\}$ on $M$ with the
  additional properties that 
  \begin{enumerate}
  \item for every $x\in U$ the image $\vfi_U( U \cap
    g_x)$ of the geodesic leaf $g_x$ through $x$  is contained in a
    projective line, and 
  \item  for every $x\in U$ the image $\vfi_U( U \cap
    f_x)$ of the (weakly) stable leaf $f_x$ through $x$  is contained in a
    projective plane.
  \end{enumerate}
\end{defi}

Two foliated projective structures on $(M,\mathcal{F}, \mathcal{G})$
are \emph{equivalent} if there is a projective homeomorphism $h: M \to M$
isotopic to the identity such that $h^*\mathcal{G} = \mathcal{G}$ and 
$h^*\mathcal{F} = \mathcal{F}$.
The space of equivalence classes of foliated projective
structures on $(M,\mathcal{F}, \mathcal{G})$ is denoted by 
$\mathcal{P}_f (M)$.

\begin{remark}
  \label{rem_folprostr}
  \begin{itemize}
  \item Note that the natural map of $\mathcal{P}_f (M)$
    into $\mathcal{P}(M)$ is not an inclusion since we do not only
    restrict to a subset of
    projective structures, but at the same time refine the equivalence
    relation. We do not know if this map is injective or
    not.
  \item It is a simple exercise to show that any homeomorphism $h$ of
    $M$ homotopic to the identity and respecting the foliation
    $\mathcal{G}$ sends every geodesic to itself. So it is of the form $m \mapsto
    \phi_{f(m)} (m)$ where $f:M \rightarrow \RR$ is a continuous
    function and $(\phi_t)_{t \in
      \RR}$ is the geodesic flow on $M$. Since $m\mapsto \phi_{f(m)} (m)$ is
    a homeomorphism, we obtain that for any $t>0$ and $m\in M$, the inequality $ f(m) -
    f( \phi_t (m)) < t$ holds. Conversely for any continuous function $f$
    satisfying this
    inequality, the map $m \mapsto \phi_{f(m)} (m)$ is a
    homeomorphism of $M$ respecting both foliations leaf by
    leaf. So, using the family $( \lambda f)_{ \lambda \in [0,1]}$, $h$ is
    isotopic to the identity through homeomorphisms respecting the foliations. 
    We observe that the homeomorphisms homotopic to the identity 
    and respecting the foliation $\mathcal{G}$ are 
    precisely the ones sending each geodesic leaf to itself.
  \end{itemize}
\end{remark}

A developing pair $(\dev, \hol)$ of a projective structure on $M$
defines a foliated projective
structure on  $(M,\mathcal{F}, \mathcal{G})$ if the following
conditions are satisfied. 
\begin{enumerate}
\item For every $g\in \widetilde{\mathcal{G}}$, the image $\dev(g)$ is 
  contained in a projective line, and 
\item for every $f \in \widetilde{\mathcal{F}}$, the image $\dev(f)$
  is contained in a projective plane.
\end{enumerate}

\subsubsection{Properly Convex Foliated Projective Structures}
\label{sec_defPCFP}

\begin{defi}
  A foliated projective structure on $(M,\mathcal{F}, \mathcal{G})$ is
  said to be \emph{convex} if for every $f\in \widetilde{\mathcal{F}}$
  the image $\dev(f)$ is a convex set in the projective plane
  containing $\dev(f)$. It is \emph{properly convex} if for every
  $f\in \widetilde{\mathcal{F}}$ the image $\dev(f)$ is a properly
  convex set.
\end{defi}

Note that we do not require that the restriction of $\dev$ to a (weakly) stable leaf is a 
homeomorphism onto its image but this will be a consequence of the other
conditions. In fact, already for convex projective structures on closed
surfaces one can show that such a condition on the image of the developing map
suffices. 

Let $\mathcal{P}_{cf} (M)$ denote the set of equivalence classes of
convex foliated projective structures on $M$ and $\mathcal{P}_{pcf}(M)$
the set of equivalences classes of properly convex
foliated projective structures. They are naturally subsets
of the moduli space $\mathcal{P}_{f}(M)$ of foliated projective
structures on $M$.

\subsection{The Hitchin Component}
Let $\Gamma$ be as above the fundamental group of a closed Riemann surface
$\Sigma$, and denote by $\rep(\G, \PGL_n(\RR))$ the set of conjugacy classes
of representations, that is
\begin{equation*}
  \rep(\G, \PGL_n(\RR)) = \hom(\Gamma, \PGL_n(\RR))/\PGL_n(\RR).
\end{equation*}
Of course this space (with the quotient topology) is not Hausdorff, but we will not worry about this 
since we are just concerned with a component of it that has
the topology of a ball.

Relying on the correspondence between 
stable Higgs bundle and irreducible representations of $\pi_1(\Sigma
)$ in $\PSL_n(\CC)$ due to N.~Hitchin \cite{Hitchin_1987}, C.~Simpson
\cite{Simpson_1988, Simpson_1992}, K.~Corlette \cite{Corlette_1988} and
S.~Donaldson \cite{Donaldson_1985},
N.~Hitchin proved in \cite{Hitchin_1992} that the
connected component of $\hom(\G, \PSL_n(\RR))/ \PGL_n(\RR)$ containing the
representation 
\begin{equation*}
  \rho_n\circ \iota: \Gamma \to \PSL_n(\RR) \subset \PGL_n( \RR), 
\end{equation*}
is a ball of dimension $(2g-2) (n^2 -1)$. Here
$\iota: \Gamma \to \PSL_2(\RR)$ is a uniformization and
$\rho_n$ is the n-dimensional
irreducible representation of $\PSL_2(\RR)$.
\begin{nota}
  This component is called the \emph{Teichm\"uller component} or 
\emph{Hitchin component} and is denoted by $\mathcal{T}^n( \Gamma)$ or 
$\mathcal{T}^n( \Sigma)$.
\end{nota}
The Hitchin component naturally embeds
\begin{equation*}
  \mathcal{T}^n( \Gamma) \subset \rep(\G, \PGL_n(\RR)) \subset
  \rep(\overline{\G}, \PGL_n(\RR))
\end{equation*}
where the last inclusion comes from the projection $\overline{ \Gamma}
\rightarrow \Gamma$. Even though we do not use it, the
Hitchin component can be seen as a connected component of $\rep(\overline{\G}, \PGL_n(\RR))$:
\begin{lem}
  The Hitchin component $\mathcal{T}^n ( \Gamma)$ is a
  connected component of $\rep(\overline{\G}, \PGL_n(\RR))$.
\end{lem}

\begin{proof}
  Clearly $\mathcal{T}^n ( \Gamma)$ is closed in $\rep(\overline{\G},
  \PGL_n(\RR))$. By \cite{Hitchin_1992} any representation $\rho \in
  \mathcal{T}^n ( \Gamma)$ is strongly irreducible (e.g.
  \cite[Lemma~10.1]{Labourie_2003} explains this fact). Therefore
  $\mathcal{T}^n ( \Gamma)$ is contained in the open subset of irreducible
  representations $\rep_{irr}(\overline{\G}, \PGL_n(\RR))$. Since any
  irreducible representation $\overline{\rho}: \overline{\G} \to \PGL_n(\RR)$
  necessarily factors through a representation $\rho: \Gamma \to \PGL_n(\RR)$,
  this shows that $\mathcal{T}^n ( \Gamma)$ is open in
  $\rep_{irr}(\overline{\G}, \PGL_n(\RR))$ and hence open in
  $\rep(\overline{\G}, \PGL_n(\RR))$.
\end{proof}

With this we are now able to restate our main result:
\begin{thm}\label{thm:main}
  The holonomy map
\begin{equation*}
  \hol: \mathcal{P}_{pcf} (M) \longrightarrow \rep(\overline{\G}, \PGL_4(\RR))
\end{equation*}
is a homeomorphism onto the Hitchin component $\mathcal{T}^4 ( \Gamma)$.
\end{thm} 

\section{Examples}

\label{sec_examples}
\begin{quote}
  In this section we define several families of projective structures on $M$.
  These families will provide some justification for the definitions given in
  the previous section and make the reader acquainted with some geometric
  constructions which will be used in the following sections.
\end{quote}

We can summarize this section in the following
\begin{prop}
  \label{prop_incstrict}
  All the inclusions
  \begin{equation*}
     \mathcal{P}_{pcf}(M) \subset \mathcal{P}_{cf} (M) \subset
  \mathcal{P}_{f}(M)
  \end{equation*}
  are strict and the projection from $\mathcal{P}_{f} (M)$ to
  $\mathcal{P}(M)$, defined by forgetting the foliations, is not onto.
\end{prop}

In order to define a $( \textup{PGL}_4( \RR), \PTR)$-structure, it is
sufficient to give a pair $(\dev, \hol)$ consisting of the holonomy
representation $\hol:\overline{\G} \to \PGL_4(\RR)$ and the $\hol$-equivariant
developing map $\dev: \widetilde{M} \to \PTR$. In all examples given below the
holonomy factors through $\Gamma = \overline{ \Gamma}/ \langle \tau \rangle$
and the developing map factors through the quotient $\overline{M} =
\widetilde{ M}/ \langle\tau\rangle = S \widetilde{ \Sigma}$.  In particular we
will specify developing pairs $( \dev, \hol)$ with $\hol:\G \to \PSL_4(\RR)$
and $\dev: \overline{M} \to \PTR$.

\subsection{Homogeneous Examples}
\label{sec_homoexamples}
We first construct families of projective structures on $M$ which are induced
from homogeneous projective structures on $\PSL_2( \RR)$ when we realize $M$ as
a quotient of $\PSL_2( \RR)$ (see Paragraph~\ref{sec_Masquotient}). Since this
procedure will be used several times we state it in the following lemma.

\begin{lem}\label{lem_examples}
  Let $\iota : \Gamma \rightarrow \PSL_2( \RR)$ be a uniformization, $\rho:
  \PSL_2( \RR) \rightarrow \PSL_4( \RR)$ a homomorphism and $x$ a point in
  $\PTR$ with
  \begin{equation*}
    \textup{Stab}_{ \PSL_2 ( \RR)}( x) = \big \{ g \in \PSL_2( \RR) \mid \rho(
    g) \cdot x = x \big \} \text{ being finite}.
  \end{equation*}
  Then the following assignment
  \begin{eqnarray*}
    \dev: \overline{M} & = &\PSL_2 (\RR) \rightarrow \PTR \, ,  \, g
    \mapsto \rho( g) \cdot x \\
    \hol & = &\rho \circ \iota
  \end{eqnarray*}
  defines a projective structure on $M$. Furthermore
  \begin{itemize}
  \item The images of geodesics are contained in projective lines if and only
    if $\dev( A)$ is.
  \item The images of (weakly) stable leaves are contained in projective planes if and
    only if $\dev( P)$ is.
  \item The restriction of $\dev$ to every (weakly) stable leaf is a homeomorphism onto
    a (properly) convex set if and only if the restriction to $P$ is so.
  \end{itemize}
\end{lem}

\begin{proof}
  The fact that $\dev$ is a $\hol$-equivariant local homeomorphism is
  clear. The listed properties follow from the $\PSL_2(
  \RR)$-equivariance 
  and the description of the leaves (Paragraph~\ref{sec_Masquotient}).
\end{proof}

\begin{remark}
  This lemma could also be stated by saying that given a locally
  homogeneous $( \PSL_2( \RR) ,
  \PSL_2( \RR))$-struc\-ture on $M$ and a $( \PSL_2(
  \RR), U)$ geometry with $U$ an open in $\PTR$ we automatically obtain a $( \PSL_2(
  \RR), U)$-structure and hence a projective structure on $M$.
\end{remark}

\subsubsection{The Diagonal Embedding}
\label{sec_diagemb}
We keep the notation of Lemma~\ref{lem_examples}. If we choose
\begin{center}
  $\displaystyle 
  \begin{array}{rcl}
    \rho: \PSL_2( \RR) & \longrightarrow & \PSL_4( \RR) \\
    g & \longmapsto & \Big( 
    \begin{array}{cc}
      g & 0 \\ 0 & g
    \end{array} \Big)
  \end{array}$ and $\displaystyle x = [ 1, 0 , 0, 1 ]$,
\end{center}
we obtain a convex foliated projective structure on $M$ with non properly
convex (weakly) stable leaves. Indeed $\textup{Stab}_{ \PSL_2 ( \RR)}( x) = \{
\textup{Id} \}$ shows that $\dev$ is a homeomorphism onto its image
and:
\begin{eqnarray*}
  \dev( A) & = & \{ [ ( e^t, 0,0,1) ] \mid t \in \RR \}, \\
  \dev( P) & = & \{ [ ( 1, 0,u,v )] \mid u \in \RR, v > 0 \}.
\end{eqnarray*}
So $\dev(A)$ is contained in a projective line, $\dev( P)$ is convex in a
projective plane but its closure contains a projective line.

In Section~\ref{sec_quasifuchs} and Section~\ref{sec_geomdesc} we
generalize this example in two different ways. 

\subsubsection{The Irreducible Example}
\label{sec_exampleirr}
Let $\rho_4: \PSL_2(
\RR) \rightarrow \PSL_4( \RR)$ be the $4$-dimensional irreducible 
representation. In other words this is the representation on the $3$-fold
symmetric power $\RR^4 \simeq \textup{Sym}^3 \RR^2$. Hence we will
consider elements of $\RR^4$ as homogeneous polynomials of degree three
in two 
variables $X$ and $Y$, so it will make sense to speak of the roots of
an element of $\RR^4$. 
We choose $x = [R] \in \PTR$ such that $R$ has 
only one real root counted with multiplicity.

A direct calculation shows that for the projective structure defined
by Lemma~\ref{lem_examples} $\dev(A)$ and $\dev(P)$ are contained
respectively in a projective line and in a projective plane. Actually,
if $R = X( X^2 + Y^2) = ( 1,0,1,0)$,
\begin{equation*}
  \dev(P) = \big\{ [ a^4+ a^2b^2, 2ab, 1,0 ] \in \PTR \mid a,b
    \in \RR \big\}
\end{equation*}
is the projection of the properly convex cone of $\RR^4 - \{0\}$:
\begin{equation*}
  \big\{ ( \alpha, \beta, \gamma, 0) \in \RR^4 \mid \beta^2 - 4 \alpha
  \gamma < 0 \big \}.
\end{equation*}

So this defines a properly convex foliated
projective structure on $M$. In Section~\ref{sec_curvehomo} we will
give a more geometric description of this example.

\begin{remark}
  Note that the choice of $R$ is related to the previous choice of a Cartan subgroup $A$. 
Any other polynomial $g\cdot R$ in the orbit of $R$ defines the same structure when parametrizing the geodesic flow 
by $g A g^{-1}$ instead of $A$. 
\end{remark}

\subsubsection{The Irreducible Example Revisited}
\label{sec_irrevis}
Taking the same representation $\rho_4: \PSL_2(\RR) \to \PSL_4(\RR)$
but choosing $x=[Q]$ where $Q$ is a polynomial having three distinct real
roots, the projective structure on $M$ defined by
Lemma~\ref{lem_examples} is foliated but not
convex. In this case, the image of $P$ is the complementary
subset (in the projective plane containing it) of the closure of convex set
described above and cannot be convex.

\subsection{Nonfoliated Structures with Quasi-Fuchsian
  Holonomy}
\label{sec_quasifuchs}

In this paragraph we generalize the example of
Section~\ref{sec_diagemb} for any
quasi-Fuchsian representation $q : \Gamma \rightarrow \PSL_2(
\CC)$. We thank Bill Goldman for explaining us this construction.

Recall that any quasi-Fuchsian representation $q$ is a deformation
of a Fuchsian representation $\Gamma \rightarrow \PSL_2( \RR)$, and
there exists a $q$-equivariant local orientation preserving
homeomorphism $u: \widetilde{ \Sigma} \rightarrow \POC$.

Fixing an identification $\CC \simeq \RR^2$ we have an embedding
$\PSL_2( \CC) \subset \PSL_4( \RR)$ such that the Hopf fibration $\PTR
\rightarrow \POC$ is a $\PSL_2( \CC)$-equivariant fibration by
circles.

\begin{prop}
  \label{prop_qfprojstru}
  Let $q: \G \to \PSL_2(\CC)$ be a quasi-Fuchsian representation
  and $u : \widetilde{ \Sigma} \rightarrow \POC$ a $q$-equivariant
  local orientation preserving homeomorphism.
  \begin{enumerate}
  \item \label{item1lemqf} Then the pull back $u^* \PTR$ of the Hopf
    fibration $\PTR \rightarrow \POC$ admits a $\Gamma$-invariant
    projective structure and the quotient of $u^* \PTR$ by $\Gamma$ is
    homeomorphic to $M$.
  \item \label{item2lemqf} The induced projective structure on $M$ is
    foliated if and only if the representation $q$ is (conjugate
    to) a Fuchsian representation $\Gamma \rightarrow \PSL_2( \RR)$.
  \end{enumerate}
\end{prop}

\begin{proof}
  For (\ref{item1lemqf}), the projective structure on the pull back
  $u^* \PTR$ is tautological since $u^* \PTR \rightarrow \PTR$ is a
  local homeomorphism. From the homotopy invariance of fiber bundles
  (see \cite[p.53]{Steenrod_1951}) it is enough to show that the
  quotient of $u^* \PTR$ by $\Gamma$ is homeomorphic to $M$ when the
  representation $q: \Gamma \rightarrow \PSL_2(\RR)$ is Fuchsian.
  In this case $u$ is the composition $\widetilde{\Sigma} \simeq \HH^2
  \hookrightarrow \POC$. Also the embedding $\PSL_2( \RR) \subset
  \PSL_2( \CC) \subset \PSL_4( \RR)$ is the diagonal embedding of
  Section~\ref{sec_diagemb} and the map $\dev: \overline{M} \simeq
  \PSL_2(\RR) \rightarrow \PTR$ defined in Section~\ref{sec_diagemb}
  fits into the commutative diagram:
  \begin{equation*}
    \begin{array}{rcccc}
      \dev: \overline{M} & \stackrel{\sim}{ \longrightarrow} & \PSL_2( \RR) &
      \stackrel{\sim}{ \longrightarrow} & \PTR \\
      & & \big \downarrow & & \big \downarrow \\
      u: \widetilde{\Sigma} & \stackrel{\sim}{ \longrightarrow} & \HH^2 &
      \stackrel{\sim}{ \longrightarrow}  & \POC.
    \end{array}
  \end{equation*}
  So $\overline{M}$ is homeomorphic to $u^* \PTR$ as $\G$-space.

  (\ref{item2lemqf}) If the structure is foliated there exists a
  $q$-equivariant map $\xi^3: \partial \Gamma \rightarrow \PTRd$
  (see Proposition~\ref{prop:fol_maps}). In particular any element
  $\gamma$ in $\Gamma -\{1\}$ will have an eigenline $\xi^3( t)$ in
  $\RR^{4*}$, with $t$ a fixed point of $\gamma$ in $\partial \Gamma$,
  and hence a real eigenvalue for its action on $\RR^4$. Since the
  eigenvalues of an element of $\PSL_2( \CC) \subset \PSL_4( \RR)$ are
  $\{ \lambda, \lambda, \lambda^{-1}, \lambda^{-1} \}$ we deduce that
  $q( \gamma)$ has only real eigenvalues. In particular
  $\textup{tr}\left( q( \gamma)\right) = \textup{tr}
  \left(\overline{ q(\gamma)}\right)$. So the representation
  $(q, \overline{q}): \Gamma \to \PSL_2(\CC) \times
  \PSL_2(\CC)$ cannot have Zariski-dense image. Therefore, by
  Lemma~\ref{lem_adhzarofpair}, there exists $g$ in $GL_2( \CC)$ such
  that, $q( \gamma) = g \overline{ q(\gamma)} g^{-1}$ for all
  $\gamma \in \Gamma$. Since $q(\Gamma)$ is Zariski dense the
  element $g \overline{g}$ is central. Up to multiplying $g$ by a
  scalar we have $g\overline{g} = \pm \textup{Id}$.

  If we have $g \overline{g} = \textup{Id}$ then, for some $\beta$ in
  $\CC$, $h= \beta g + \overline{\beta} \textup{Id}$ belongs to
  $\GL_2( \CC)$ and satisfies $g \overline{h} =h$. The discrete and
  faithful representation $\gamma \rightarrow h q(\gamma) h^{-1}$
  takes values in $\PSL_2( \RR)$ so $q$ is conjugate to a Fuchsian
  representation.

  If $g \overline{g} = - \textup{Id}$, setting $T= \Big(
  {\begin{array}{cc} \scriptstyle 0 & \scriptstyle 1\\ \scriptstyle -1 &
  \scriptstyle 0
    \end{array}}\Big)$ then, for $\beta \in \CC$, $h= \beta g +
  \overline{\beta}T$ is invertible. The representation $h q
  h^{-1}$ has values in $\textup{PSU}_2(\RR)$ but this is impossible since
  $q$ is faithful and discrete.
\end{proof}

\subsection{Geometric Description of the Diagonal Embedding}
\label{sec_geomdesc}

Recall that in this example the holonomy is
\begin{eqnarray*}
  \hol: \Gamma & \longrightarrow & \PSL_4( \RR)\\
  \gamma & \longmapsto & \rho( \iota( \gamma)) = \left(
    \begin{array}{cc}
      \iota(\gamma) & 0 \\ 0 & \iota(\gamma)
    \end{array}\right)
\end{eqnarray*}
and the developing map is
\begin{eqnarray*}
  \dev: \PSL_2( \RR) \simeq \overline{M} & \longrightarrow & \PTR \\
  g & \longmapsto & \rho(g) \cdot x
\end{eqnarray*}
where $x = [(1,0,0,1)]$. We wish to describe $\dev$ as a map $\partial
\Gamma^{3+} \simeq \overline{M} \rightarrow \PTR$.

It will be useful to have a lift of $\hol$ to $\SL_4( \RR)$
\begin{eqnarray*}
  \widehat{ \hol} : \Gamma & \longrightarrow & \SL_4( \RR) \\
  \gamma & \longmapsto & \left(
    \begin{array}{cc}
      \hat{\iota}( \gamma) & 0 \\ 0 & \hat{\iota}( \gamma)
    \end{array}\right)
\end{eqnarray*}
where $\hat{ \iota}: \Gamma \rightarrow \SL_2(\RR)$ is one of the
$2^g$ lifts of $\iota$. We also remind the reader that $\partial
\Gamma$ is being equivariantly 
identified with $\POR$. The sphere $\mathbb{S}( \RR^2) =
\RR^2-\{0\}/\{ x \sim \lambda^2 x\}$  is a $\Gamma$-space that
projects onto $\partial \Gamma= \POR = \mathbb{S}( \RR^2)/\{\pm1\}$, we
denote it by $\widehat{\partial \Gamma}$. It has the following property:

for any $(t_+, t_-)$ in $\partial \Gamma^{(2)}$, there are exactly
two lifts $(\hat{t}_+, \hat{t}_-)$ in $\widehat{ \partial \Gamma}^2$
of this pair so that $(-\hat{t}_-,\hat{t}_+,\hat{t}_-)$ is oriented,
these two lifts are exchanged by the action of $-1$.

A straightforward calculation shows that the image $\dev( t_+, t_-)$ of any
geodesic is a segment in $\PTR$ with endpoints $\xi^+( t_+)$ and $\xi^-(
t_-)$ where $\xi^\pm: \partial \Gamma \rightarrow \PTR$ are the two
$\rho$-equivariant maps $\POR \hookrightarrow \PTR$ coming from the
decomposition $\RR^4 = \RR^2 \oplus \RR^2$. To describe $\dev$ we need
lifts of $\xi^\pm$ to $\RR^4$. No continuous lift exists, so we rather
choose $\eta^{\pm}: \widehat{ \partial \Gamma} \rightarrow \RR^4$ two
continuous maps, equivariant by $-1$, lifting $\widehat{ \xi}^\pm :
\widehat{ \partial \Gamma} \rightarrow \partial \Gamma \rightarrow
\PTR$. Then there exists a continuous function $\vfi : \partial \Gamma^{3+}
\rightarrow \RR$ such that, for all $( t_+, t_0, t_-) \in \partial \Gamma^{3+}$
\begin{equation}
  \label{eq_phidev}
  \dev( t_+, t_0, t_-) = [ \eta^+( \hat{t}_+) + \vfi( t_+, t_0, t_-)
  \eta^-( \hat{t}_-) ] \in \PTR,
\end{equation}
where $( \hat{t}_+, \hat{t}_-)$ is one of the two lifts of $( t_+,
t_-)$ with $( - \hat{t}_-, \hat{t}_+, \hat{t}_-)$ oriented. Observing
that $\vfi$ never vanishes, we can suppose $\vfi>0$ up to changing
$\eta^+$ in $-\eta^+$. Since $\dev$ is a local homeomorphism, $\vfi$
must be monotonely decreasing along geodesics. 

To state the condition on $\vfi$ coming from the equivariance of
$\dev$ we consider the two maps $f^\pm :
\Gamma \times \partial \Gamma \rightarrow \RR^*$ measuring the lack of
equivariance of $\eta^\pm$: for all $\hat{t} \in \widehat{ \partial
  \Gamma}$ projecting on $t \in \partial \Gamma$
\begin{equation}
  \label{eq_equivf}
  \eta^\pm( \gamma \cdot \hat{t}) = f^\pm( \gamma, t) \widehat{\hol}(
  \gamma) \cdot \eta^\pm( \hat{t}).
\end{equation}
With this, for all $( t_+, t_0, t_-) \in \partial \Gamma^{3+}$ and
$\gamma \in \Gamma$, the developing map $\dev$ is equivariant if and
only if 
\begin{equation}
  \label{eq_equivphi}
  \vfi( \gamma \cdot ( t_+, t_0, t_-)) = \frac{ f^-( \gamma, t_-) }{
    f^+( \gamma, t_+) } \vfi( t_+, t_0, t_-).
\end{equation}

Given maps with those conditions, we can construct a foliated projective
structure on $M$. 

\begin{prop}
  \label{prop_xiandfidefproj}
  Let $\hol: \Gamma \rightarrow \PSL_4( \RR)$ be a representation and
  $\widehat{ \hol}: \Gamma \rightarrow \SL_4( \RR)$ a lift of
  $\hol$. Suppose that
  \begin{enumerate}
  \item \label{itpropi} there exist two continuous $\hol$-equivariant maps
    $\xi^\pm: \partial \Gamma \rightarrow \PTR$, and
  \item the image of $\xi^-$ is contained in a projective line,
  \item there exist two lifts $\eta^\pm: \widehat{ \partial \Gamma}
    \rightarrow \RR^4$ of $\widehat{ \xi}^\pm : \widehat{ \partial \Gamma}
    \rightarrow \partial \Gamma \rightarrow \PTR$ and functions $f^\pm:\Gamma
    \times \partial \Gamma \rightarrow \RR^*$ satisfying \eqref{eq_equivf},
  \item \label{itpropiv} there exists a continuous function $\vfi: \partial
    \Gamma^{3+} \rightarrow \RR_{>0}$ satisfying the identity
    \eqref{eq_equivphi}
  \item \label{itpropv} and $\vfi$ satisfies the limit condition: for all $(
    t_+, t_-)$ in $\partial \Gamma^{(2)}$
  \begin{equation}
    \label{eq_limphi}
    \lim_{
    t_0 \rightarrow t_+} \vfi( t_+, t_0, t_-) = 0\text{ and } \lim_{
    t_0 \rightarrow t_-} \vfi( t_+, t_0, t_-) = \infty.
  \end{equation}
  \end{enumerate}
  Moreover, suppose that the map $\dev: \partial \Gamma^{3+} \rightarrow \PTR$
  defined by \eqref{eq_phidev} is a local homeomorphism.
  
  Then the pair $(\dev, \hol)$ defines on $M$ a convex foliated projective
  structure which is not properly convex.
\end{prop}

\begin{proof}
  By construction $\dev$ is a $\hol$-equivariant local homeomorphism and
  images of geodesics are contained in projective lines. The limit condition
  on $\vfi$ and the fact that $\xi^-$ is contained in a projective line $L$
  imply that the image of the (weakly) stable leaf $t_+$ is a sector whose tip is
  $\xi^-( t_+)$ and which is bounded by the projective lines $L$ and
  $\overline{ \xi^+( t_+) \xi^-( t_+)}$.
\end{proof}

A few remarks need to be made about this construction.

Instead of imposing that $\dev$ is a local homeomorphism, we could have stated
conditions on $\xi^\pm$ and $\vfi$ which imply it. For example a condition of
the type:
\begin{itemize}
\item The function $\vfi$ is monotonely decreasing along geodesics, $\xi^\pm$
  are homeomorphisms onto $\mathcal{C}^1$-submanifolds of $\PTR$ satisfying
  the condition that if $L^\pm$ are projective lines tangent to $\xi^\pm(
  \partial \Gamma)$ at $\xi^\pm( t_\pm)$ with $t_+ \neq t_-$, then $L^+$ and
  $L^-$ do not intersects.
\end{itemize}
would suffice. It is easy to see that $\xi^\pm$ cannot be locally constant.

In a decomposition $\RR^4 = \RR^2 \oplus \RR^2$ adapted to the image of
$\xi^-$, $\hol$ (and $\widehat{ \hol}$) has the form
\begin{equation*}
  \gamma \longrightarrow \left(
    \begin{array}{cc}
      \rho_Q( \gamma) & 0 \\ c(\gamma) & \rho_L( \gamma)
    \end{array}\right).
\end{equation*}
It can be shown that $\rho_Q$ and $\rho_L$ are Fuchsian (see
Lemma~\ref{lem_curveteich}) and that $\xi^\pm$ are uniquely determined by
$\hol$. In fact $\xi^-$ is a homeomorphism onto its image; this was implicitly
used in the above proof.

Our analysis shows that the holonomy of non-properly convex foliated
projective structures have to be of the above form with $\rho_Q$ and $\rho_L$
Fuchsian. Unfortunately we do not have a construction of a nontrivial example
besides the diagonal example.

Condition \eqref{eq_limphi} does indeed depend only on $\hol$ and not on
$\vfi$.  If $\vfi_1$, $\vfi_2$ satisfy this condition then their ratio $\vfi_1
/ \vfi_2$ descends to a continuous function on $M$ and hence is bounded. If we
change the lift $\eta^+$ then the function $\vfi$ will change to $( t_+, t_0,
t_-) \mapsto \lambda(t_+) \vfi( t_+, t_0, t_-)$ for some continuous function
$\lambda : \partial \Gamma \rightarrow \RR_{>0}$. So the behavior at infinity
of the continuous functions satisfying equality \eqref{eq_equivphi} depends
only on $\hol$ and $\xi^\pm$, and the curves $\xi^\pm$ are uniquely determined
by $\hol$.

The above description of the diagonal embedding and the following lemma will be used in the proof of
Lemma~\ref{lem_interxi3notline}.

\begin{lem}
  \label{lem_notpropconv}
  Suppose that $\hol$ is the holonomy representation of a projective structure
  on $M$ such that Conditions (\ref{itpropi})-(\ref{itpropiv}) of
  Proposition~\ref{prop_xiandfidefproj} are satisfied. Assume that its
  semisimplification $\hol_0$ satisfies Conditions
  (\ref{itpropi})-(\ref{itpropv}). Then all (weakly) stable leaves of the projective
  structure associated to $\hol$ are developed into sectors.
\end{lem}

\begin{proof}
  First note that the images of all (weakly) stable leaves are sectors if and only if
  Condition \eqref{eq_limphi} is satisfied. We write $\hol$ as
  above with respect to a decomposition $\RR^4 = \RR^2 \oplus \RR^2$ adapted to
  $\xi^-$. Let $\eta^\pm: \widehat{ \partial \Gamma} \rightarrow
  \RR^4$ be the lifted curves for $\hol$. We decompose
  \begin{equation*}
    \eta^+( t) = \eta^{+}_{Q}(t) + \eta^{+}_{L}(t) \in \RR^2 \oplus \RR^2.
  \end{equation*}
  Then the curves for
  \begin{eqnarray*}
    \hol_0 : \Gamma & \longrightarrow & \PSL_4( \RR) \\
    \gamma & \longmapsto &  \left( \begin{array}{cc}
      \rho_Q( \gamma) & 0 \\ 0 & \rho_L( \gamma)
    \end{array}\right)
  \end{eqnarray*}
  are $\eta^{+}_{Q}$ and $\eta^-_{L}$. The functions $f^\pm$ (Equation
  \eqref{eq_equivf}) are the same for $( \eta^+, \eta^-)$
  and $(\eta^{+}_{Q}, \eta^-_{L})$. Therefore the equivariance condition for
  $\vfi$ is the same for both representations $\hol$ and $\hol_0$. So 
  by the above remark their behavior are infinity is the same. In
  particular Condition \eqref{eq_limphi} holds for $\hol$.
\end{proof}

\section{From Convex Representations to Properly Convex Foliated
  Projective Structures}
\label{sec_reptostruct}

\begin{quote}
  In this section we first describe the properly convex foliated
  projective structure on $M= S\Sigma$ defined in Section~\ref{sec_exampleirr}
  geometrically. 
  This geometric description of the developing
  map enables us to construct properly convex foliated projective
  structure with $\hol = \rho$ for any representation $\rho$ in the
  Hitchin component. 
\end{quote}

\subsection{A Different Description of the Homogeneous Example}
\label{sec_curvehomo}
Our aim is now to describe the developing map of
Example~\ref{sec_exampleirr} as a
map $\partial \Gamma^{3+} \rightarrow \PTR$. Recall that the holonomy was $\hol
= \rho_4 \circ \iota$ with $\iota$ a Fuchsian representation and $\rho_4$ the
$4$-dimensional irreducible representation of $\PSL_2( \RR)$, that is the
representation on $\RR^4= \textup{Sym}^3 \RR^2$ the space of homogeneous
polynomials of degree three in two variables $X$ and $Y$. The developing map
was
\begin{eqnarray}
  \label{eq_dev}
  \notag \dev : \overline{M} \simeq \PSL_2( \RR) & \longrightarrow & \PTR \\
  g & \longmapsto & \rho_4(g) \cdot [R],
\end{eqnarray}
where $R$ is  $X( X^2 +
Y^2)$. In fact the description for the nonconvex example given in Section~\ref{sec_irrevis},
where the holonomy is the same but the developing map is
\begin{eqnarray*}
  \dev' : \overline{M} \simeq \PSL_2( \RR) & \longrightarrow & \PTR \\
  g & \longmapsto & \rho_4(g) \cdot [Q],
\end{eqnarray*}
with $Q=XY(X+Y)$, is easier to obtain.

It is convenient to consider $\RR^2 = \textup{Sym}^1\RR^2$ as the
space of homogeneous polynomials of degree one in $X$ and $Y$.  
The Veronese embedding 
\begin{eqnarray*}
  \xi^1: \partial \Gamma \simeq \POR &\longrightarrow& \PTR\\ 
  t = {[S]} & \longmapsto & [S^3]
\end{eqnarray*}
is a $\rho_4$-equivariant map, so certainly $\hol$-equivariant, and extends to
an equivariant map into the flag variety 
\begin{equation*}
  \xi = (\xi^1, \xi^2,\xi^3):  \partial \Gamma \longrightarrow \Flag(\RR^4), 
\end{equation*}
which is also $\rho$-equivariant.

The maps $\xi^i$ can be described as follows:
\begin{itemize}
\item $\xi^1([S])$ is the line of polynomials divisible by $S^3$,
\item $\xi^2([S])$ is the plane of polynomials divisible by $S^2$,
\item $\xi^3([S])$ is the hyperplane of polynomials divisible by $S$.
\end{itemize}
The four orbits of $\PSL_2(\RR)$ in $\PTR$ can be described in term of $\xi$: 
\begin{itemize}
\item one open orbit $\Lambda_{\xi}$ which is the set of polynomials having
  three distinct real roots, \emph{i.e.} the points in $\PTR$ that are
  contained in (exactly) three pairwise distinct
  $\xi^3(t)$,
\item the other open orbit $\Omega_{\xi}$ which is the set of
  polynomials having a pair of complex conjugate roots, \emph{i.e} points of $\PTR$
  contained in exactly one $\xi^3(t)$,
\item the relatively closed orbit is the surface $\bigcup_{t \in \partial
    \Gamma} \xi^2(t) \backslash \xi^1(\partial \Gamma)$,
\item and the closed orbit is the curve $\xi^1(\partial \Gamma)$.
\end{itemize}

Remark that the two open orbits are the connected components of the
complementary of the surface, called \emph{discriminant}, $\xi^2( \partial
\Gamma) = \bigcup_{t \in \partial \Gamma} \xi^2(t) \subset \PTR$.

\subsubsection{The Nonconvex Example}
Let us first describe the developing map of the nonconvex foliated
projective 
structure $(\dev', \hol)$. The open orbit $\Lambda_{\xi}$ coincides
with the image of $\dev'$. The developing map $\dev'$ can be described as
\begin{eqnarray}
  \label{eq_dev'curve}
  \notag \dev': \partial\Gamma^{3+} \simeq \overline{M} &\longrightarrow & \PTR
  \\
  (t_+, t_0, t_-) &\longmapsto& \xi^3(t_+) \cap \xi^3(t_0) \cap \xi^3(t_-).
\end{eqnarray}
The image of the geodesic $(t_+, t_-) \in \partial \Gamma^{(2)}$ is an open
segment in the projective line $\xi^3(t_+) \cap \xi^3(t_-)$. The image of the 
(weakly) stable leaf $t_+$ is contained in the projective plane
$\xi^3(t_+)$.

Note that the only property of the $\hol$-equivariant map $\xi:
\partial \Gamma \to \Flag(\RR^4)$ needed to define the developing map is
that, for any pairwise distinct $t_1, t_2, t_3\in \partial \Gamma$, the three
projective planes $\xi^3(t_1)$, $\xi^3(t_2)$ and $\xi^3(t_3)$
intersect in a unique point in $\PTR$. 
To ensure that the developing map is a local homeomorphism we need
that for $t_1, t_2, t_3, t_4$ pairwise distinct the intersection
$\bigcap_{i=1}^4 \xi^3(t_i) = \emptyset$.

\subsubsection{The Properly Convex Example}
Let us now describe the developing map of the properly convex foliated
structure $(\dev, \hol)$ defined above in terms of the $\hol$-equivariant map
$\xi: \partial \Gamma \to \Flag(\RR^4)$. First note that the image of $\dev$
is the open orbit $\Omega_\xi$.

The map $\xi:\partial \Gamma \to \Flag(\RR^4)$ has the property that, for every
distinct $t,t' \in \partial \Gamma$, the intersection $\xi^3(t) \cap
\xi^2(t')$ is a point. 
We can therefore define an equivariant map, or rather an equivariant family of
maps $(\xi_t)_{t\in \partial\Gamma}$ 
\begin{eqnarray*}
 \partial \Gamma \times \partial \Gamma &\longrightarrow& \PTR\\
  (t,t') &\longmapsto& \xi_t(t') = \Big\{ 
    \begin{array}{cl}
      \xi^3(t) \cap \xi^2(t') & \text{if } t \neq t'\\
      \xi^1(t) & \text{otherwise}
    \end{array}. 
\end{eqnarray*}
Note that for every $t$ the image of $\xi^1_t$ in $\xi^3(t)$ bounds the
properly convex domain $C_t= \dev(t)$.

Having this family of maps, we see that the image of the geodesic leaf $g =
(t_+, t_-) \in \partial \Gamma^{(2)}$ under the developing map
$\dev$ of \eqref{eq_dev} is the intersection of the projective line
$\overline{ \xi^1(t_+) \xi^1_{t_+}(t_-)} = \PP(\xi^1(t_+) \oplus
\xi^1_{t_+}(t_-))$ with the convex $C_{t_+}$. The endpoint
at $+\infty$ of the open segment $\dev(g)$ is $\xi^1(t_+)$ and
$\xi^1_{t_+}(t_-)$ is the endpoint at $-\infty$. 
\begin{figure}[htbp]
  \centering
\input{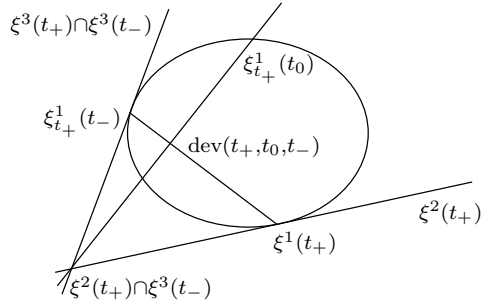}
  \caption{The image of the developing map}
  \label{fig_devmap}
\end{figure}
Moreover the projective line $\xi^2(t_+)$ is tangent to the convex
$C_{t_+}$ at $\xi^1(t_+)$. The tangent line at the point $\xi^1_{t_+}(t_-)$ is
the projective line $\xi^3(t_+) \cap \xi^3(t_-)$. These two projective lines
intersect in the point $\xi^2(t_+) \cap \xi^3(t_-) = \xi^1_{t_-}(t_+)$.
 
Note that given $t_0 \in \partial\Gamma$ distinct from $t_+$ and $t_-$
the two projective lines $\overline{\xi^1(t_+) \xi^1_{t_+}(t_-)}$ and 
$\overline{\xi^1_{t_-} (t_+) \xi^1_{t_+}(t_0)}$ intersect in a unique
point that belongs to $C_{t_+}$.
With this we can now give an explicit formula for the developing map
\begin{eqnarray}
  \label{eq_devcurv}
  \notag \dev: \partial\Gamma^{3+} \simeq \overline{M}
  &\longrightarrow& \PTR \\
  (t_+, t_0, t_-) &\longmapsto& \overline{\xi^1(t_+)
  \xi^{1}_{t_+}(t_-)} \cap 
\overline{\xi^{1}_{t_-} (t_+)  \xi^{1}_{t_+}(t_0)}.
\end{eqnarray}

This gives a description of the developing map of the homogeneous
properly convex foliated structure in geometric terms using only
``convexity'' properties of the curve $\xi: \POR \to \Flag(\RR^4)$. 

\subsection{Convex Curves and Convex Representations}
\begin{defi}\label{def_curvconv}
  A curve 
  \begin{equation*}
    \xi^1: S^1 \longrightarrow \PkR{n-1}
  \end{equation*}
  is said to be \emph{convex} if for, every $n$-tuple $(t_1, \dots,
  t_n)$ of pairwise distinct points $t_i \in \POR$, we have 
  \begin{equation*}
    \bigoplus_{i=1}^n \xi^1(t_i) = \RR^n.
  \end{equation*}
\end{defi}
In \cite{Labourie_2003, Guichard_2005, Guichard_2006} convex curves are called
hyperconvex. They were previously known and studied under the name of convex
curves (see e.g. \cite{Schoenberg_1954}) and we stick to this terminology.  A
convex curve $\xi^1: S^1 \to \PkR{2}$ is precisely an injective curve which
parametrizes the boundary of a strictly convex domain in $\PkR{2}$.

We are interested in convex curves $\partial \Gamma \rightarrow \PkR{n-1}$
which are equivariant with respect to some representation $\rho: \Gamma \to
\PSL_n(\RR)$.
\begin{defi}
  A representation $\rho: \Gamma \to \PSL_n(\RR)$ is said to be \emph{convex}
  if there exists a $\rho$-equivariant (continuous) convex curve $\xi^1:
  \partial\Gamma \to \PkR{n-1}$.
\end{defi}

Convex representation are deeply related with the Hitchin component:
\begin{thm}[Labourie \cite{Labourie_2003}, Guichard \cite{Guichard_2006}]
  \label{thm_hiteqconvex}
  The Hitchin component $\mathcal{T}^n(\Gamma)$ is the set of (conjugacy
  classes of) convex representations $\rho: \Gamma \to \PSL_n(\RR)$.
\end{thm}

Let us recall some facts and properties about $\rho$-equivariant convex
curves.

\begin{defi}
  A curve $\xi = (\xi^1, \dots, \xi^{n-1}): S^1 \to \Flag(\RR^n)$ is
  \emph{Frenet} if
  \begin{enumerate}
  \item For every $(n_1, \dots, n_k)$ with $\sum_{i=1}^k n_i = n$ and every
    $x_1,\dots, x_k \in S^1$ pairwise distinct, the following sum is direct:
    \begin{equation*}
      \bigoplus_{i=1}^{k} \xi^{n_i}(x_i) = \RR^n.
    \end{equation*}
  \item For every $(m_1, \dots, m_k)$ with $\sum_{i=1}^{k} m_i =m \leq n$ and
    for every $x\in S^1$
    \begin{equation*}
      \lim_{(x_i) \to x} \bigoplus_{i=1}^k \xi^{m_i}(x_i) = \xi^m(x),
    \end{equation*}
    the limit is taken over $k$-tuples $(x_1, \dots, x_k)$ of pairwise
    distinct $x_i$.
  \end{enumerate}
\end{defi}
If $\xi$ is Frenet, the curve $\xi^1$ is convex and entirely determines $\xi$
by the limit condition.

\begin{thm}[Labourie \cite{Labourie_2003}]
  \label{thm_convexFrenet}
  If a representation $\rho$ is a convex, then there exists a (unique)
  $\rho$-equivariant Frenet curve $\xi = (\xi^1, \dots, \xi^{n-1}): \partial
  \Gamma \to \Flag(\RR^n)$.
\end{thm}

Frenet curves satisfy a duality property
\begin{thm}[\cite{Guichard_2005} Th\'eor\`eme~5]
  \label{thm_convdual}
  Let $\xi = (\xi^1, \dots, \xi^{n-1}) : S^1 \rightarrow \Flag(V)$ be a Frenet
  curve. Then the dual curve $\xi^\perp : S^1 \rightarrow \Flag( V^*)$, $t
  \mapsto ( \xi^{n-1, \perp}(t), \dots, \xi^{1,\perp}(t))$ is also Frenet.
\end{thm}
\begin{remark}
  This duality of Frenet curves is more natural in the context of positive
  curves into flag variety for which we refer the reader to
  \cite{Fock_Goncharov} where the connection between convex curves into
  $\PkR{n}$ and positive curves into $\Flag(\RR^n)$ is discussed.
\end{remark}

This means that we can check if a curve is Frenet indifferently by
investigating sums or intersections of vectors spaces.  From this we can
deduce
\begin{prop}[\cite{Anisov_1998}]
  \label{prop_hierarchy}
  Let $\xi = (\xi^1, \xi^2, \xi^3): S^1 \to \Flag(\RR^4)$ be a Frenet curve.
  For $t \in S^1$, let $\xi_t: S^1 \to \Flag(\xi^3(t))$ be the curve defined
  by
  \begin{eqnarray*}
    \xi_t: S^1 &\longrightarrow& \Flag(\xi^3(t)) \\
    t' &\longmapsto& \Big\{
    \begin{array}{cl}
      ( \xi^3(t) \cap \xi^2(t') , \xi^3(t) \cap \xi^3(t') ) &
      \text{if } t'\neq t\\
      ( \xi^1(t) , \xi^2(t) ) & \text{otherwise.} 
    \end{array}
  \end{eqnarray*}
  Then $\xi_t$ is a Frenet curve.
\end{prop}

\begin{proof}
  By the duality property, it suffices to check:
  \begin{itemize}
  \item for all $t^{\prime}_{1} \neq t^{\prime}_{2}$, $\displaystyle
    \xi^{1}_{t} (t^{\prime}_{1} ) \cap \xi^{2}_{t}( t^{\prime}_{2} ) = \{0\}$,
  \item for all pairwise distinct $t^{\prime}_{1}$, $t^{\prime}_{2}$ and
    $t^{\prime}_{3}$, $\displaystyle \xi^{2}_{t} (t^{\prime}_{1} ) \cap
    \xi^{2}_{t}( t^{\prime}_{2} ) \cap \xi^{2}_{t}( t^{\prime}_{3} )= \{0\}$,
  \item and for all $t'$, $\displaystyle \lim_{ t^{\prime}_{1} \neq
      t^{\prime}_{2} \to t'} \xi^{2}_{t} (t^{\prime}_{1} ) \cap \xi^{2}_{t}(
    t^{\prime}_{2} ) = \xi^{1}_{t} ( t') $.
  \end{itemize}
  These three properties follow from the corresponding properties for $\xi$.
\end{proof}

\subsection{The Properly Convex Foliated Structure of a Convex
  Representation}
\begin{thm}\label{thm:reptostruct}
  Let $\rho: \Gamma \to \PSL_4(\RR)$ be a convex representation. Then there
  exists a developing pair $(\dev, \hol)$, with holonomy homomorphism $\hol
  =\rho\circ p: \overline{\Gamma} \to \PSL_4(\RR)$ defining a properly convex
  foliated projective structure on $M$.
\end{thm}

In fact we construct a section $s: \mathcal{T}^4( \Gamma) \rightarrow
\mathcal{P}_{pcf}(M)$ of the holonomy map. This section will be automatically
continuous, injective with image a connected component of $\mathcal{P}_{pcf}(
M)$.

\begin{proof}
  We want to use formula~\eqref{eq_devcurv}.  Since the holonomy homomorphism
  $\hol$ factors through $\Gamma$, the developing map $\dev$ will be defined by a
  $\rho$-equivariant map $\dev: \overline{M} \to \PTR$.  Let
  \begin{equation*}
    \xi = (\xi^1, \xi^2, \xi^3): \partial \Gamma \to \Flag(\RR^4)
  \end{equation*}
  be the $\rho$-equivariant Frenet curve.  Let $( \xi^{1}_{t})_{ t \in
    \partial \Gamma }$ be the $\rho$-equivariant family of continuous curves
  \begin{eqnarray*}
    \xi^{1}_{t}: \partial \Gamma & \longrightarrow& \xi^3(t) \subset \PTR\\
    t' &\longmapsto& \xi^3(t) \cap \xi^2(t') \quad \text{if } t'\neq t\\
    t &\longmapsto& \xi^1(t). 
  \end{eqnarray*}
  By Proposition~\ref{prop_hierarchy}, for every $t$, the curve $\xi^{1}_{t}$
  is convex and hence bounds a properly convex domain $C_t \subset \xi^3(t)$.
  Note that, as in the homogeneous example above, the tangent line to $C_t$ at
  $\xi^1(t)$ is $\xi^2(t)$ and the tangent line to $C_t$ at $\xi^{1}_{t}(t')$,
  $t'\neq t$, is $\xi^3(t') \cap \xi^3(t)$ (see
  Figure~\ref{fig_devmap}). In particular, the
  point
  \begin{equation*}
    \xi^{1}_{t'}(t) = \xi^2(t) \cap \xi^3(t')  = \xi^2(t) \cap
    \left(\xi^3(t') \cap \xi^3(t)\right) 
  \end{equation*}
  is the intersection of the two tangent lines.
  
  We define the developing map by
  \begin{eqnarray*}
    \dev: \partial \Gamma^{3+} \simeq \overline{M} &\longrightarrow & \PTR \\
    (t_+, t_0, t_-) &\longmapsto& \overline{\xi^1(t_+)
      \xi^1_{t_+}(t_-)} \cap \overline{\xi^1_{t_-} (t_+) \xi^1_{t_+}(t_0)}
  \end{eqnarray*}
  Then $\dev$ is a $\rho$-equivariant, it is continuous and injective so it is
  a homeomorphism onto its image. The image of the (weakly) stable leaf $t$ is the
  proper convex set $C_t$ and the image of the geodesic $(t_+, t_-)$ is
  contained in the projective line $\overline{\xi^1(t_+) \xi^1_{t_+}(t_-) }$.
  Therefore, the pair $(\dev, \hol)$ defines a properly convex foliated
  projective structure on $M$.
\end{proof}

Similar to the foliated structure with non-convex leaves in the homogeneous
example above we also get the following: 
\begin{thm}
  \label{thm:nonprop}
  Let $\rho: \Gamma \to \PSL_4(\RR)$ be a convex representation. Then
  there exists a developing pair $(\dev', \hol)$, with holonomy
  homomorphism $\hol =\rho\circ p: \overline{\Gamma} \to \PSL_4(\RR)$
  defining a foliated projective structure on $M$ which is not convex.
\end{thm}

As above this could be stated as the existence of a section of $\hol$
from $\mathcal{T}( \Gamma)$ to $\mathcal{P}_{f}(M)$.

\begin{proof}
  Let $\xi = (\xi^1, \xi^2, \xi^3) : \partial \Gamma \to \Flag(\RR^4)$ be the
  $\rho$-equivariant Frenet curve.  We define the $\rho$-equivariant
  developing map
  \begin{eqnarray*}
    \dev': \partial\Gamma^{3+} \simeq \overline{M} &\longrightarrow & \PTR\\
    (t_+, t_0, t_-) & \longmapsto& \xi^3(t_+) \cap \xi^3(t_0) \cap \xi^3(t_-).
  \end{eqnarray*}
  By Proposition~\ref{prop_hierarchy} $\dev'$ is well defined. It is also
  continuous and locally injective, hence a local homeomorphism.  The image of
  a geodesic leaf $(t_+, t_-) \in \partial \Gamma^{(2)}$ is contained in the
  projective line $\xi^3(t_+) \cap \xi^3(t_-)$.  The image of the (weakly) stable leaf
  $t_+$ is contained in the projective plane $\xi^3(t_+)$, but it cannot be
  convex since it is the complementary subset in $\xi^3(t_+)$ of the closure
  of the convex set $\dev(t_+)$ given in the preceding Theorem.
\end{proof}

Note that $\dev$ is a global homeomorphism whereas $\dev'$ is only a
local homeomorphism, two points of $\partial \Gamma^{3+}$ have the
same images under $\dev'$ if, and only if, they differ by a
permutation.

\subsection{Domains of Discontinuity}\label{sec_domdisc}
The above foliated projective structures $(\dev,\hol)$ and $(\dev',\hol)$
appear naturally when we consider domains of discontinuity for the action of
$\rho(\Gamma)$ on $\PTR$.  The action of $\rho(\Gamma)$ on $\PTR$ is not free
or proper, since $\rho(\gamma)$ has fixed points for every $\gamma \in
\Gamma$. 

But if we remove the ruled surface (discriminant) $\xi^2(\partial
\Gamma) \subset \PTR$, the complement $\PTR \backslash \xi^2(\partial \Gamma)$
has two connected components $\Lambda_\xi = \dev'(\overline{M})$ and
$\Omega_\xi = \dev(\overline{M})$. Namely, the image
$\dev^\prime( \overline{M})$ is contained in $\PTR \backslash \xi^2( \partial \Gamma)$  
and, using the Frenet property of $\xi$, the boundary of 
$\dev^\prime( \overline{M})$ is 
$\xi^2( \partial \Gamma)$. This implies that $\dev^\prime( \overline{M})$ is 
one connected component of $\PTR \backslash \xi^2(\partial \Gamma)$. Furthermore, 
by Proposition~\ref{prop_hierarchy} and Figure~\ref{fig_devmap},  
$\dev( \overline{ M})$ is the complementary of the closure of 
$\dev^\prime( \overline{M})$, hence the other connected component of 
$\PTR \backslash \xi^2(\partial \Gamma)$.

In particular, $\rho(\Gamma)$ acts freely
and properly discontinuously on $\Lambda_\xi$ and on $\Omega_\xi$.  If $t$ is
a (weakly) stable leaf , then $\dev(t) = \Omega_\xi \cap \xi^3(t) = C_t$ and $\dev'(t)
= \Lambda_\xi\cap \xi^3(t) = \xi^3(t)\backslash \overline{C_t}$.


\section{From Properly Convex Foliated Structures to Convex
  Representations}
\label{sec_proptorep}

In this section we will prove the following
\begin{thm}\label{thm:struct_hol}
  The holonomy representation $\hol: \overline{\G} \to \PGL_4(\RR)$ of a
  properly convex foliated projective structure factors through a convex
  representation $\rho: \G \to \PSL_4(\RR)$ and the foliated projective
  structure on $M$ is equivalent to the one associated to $\rho$ in
  Theorem~\ref{thm:reptostruct}.
\end{thm} 
Basically, we will construct an equivariant continuous curve $ \partial \Gamma
\rightarrow \PP^3(\RR)^*$ and show that it is convex.

\subsection{Maps Associated to Foliated Projective
  Structures}\label{sec_firstconstr}
\begin{prop}
  \label{prop:fol_maps}
  Let $(\dev, \hol)$ be the developing pair of a foliated projective structure
  on $M$.  Then the two maps
  \begin{eqnarray*}
    \xi^3 : \widetilde{\mathcal{F}} \simeq \widetilde{ \partial \Gamma}
    &\longrightarrow& \PTRd\\
    \mathcal{ D} : \widetilde{ \mathcal{G}} \simeq \widetilde{\partial
      \Gamma}^{(2)}_{[0]} &\longrightarrow& \textup{Gr}_2^4 (\RR),
\end{eqnarray*}
defined by 
\begin{center}
$ \xi^3 (f) = \xi^3( t_+)$ is the projective plane
containing $\dev(f)$,  and \\
$\mathcal{D}( g)= \mathcal{D}( t_+, t_-)$ is the projective line
containing $\dev(g)$ 
\end{center}
are continuous and $\hol$-equivariant. 
Moreover 
\begin{enumerate}
\item the map $\xi^3$ is locally injective.
\item for all $( t_+ , t_- )$ in $\widetilde{ \mathcal{G}}$, $\mathcal{D}
  (t_+, t_-) \subset \xi^3( t_+)$.
\item for all $t_+$, the map $t_- \mapsto \mathcal{D}( t_+, t_-)$ is not
  locally constant.
  \end{enumerate}
\end{prop} 
\begin{proof}
  The developing pair $( \dev, \hol)$ defining a foliated projective structure
  on $M$ consists of the homomorphism
\begin{equation*}
  \hol : \overline{ \Gamma} = \pi_1 (M) \longrightarrow \PGL_4 (\RR)
\end{equation*}
and the $\hol$-equivariant local homeomorphism
\begin{equation*}
  \dev : \widetilde{ M} \longrightarrow \PTR
\end{equation*}
satisfying the following properties
\begin{enumerate}
\item for any geodesic $g$ in $\widetilde{ \mathcal{ G}}$, $\dev( g)$ is
  contained in a projective line, and
\item for every (weakly) stable leaf $f$ in $\widetilde{ \mathcal{ F}}$, $\dev(f)$ is
  contained in a projective plane.
\end{enumerate}

This shows that the above definitions of $\xi^3$ and $\mathcal{D}$ are meaningful. 
The properties of $\xi^3$ and $\mathcal{D}$ are direct consequences of the
fact that $\dev$ is a $\hol$-equivariant local homeomorphism.
\end{proof}

In fact we can be a little more precise about the injectivity of the map
$\xi^3$.

\begin{lem}
  \label{lem_xi3inj}
  For every $( t_+, t_-) \in \widetilde{ \partial \Gamma}^{(2)}_{[0]}$ we have
  $\xi^3 (t_+) \neq \xi^3( t_-)$.
\end{lem}

\begin{proof}
  Take $(t_+, t_-) \in  \widetilde{ \partial \Gamma}^{(2)}_{[0]}$. By 
  Lemma~\ref{lem_closedgeo}
  there exists an element $\g \in \overline{\Gamma}$ of zero translation
  such that $(\g^n t_+), (\g^n t_-)$ converge to $\tilde{t}_{+,\gamma}$. By the 
  above local injectivity, for big enough $n$
  \bqn
  \rho( \g)^n \xi^3(t_+) = \xi^3(\g^n t_+)\neq \xi^3(\g^n t_-) = 
  \rho( \g)^n \xi^3(t_-).
  \eqn
\end{proof}

\subsection{The Holonomy Action on Convex Sets}
\label{sec_actiononconv}
\begin{prop}
\label{prop:eigenvalues}
Suppose that $(\dev, \hol)$ is a developing pair defining a properly convex
foliated projective structure on $M$.  Then for every element $\overline{\g}
\in \overline{\G}-\left<\tau\right>$ of zero translation $\hol(\overline{\g})$
is diagonalizable over $\RR$ with all eigenvalues being of the same sign.
Moreover if $t_+ \in \widetilde{\partial \Gamma}$ is an attractive fixed point
of $\overline{\g}$, the eigenvalues of the action of $\hol(\overline{\g})$
restricted to $\xi^3(t_+)$ satisfy $|\lambda_+|\geq |\lambda_0|> |\lambda_-|$.
\end{prop}

\begin{remark}
  Eigenvalues of an element of $\hol(\overline{\g}) \in \PGL_4( \RR)$ are of
  course only defined up to a common multiple.  Our statements about
  eigenvalues will clearly be invariant by scalar multiplication.
  
  From Proposition~\ref{prop:eigenvalues} one could already conclude some
  properties of the action of the central element $\tau\in \overline{ \Gamma}$
  as it commutes with an $\RR$-diagonalizable element with at least two
  distinct eigenvalues.  We will not state any of these properties until we
  are indeed able to prove that $\hol( \tau)$ is trivial.
\end{remark}

\subsubsection{Some Observations}
Since the developing pair $(\dev, \hol)$ defines a properly convex foliated
projective structure we have that the image of each (weakly) stable leaf is convex.
\begin{nota}
  We denote by $C_{ t_+}$ the properly convex subset in $\xi^3 (t_+)$ equal to
  $\dev( t_+)$.
\end{nota}

Let $\overline{\g} \in \overline{\G}-\left<\tau\right>$ be an element of zero
translation and let $( t_+, t_-) \in \widetilde{\partial \Gamma}^{(2)}_{[0]}$
be a pair of an attractive and a repulsive fixed point of $\overline{\gamma}$.
Then $\hol(\overline{\g}) \in \PGL_4(\RR)$ stabilizes the projective plane
$\xi^3(t_+)$ and also the open properly convex set $C_{t_+} = \dev(t_+)
\subset \xi^3(t_+)$. Moreover $\hol(\overline{\g})$ stabilizes the projective
line $\mathcal{D}( t_+, t_-)$ containing the image $\dev( t_+, t_-)$ of the
geodesic $( t_+, t_-) \in \widetilde{\mathcal G}$.

\begin{lem}
  The action of $\hol(\overline{\g})$ on $\mathcal{D}( t_+, t_-) \subset
  \xi^3(t_+)$ has two eigenlines $x_+, x_-$ with eigenvalues $\lambda_+,
  \lambda_-$ satisfying:
\begin{equation*}
  \lambda_+ \lambda_- \geq 0 \text{ and } |\lambda_+|> |\lambda_-|.
\end{equation*}
\end{lem}

\begin{proof}
  Since $C_{t_+}$ is properly convex, $\overline{C}_{ t_+}$ does not contain
  any projective line; hence the intersection
  \begin{equation*}
    \mathcal{D}( t_+, t_-) \cap \partial C_{t_+}
  \end{equation*}
  consists of two points $x_+ = x_+( \overline{\gamma})$ and $x_-= x_-(
  \overline{\gamma})$ which are fixed by $\hol( \overline{\gamma})$.  The
  points $x_+$ and $x_-$ are the endpoints at $+\infty$ and at $-\infty$ of
  the segment $\dev( t_+, t_-)$.  So $x_+$ and $x_-$ are eigenlines for $\hol(
  \overline{\gamma})$ corresponding to real eigenvalues $\lambda_+ =
  \lambda_+( \overline{\gamma})$ and $\lambda_- = \lambda_-(
  \overline{\gamma})$ The segment $\dev( t_+, t_-)$ with endpoints $x_+$ and
  $x_-$ is stable by $\hol( \overline{\gamma})$. This implies that the two
  eigenvalues $\lambda_+$ and $\lambda_-$ are of the same sign,
  \begin{equation*}
    \lambda_+ \lambda_- >0.
  \end{equation*}
  The action of $\overline{\g}$ on the geodesic $( t_+, t_-) \subset
  \widetilde{M}$ corresponds to a positive time map of the geodesic flow. This
  implies that for every point $x$ in $\dev( t_+, t_-)$, the limit of the
  sequence $( \hol( \overline{\gamma})^n x)_{n \in \NN}$ is equal to $x_+$.
  This gives the inequality $|\lambda_+| > | \lambda_-|.$
\end{proof}

\begin{lem}
  \label{lem_notexchange}
  The action of $\hol( \overline{\gamma})$ does not interchange the two
  components of $C_{ t_+} - \mathcal{D}( t_+, t_-)$.
\end{lem}

\begin{proof}
  By proper convexity, the set $C_{ t_+} - \mathcal{D}( t_+, t_-)$ has indeed
  two connected components, $C_1$ and $C_2$ and $\hol( \overline{\gamma})$
  either exchanges them or send them into themselves.
  
  Since the image of a geodesic $(t_+, t_-)$ is an open segment, the
  restriction of $\dev$ to this geodesic is necessarily a homeomorphism onto
  its image. Fix a point $x$ in this geodesic $(t_+, t_-) \subset t_+$. In
  the leaf $t_+$ there is a neighborhood $U$ of the geodesic segment $[ x,
  \overline{ \gamma} \cdot x]$ such that $\dev_{|U}$ is a homeomorphism onto
  its image.
  
  The complementary of $(t_+, t_-)$ in the leaf $t_+$ has also two connected
  components, $D_1$ and $D_2$ and we have (up to reindexing) $\dev(U \cap D_1
  ) = \dev(U) \cap C_1$ and $\dev(U \cap D_2) = \dev(U) \cap C_2$. Moreover,
  as the action of $\overline{ \gamma}$ does not exchange $D_1$ and $D_2$,
  there are points $y$ (close to $x$) in $U \cap D_1$ such that $\overline{
    \gamma} \cdot y$ is in $U \cap D_1$. Hence the point $m = \dev(y)$ is in
  $C_1$ and its image $\hol( \overline{ \gamma}) \cdot m = \dev(
  \overline{\gamma} \cdot y)$ is also in $C_1$. This implies that $\hol(
  \overline{ \gamma}) \cdot C_1 = C_1$ and the same for $C_2$.
\end{proof}

\begin{lem}
  The action of $\hol(\overline{\g})$ on $\xi^3(t_+)$ has a third eigenline
  $x_0$ with eigenvalue $\lambda_0$ being of the same sign as $\lambda_+$ and
  $\lambda_-$.
\end{lem}
\begin{proof}
  By Lemma~\ref{lem_notexchange}, the left tangent line to $C_{t_+}$
  at the point $x_+$ and the right tangent line to $C_{t_+}$ at $x_-$ are
  preserved by $\hol( \overline{\gamma})$.  Their intersection is a third
  eigenline $x_0=x_0(\overline{\g})$ in $\xi^3( t_+)$, fixed by $\hol(
  \overline{\gamma})$. The point $x_0$ does not lie on the projective line
  $\mathcal{D}( t_+, t_-)$. Therefore the restriction of $\hol(
  \overline{\gamma})$ to $\xi^3( t_+)$ is diagonalizable over $\RR$. Moreover,
  since $\hol(\overline{\gamma})$ does not interchange the two connected
  components of $C_{ t_+} - \mathcal{D}( t_+, t_-)$, the third eigenvalue
  $\lambda_0 = \lambda_0( \overline{\gamma})$ is of the same sign as
  $\lambda_+$ and $\lambda_-$.
\end{proof}
Summarizing we have the following
\begin{lem}
  \label{lem_gammadiag}
  For all $\overline{\gamma}\in \overline{ \Gamma}$ of zero translation, the
  element $\hol( \overline{\gamma})$ is diagonalizable over
  $\RR$ with all eigenvalues being of the same sign. In particular
  $\hol( \overline{\gamma}) \in \PSL_4(\RR)$.
\end{lem}
\begin{proof}
  Let $(t_+, t_-) \in \widetilde{\partial \Gamma}^{(2)}_{[0]}$ be the chosen
  pair of an attractive and a repulsive fixed point for $\overline{\gamma}$.
  Then the restriction of $\hol( \overline{\gamma})$ to $\xi^3(t_+)$ is
  diagonalizable over $\RR$ with three eigenvalues of the same sign. Applying
  the above arguments to $\overline{\gamma}^{-1}$ the same holds for the
  restriction of $\hol( \overline{\gamma})$ to $\xi^3(t_-)$.  Note that the
  inclusion $\xi^3( t_+) \cap \xi^3 (t_-) \subset \xi^3( t_-)$ is strict by
  Lemma~\ref{lem_xi3inj}, so the fourth eigenvalue of $\hol(\overline{\g})$ is
  the one corresponding to the eigenline of $\hol(\overline{\gamma})$ in
  $\xi^3(t_-)$ which is not contained in the intersection $\xi^3( t_+) \cap
  \xi^3( t_-)$.
\end{proof}

Keeping the same notation as above we get the following
\begin{lem}
  \label{lem_ineq+0}\label{lem_ineq-0}
  For every nontrivial element $\overline{\gamma}\in \overline{\Gamma}$ of
  zero translation, the inequalities $| \lambda_- | < | \lambda_0 | \leq | \lambda_+ |$ hold.
\end{lem}

\begin{proof}
  Suppose on the contrary that $| \lambda_0 | > | \lambda_+ |$ and let $x_0
  \in \xi^3( t_+)$ be the corresponding eigenline for $\lambda_0$. The convex
  $C_{ t_+}$ contains a neighborhood of the segment $\dev( t_+, t_-)$. The
  image of this neighborhood under $\hol(\overline{\gamma}^n)_{n\in \ZZ}$ will
  be a sector bounded by the two projective lines $\overline{ x_0 x_+}$ and
  $\overline{ x_0 x_-}$. Therefore the $\hol(\overline{\g})$-invariant convex
  $C_{ t_+}$ has to be this sector, contradicting the hypothesis that $C_{
    t_+}$ is a properly convex set.

%


  For any geodesic $(t_+,t)$ in the (weakly) stable leaf $t+$ consider the intersection of $\mathcal{D}(t_+,t)$ with a tangent line $L$ to
  $C_{t_+}$ in $x_-$.  By continuity and equivariance of $\mathcal{D}$
  \begin{equation*}
    \lim_{ n\to+\infty} \hol( \overline{\gamma})^{-n} (\mathcal{D}(t_+,t)\cap L )
    = \mathcal{D}(t_+,t_-) \cap L =x_-.
  \end{equation*}
 The restriction of $\hol( \overline{\gamma})$ to $L$ is diagonalizable with
  two real eigenvalues $\lambda_-$ and $\lambda_0$, which satisfy hence the inequality
  $| \lambda_- | < | \lambda_0 |$
\end{proof}
With this Proposition~\ref{prop:eigenvalues} is proved.

\subsubsection{Two Cases for the Action of $\hol(\overline{\gamma})$}
\label{sec_twocases}
The two possible cases $\lambda_+ = \lambda_0$ and $|\lambda_+|>|\lambda_0|$
for the eigenvalues of $\hol(\overline{\g})$ of a nontrivial element
$\overline{\gamma}\in \overline{\Gamma}$ of zero translation have different
consequences for the convex sets.
We continue to use the notation from the previous paragraph.

\begin{figure}[htbp]
  \centering \input{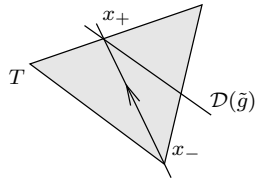}
  \caption{A triangle}
  \label{fig_triangle}
\end{figure}

\smallskip
\underline{\textbf{Case (T)}:} $\lambda_+ = \lambda_0$.\\
There is a unique $\hol( \overline{\gamma})$-invariant projective line $T=
\overline{x_0 x_+}$ through the point $x_+$ and distinct from $\mathcal{D}(
t_+, t_-)$, which is the unique tangent line to $C_{ t_+}$ at $x_+$.  Any
projective line $L$ containing the point $x_-$ is invariant by $\hol(
\overline{\gamma})$. Since $\hol( \overline{\gamma})$ fixes $T$ pointwise, by
Lemma~\ref{lem_ineq-0} it acts on $L$ with one repulsive fixed point $x_-$ and
the attractive fixed point being $L \cap T$. Therefore the only proper
invariant convex sets in $L$ with non empty interior relative to $L$ are the
two segments with endpoints $x_-$ and $L \cap T$. In particular $L \cap
C_{t_+}$ is such a segment for any $L$, and $C_{ t_+}$ is a union of such
segments. This means that $C_{t_+}$ is a triangle with one side supported by
$T$ and the third vertex being $x_-$ (Figure~\ref{fig_triangle}).

We will finally show that case (T) will not occur.

\smallskip
\underline{\textbf{Case (C)}:} $| \lambda_0 | < | \lambda_+  |$.\\
In this case the third eigenline $x_0$ in the projective plane $\xi^3( t_+)$
is the eigenspace corresponding to $\lambda_0$.  The projective line
$\overline{x_0 x_+}$ is the unique $\hol(\overline{\g})$-invariant projective
line through $x_+$ which is contained in $\xi^3(t_+)$ and different from
$\mathcal{D}( t_+, t_-)$.  Similarly $\overline{x_0 x_-}$ is the unique
$\hol(\overline{\g})$-invariant projective line through $x_-$ which is
contained in $\xi^3(t_+)$ and different from $\mathcal{D}( t_+, t_-)$.  In
particular the tangent lines to $C_{ t_+}$ at $x_+$ and $x_-$ are unique.
\begin{figure}[htbp]
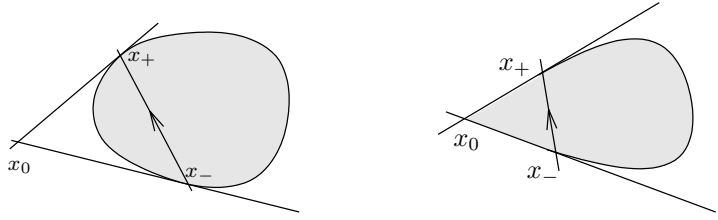

  \centering \input{niceconv.tex} \hspace*{4em}
  \input{niceconvdeg.tex}
  \caption{A ``nice'' convex and a degenerate proper convex}
  \label{fig_convtriang}
\end{figure}
We might have degenerate cases as in Figure~\ref{fig_convtriang}.  But as we
will see the degenerate case actually never occurs.

\subsection{Defining the Map $\xi^1$}\label{sec_defxi1}
\begin{prop}\label{prop:xi1}
  Let $(\dev, \hol)$ be a developing pair defining a properly convex foliated
  projective structure on $M$. Then there exists a continuous
  $\hol$-equivariant map
  \begin{equation*}
    \xi^1 : \widetilde{ \partial \Gamma} \rightarrow \PTR
  \end{equation*}
  such that $\xi^1(t_+) \in {\mathcal D}(g)$ for all $t_+
  \in \widetilde{ \partial \Gamma}$ and all $g\subset t_+$.
\end{prop}
We continue to use the notations from the previous section.
\begin{lem}
  \label{lem_alterxi1}
  We have the following alternative:
  \begin{itemize}
  \item for all $t_+\in\widetilde{ \partial \Gamma}$ the intersection
    \begin{equation*}
      \bigcap_{ g \subset t_+} \mathcal{D}(g) = \underset{ (\tau t, t_+, t)
      \textrm{ oriented} }{\bigcap_{ t \in \widetilde{ \partial \Gamma}}}
      \mathcal{D}( t_+, t) = \emptyset 
    \end{equation*}
    is empty, or
  \item for all $t_+\in\widetilde{ \partial \Gamma}$ the intersection 
$\bigcap_{ g \subset t_+} \mathcal{D}(g)$ is a point in $\PTR$.
  \end{itemize}
\end{lem}
\begin{proof}
  The injectivity of the map $\mathcal{D}$ (Proposition~\ref{prop:fol_maps})
  implies that the intersection $\bigcap_{ g \subset t_+} \mathcal{D}(g)$ is
  either empty or a point.  The continuity of $\mathcal{D}$ implies that the
  $\overline{ \Gamma}$-invariant set
 \begin{equation*}
    \{ t_+ \in \widetilde{ \partial \Gamma} \mid \bigcap_{ g \subset t_+}
    \mathcal{D}(g) = \emptyset \}
  \end{equation*}
  is open. Since the action of $\overline{ \Gamma}$ on $\widetilde{ \partial
    \Gamma}$ is minimal (Lemma~\ref{lem_minstable}) this set is either empty or equal to 
$\widetilde{\partial
    \Gamma}$ this proves the claim.
\end{proof}

\begin{lem}\label{lem_Dintersects}
  Let $\overline{ \gamma}$ be of zero translation and $t_+$ an attractive 
  fixed point. For all geodesics $g_t = ( t_+, t)$ in the (weakly) stable leaf $t_+$, 
  the projective
  line $\mathcal{D}(g_t)$ contains the point $x_+$.
\end{lem}

\begin{proof}
  As observed in the proof of Lemma~\ref{lem_ineq+0} the sequence 
$(\hol(
  \overline{\gamma})^{-n} \mathcal{D}( g_t) )_{ n\in \NN}$ converges to 
$\overline{ x_+ x_-}$. 
Let $L$ be a $\hol( \overline{\gamma})$-invariant projective line
  tangent to $C_{ t_+}$ at $x_+$, then
  \begin{equation*}
    \lim_{ n \to + \infty } \hol( \overline{\gamma})^{-n} ( \mathcal{D}(g_t)
    \cap L) = L \cap \overline{x_+ x_-}= x_+ 
  \end{equation*}
  Since the restriction of $\hol( \overline{\gamma})$ to $L$ is diagonalizable
  with two real eigenvalues $\lambda_+$ and $\lambda_0$ satisfying the
  inequality $| \lambda_0 | \leq | \lambda_+ |$ and with $x_+$ corresponding
  to $\lambda_+$, this implies
  \begin{equation*}
    \mathcal{D}(g_t) \cap L = x_+,
  \end{equation*}
  and in particular, $x_+ \in \mathcal{D}(g_t)$.
\end{proof}
Lemma~\ref{lem_Dintersects} implies that we are in the second case of
Lemma~\ref{lem_alterxi1}.
\begin{defi}
  \label{def_xi1}
  For all $t$ in $\widetilde{ \partial \Gamma}$ we define $\xi^1(t)$ to be the
  common intersection of the projective lines $\mathcal{D}(g)$ for $g$ in the
  leaf $t$:
  \begin{equation*}
    \xi^1(t):= \bigcap_{ g \subset t} \mathcal{D}(g).
  \end{equation*}
\end{defi}

\begin{proof}[Proof of Proposition~\ref{prop:xi1}]
  Continuity and equivariance of $\xi^1$ follow from the corresponding
  properties of $\mathcal{D}$. 
From the very definition, for all $t\in \widetilde{ \partial \Gamma}$ and all $g \subset t$, 
$\xi^1( t) \subset \mathcal{D}(g) \subset \xi^3(t)$.
\end{proof}

\subsection{Intersections of $\xi^3$}
\label{sec_posexis}
As we already indicated, our ultimate goal is to show that the representation
$\hol : \overline{ \Gamma} \rightarrow \PGL_4( \RR)$ factors through a
representation of $\Gamma$ and that the induced curve $\xi^3 : \partial \Gamma
\rightarrow \PTRd$ is convex. In this paragraph we establish some facts about
the possible intersections of $\xi^3$.
Since $\xi^3$ is not constant, the intersection
\begin{equation*}
  \bigcap_{ t \in \widetilde{ \partial \Gamma}} \xi^3( t) \subset \PTR
\end{equation*}
is either empty, a point or a projective line.

\begin{lem}
  \label{lem_interxinotdev}
  If the intersection $\bigcap_{ t \in \widetilde{ \partial \Gamma} } \xi^3(t)$ is a projective line, this line  
  does not meet the image $\dev( \widetilde{M})$.
\end{lem}

\begin{proof}
  Suppose that the intersection $L =\bigcap_{ t \in \widetilde{ \partial
      \Gamma} } \xi^3(t)$ is a projective line and that there is a point $m$ in
  $\widetilde{M}$ such that $\dev(m)$ belongs to $L$. Then, for a small enough
  neighborhood $U$ of $m$, $\dev(U)$ will, in some affine chart, be contained in
  one of the sectors bounded by the two planes $\xi^3(t)$, $\xi^3(t')$
  (Figure~\ref{fig_devU}), contradicting that $\dev$ is a local
  homeomorphism.
  \begin{figure}[htbp]
    \centering \input{devU.tex}
    \caption{$\dev(U)$.}
    \label{fig_devU}
  \end{figure}
\end{proof}

\begin{lem}
  \label{lem_xilocline}
  Suppose that there is a non-empty open subset $U$ in $ \widetilde{ \partial
    \Gamma}$ such that $ L= \bigcap_{ t \in U } \xi^3(t)$ is a projective
  line, then
  \begin{equation*}
    \bigcap_{ t \in  \widetilde{ \partial \Gamma} } \xi^3(t) = L.
  \end{equation*}
\end{lem}

\begin{proof}
  The subset
  \begin{equation*}
    \big \{ t \in \widetilde{ \partial \Gamma} \mid \exists \text{ open } U_t
    \ni t \text{ such that } \bigcap_{ s \in  U_t } \xi^3(s) \text{ is a
    projective line} \big\}
  \end{equation*}
  is a non-empty, open and $\overline{ \Gamma}$-invariant subset of
  $\widetilde{ \partial \Gamma}$, so by minimality it equals
  $\widetilde{\partial \Gamma}$. Note that by the local injectivity of
  $\xi^3$ the intersection $\bigcap_{ s \in U_t } \xi^3(s)$ is independent of
  the choice of the open $U_t$. Therefore we have a locally constant
  continuous map
  \begin{eqnarray*}
    \widetilde{ \partial \Gamma} \longrightarrow \Gr_2^4 ( \RR), \quad
    t  \longmapsto  \bigcap_{ s \in  U_t } \xi^3(s).
  \end{eqnarray*}
  This map is constant equal to $L$.
\end{proof}

\begin{lem}
  \label{lem_xilocpoint}
  If there exists a non-empty open subset $U$ of $\widetilde{ \partial
    \Gamma}$, such that
  \begin{equation*}
    x=    \bigcap_{ t \in  U } \xi^3(t)
  \end{equation*}
  is a point in $\PTR$, then
  \begin{equation*}
    \bigcap_{ t \in  \widetilde{ \partial \Gamma} } \xi^3(t) \text{ is equal
    to }x.
  \end{equation*}
\end{lem}

\begin{proof}
  We argue along the lines of the proof of the preceding lemma. The set
  \begin{equation*}
    \big \{ t \in \widetilde{ \partial \Gamma} \mid \exists \text{ open } U_t
    \ni t \text{ such that } \bigcap_{ s \in  U_t } \xi^3(s) \text{ is a
    point} \big\}
  \end{equation*}
  equals $\widetilde{ \partial \Gamma}$.  Since by the previous lemma the
  intersection $\bigcap_{ t \in V } \xi^3(t)$ is never a projective line for
  $V$ a non-empty open subset of $\widetilde{ \partial \Gamma}$ we get a well
  defined locally constant map
  \begin{eqnarray*}
    \widetilde{ \partial \Gamma} \longrightarrow \PTR, \quad
    t  \longmapsto  \bigcap_{ s \in  U_t } \xi^3(s).
  \end{eqnarray*}
  This map is constant equal to $x$.
\end{proof}

\subsection{Semi-Continuity of $C_t$}
\label{sec_semicontconv}
The map that sends $t$ in $\widetilde{ \partial \Gamma}$ to the closure
$\overline{C_{t}}$ of the convex $C_t = \dev(t)$ has some semi-continuity
properties.

Fixing any continuous distance on $\PTR$ the space of compact
subsets of $\PTR$ is endowed with the Hausdorff distance.

\begin{lem}
  \label{lem_semcontCt}
  Let $(t_n)_{n \in \NN}$ be a sequence in $\widetilde{ \partial \Gamma}$
  converging to $t \in \widetilde{ \partial \Gamma}$ and such that the
  sequence of convex sets $( \overline{C}_{ t_n} )_{ n \in \NN}$ has a limit,
  then
  \begin{equation*}
    \lim_{ n \to + \infty } \overline{C}_{ t_n} \supset \overline{C}_t.
  \end{equation*}
\end{lem}

\begin{proof}
  It is enough to show that $C_t$ is contained in the limit $\lim_{ n \to +
    \infty } \overline{C}_{ t_n}$. For this, it is sufficient to show that,
  for all $x$ in $C_t$, there is a sequence $( x_n)_{n \in \NN}$ of $\PTR$
  converging to $x$ and such that $x_n$ belongs to $\overline{C}_{ t_n}$ for
  all $n$.
  
  Choose a point $m$ in $\widetilde{M}$ with $\dev(m) = x$. Since $( t_n)_{ n
    \in \NN}$ converges to $t$, there is a sequence $(m_n)_{n \in \NN}$
  converging to $m$ such that $m_n$ is contained in the the leaf $t_n$ for all
  $n$. The sequence $x_n = \dev( m_n)$ satisfies the above
  conditions.
\end{proof}

Lemma~\ref{lem_semcontCt} has the following refinement:

\begin{lem}
  \label{lem_semicont2}
  Under the same hypothesis as in the preceding lemma, suppose that $P$ is a
  projective line or a projective plane transversal to $\xi^3(t)$ (i.e. the
  intersection $P \cap \xi^3(t)$ is of the smallest possible dimension) and
  intersecting $C_t$, then
  \begin{equation*}
    \lim_{ n \to + \infty } \overline{C}_{ t_n}\cap P \supset
    \overline{C}_t\cap P.
  \end{equation*}
\end{lem}

\begin{remark}
  Instead of taking a fixed $P$ we could also work with a sequence $(P_n)_{ n
    \in \NN}$ converging to $P$.
\end{remark}

\begin{proof}
  Let $x = \dev(m)$ be in $C_t \cap P$ and let $U$, $V$ be neighborhoods of $m
  \in \widetilde{M}$ and $x \in \PTR$ respectively, such that the restriction
  of $\dev$ is a homeomorphism from $U$ onto $V$. The transversality condition
  implies that there exists a sequence $( m_n)_{ n \geq N}$, defined only for
  large enough $N$, such that $m_n \in t_n \cap \dev^{-1}( P \cap V)$. Now one
  can conclude as in the proof of Lemma~\ref{lem_semcontCt}.
\end{proof}

\subsection{Nontriviality of $\xi^1$}

\begin{prop}
  \label{prop_xi1notconstant}
  The map $\xi^1: \widetilde{\partial\Gamma} \to \PTR$ is not constant.
\end{prop}

\begin{proof}
  We argue by contradiction. Assume that $\xi^1$ is constant equal to $x\in
  \PTR$.
  
Before we give the formal argument which leads to a contradiction, let us summarize the idea of the
proof. We will show that all the convex sets $C_t$ have to be
triangles and that these triangles all share a common edge. Then we
show that for every triangle the vertex opposite to this edge is
contained in a fixed projective line. This forces the image of
the developing map to be contained in a two dimensional subspace,
which gives the desired contradiction.

  Let $\overline{\gamma} \in \overline{\Gamma}-\left<\tau\right>$ be an
  element of zero translation and $(\tilde{t}_{+,\gamma},\tilde{t}_{-,\gamma}) \in \widetilde{ \partial
    \Gamma}^{(2)}_{[0]}$ a pair of  fixed points. Since $x=\xi^1(\tilde{t}_{+,\gamma})$ is
  an eigenline for $\hol( \overline{\gamma})$ corresponding to the largest
  eigenvalue of $\hol( \overline{\gamma})_{ | \xi^3 (\tilde{t}_{+,\gamma})}$ and $x=\xi^1(\tilde{t}_{-,\gamma})$
  is an eigenline corresponding to the largest eigenvalue of $\hol(
  \overline{\gamma})^{-1}_{ | \xi^3( \tilde{t}_{-,\gamma})}$, we are necessarily in Case (T), so
  both convex sets $C_{ \tilde{t}_{+,\gamma}}$ and $C_{\tilde{t}_{-,\gamma}}$ are triangles.
  The intersection $L = \xi^3 ( \tilde{t}_{+,\gamma}) \cap \xi^3( \tilde{t}_{-,\gamma})$ is an eigenspace for
  $\hol( \overline{\gamma})$ (corresponding to the eigenvalue $\lambda_+(
  \overline{\gamma}) = \lambda_+( \overline{\gamma}^{-1})^{-1}$). Moreover $L$
  is the tangent line to $C_{ \tilde{t}_{+,\gamma}}$ at $x$ which is the same as the tangent
  line to $C_{ \tilde{t}_{-,\gamma}}$ at $x$. 

  Therefore for all $(\tilde{t}_{+,\gamma}, \tilde{t}_{-,\gamma})$ 
  we have $C_{\tilde{t}_{+,\gamma}} \cap \xi^3(\tilde{t}_{-,\gamma}) = \emptyset$.
  Since the set $\{(t_+, t_-) \in  \widetilde{\partial
   \Gamma}^{(2)}_{[0]}  \mid \xi^3(t_+) \cap C_{ t_-} =  \emptyset
   \}$ is closed,
  we deduce from Lemma~\ref{lem_closedgeo} that
  for all $t$, with $(\tilde{t}_{+,\gamma}, t) \in \widetilde{ \partial \Gamma}^{(2)}_{[0]}$
  the intersection $C_{\tilde{t}_{+,\gamma}} \cap \xi^3(t)$ is empty. 
In particular, the
  projective line $\xi^3(\tilde{t}_{+,\gamma}) \cap \xi^3( t)$ is the line tangent to $C_{
    \tilde{t}_{+,\gamma}}$ at $x$ and hence equals $L$. So there exists a non-empty open
  subset $U$ of $\widetilde{ \partial \Gamma}$ such that $L = \bigcap_{ t \in
    U} \xi^3(t)$. By Lemma~\ref{lem_xilocline} this implies $L = \bigcap_{ t
    \in \widetilde{ \partial \Gamma}} \xi^3(t)$.
\smallskip  
\noindent
{\bf Claim 1:} The segment $I \subset L$ corresponding to the side of the triangle $C_{ \tilde{t}_{+,\gamma}}$ is independent 
of $\overline{\gamma}$.

Note that, since $\tau^{-2g} = \Pi_{ i=1}^{g} [ a_i, b_i]$ and the restriction
of $\hol(a_i)$ and $\hol(b_i)$ to $L$ are trivial, the element $\hol(
\tau)^{-2g}_{ | L}$ is the identity. Therefore $\hol( \tau)_{| L} = 1$ since
$x$ is an eigenline for $\tau$. This implies that the sides of the triangles 
$C_{ \tilde{t}_{+,\gamma}}$ and $C_{ \tilde{t}_{-,\gamma}}$ in $L$ 
do not depend on the pair $(\tilde{t}_{+,\gamma}, \tilde{t}_{-,\gamma})$. 
We denote this segments
by $I_+( \overline{\gamma})$ and $I_-( \overline{\gamma})$. We now want to
show that they
are independent of $\overline{\gamma}$.  

For this consider $t \in
\widetilde{ \partial {\Gamma}}$ such that $(\tilde{t}_{+,\gamma}, t) \in \widetilde{ \partial
  \Gamma}^{(2)}_{[0]}$. Let $w\in L$ be a point and $P$ the projective plane
spanned by $w$, $x_-( \overline{\gamma})$ and $x_-(
\overline{\gamma}^{-1})$ (using the notation from
Section~\ref{sec_actiononconv}).  The plane $P$ is $\hol(
\overline{\gamma})$-invariant. If $w$ is in $I_+$, then
Lemma~\ref{lem_semicont2} implies that
  \begin{equation*}
    \lim_{ n \to + \infty} \hol( \overline{\gamma})^n \overline{C}_t \cap P
    \supset
    \overline{C}_{\tilde{t}_{+,\gamma}} \cap P.
  \end{equation*}
  This is possible if and only if $w$ belongs to $\overline{C}_t \cap P$.
  Applying this to all points $w\in I_+( \overline{\gamma})$ implies that if
  $t = t_{+,\overline{\gamma}'}$ is the attractive fixed point for another
  element $\overline{\gamma}' \in \overline{\Gamma}-\left<\tau\right>$ of zero
  translation, then $I_+( \overline{\gamma}') \supset I_+(
  \overline{\gamma})$.  Interchanging the roles of $\overline{\gamma}$ and
  $\overline{\gamma}'$ we get that $I = I_+( \overline{\gamma}) = I_+(
  \overline{\gamma}') = I_-( \overline{\gamma})$ is independent
  of $\overline{\gamma}$.\\
\smallskip
\noindent
{\bf Claim 2}: There exists a projective line $Q$ such that for every element
$\overline{\g} \in \overline{\G}-\left<\tau\right>$ of zero translation we
have $x_-( \overline{\gamma}) \subset Q$.
 
Let $y$ be the intersection point of $\xi^3( t_{+, \gamma'})$ with
the $\hol( \overline{\gamma})$-invariant projective line $Q$ spanned by $x_-(
\overline{\gamma})$ and $x_-( \overline{\gamma}^{-1})$. Then for any
projective line $D$ through $y$ which is contained in $\xi^3( t_{+, \gamma'})$, 
the intersection $\overline{ C}_{ t_{+, \gamma'}} \cap D$ is either empty or contains the intersection
point $L \cap D$. This shows that $y$ is contained in the sector with tip
$x_-( \overline{\gamma}')$, which is bounded by the projective lines
supporting the triangle $C_{ t_{+,\gamma'}}$ and which contains $C_{ t_{+,
  \gamma'}}$ (Figure~\ref{fig_ysector} shows what cannot happen).
  \begin{figure}[htbp]
    \centering \input{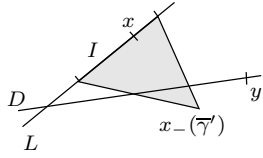}
    \caption{The point $y$ cannot be in the other sector}
    \label{fig_ysector}
  \end{figure}
  Similarly $y$ cannot be contained in the open triangle $C_{
    t_{+,\overline{\gamma}'}}$.  

  A short calculation, using this last condition
  for $ \overline{\gamma}^{\pm 1}$ and $
  \overline{\gamma}^{ \prime \pm 1}$, gives
  \begin{equation*}
    \overline{ x_-( \overline{\gamma}) x_-( \overline{\gamma}^{-1}) }
  = \overline{ x_-( \overline{\gamma}') x_-( \overline{\gamma}^{\prime
  -1}) }=Q,
  \end{equation*}
  and this holds for any $ \overline{\gamma}$ and $
  \overline{\gamma}^{ \prime }$.
  
  Hence the projective line $Q$ is invariant by every element of zero
  translation, and also by $\tau$ since $x_-( \overline{\gamma}^{\pm 1})$ are
  eigenlines for $\hol(\tau)$. Thus $Q$ is
  $\hol(\overline{\Gamma})$-invariant. We can define a $\hol$-equivariant continuous map
  \begin{eqnarray*}
    \xi_- : \widetilde{ \partial \Gamma} & \longrightarrow & Q \\
    t & \longmapsto & \xi_-(t) = \xi^3(t) \cap Q.
  \end{eqnarray*}
  Then the closed $\overline{ \Gamma}$-invariant set
  \begin{equation*}
    \{ (t, t') \in \widetilde{ \partial \Gamma}^{(2)}_{[0]} \mid
    \mathcal{D}( t, t') = x \oplus \xi_-( t) \}
  \end{equation*}
 contains the pairs $( \tilde{t}_{+,\gamma}, \tilde{t}_{-,\gamma})$ for all elements
  $\g\in \overline{\Gamma}-\{1\}$ of zero translation. Hence Lemma~\ref{lem_closedgeo}
implies that for any $(t,
  t')\in\widetilde{ \partial \Gamma}^{(2)}_{[0]} $
  \begin{equation*}
    \mathcal{D}( t, t') = x\oplus \xi_-(t)
  \end{equation*}
  which contradicts the local injectivity of the map $t' \mapsto \mathcal{D}
  (t, t')$ (see Proposition~\ref{prop:fol_maps}).
\end{proof}

%

%

\subsection{The Holonomy Factors}
In this paragraph we study the possible intersections of $\xi^3$ and show that
\begin{prop}
  \label{prop_holtautriv}
  The intersection $\bigcap_{t \in \widetilde{ \partial \Gamma}} \xi^3( t)$ is
  empty and $\hol( \tau) = \textup{Id}$.
\end{prop}

\begin{lem}
  \label{lem_interxi3notline}
  The intersection $\bigcap_{t \in \widetilde{ \partial \Gamma}} \xi^3( t)$
  cannot be a projective line.
\end{lem}

\begin{proof}
  Suppose $L=\bigcap_{t \in \widetilde{ \partial \Gamma}} \xi^3( t)$ is a
  projective line.   Then we have two representations
\begin{eqnarray*}
  \rho_L : \overline{ \Gamma} & \longrightarrow & \PGL( L)\\
  \rho_Q : \overline{ \Gamma} & \longrightarrow & \PGL( \RR^4 / L).
\end{eqnarray*}

Consider the set 
\bqn
C= \{ t\in \widetilde{\partial \Gamma} \mid \xi^1(t) \subset L\}
\eqn
This set is a $\overline{\Gamma}$-invariant closed set, which by
  minimality (Lemma~\ref{lem_minstable}) is either empty or equal to
  $\widetilde{\partial \Gamma}$. So we have to consider two cases 
\begin{itemize}
\item {\bf(Case 1):}  $\forall t\in \widetilde{\partial \Gamma} \quad
  \xi^1(t) \subset L$, 
\item {\bf(Case 2):} $\forall t\in \widetilde{\partial \Gamma} \quad
  \xi^1(t) \oplus L = \xi^3(t)$.
\end{itemize}
\smallskip
Let us first assume that we are in {\bf(Case 1)}. We show that in
this case the two representations $\rho_L$ and
$\rho_Q$ are Fuchsian with length functions  (see Appendix~\ref{sec_appendix})
$ \ell_{ \rho_L}< \ell_{\rho_Q}$. This will contradict Fact~\ref{fact_lengthteich}


Since we are in (Case 1) the curve $\xi^1$ is a continuous $\rho_L$-equivariant curve 
\begin{equation*}
  \xi^1 : \widetilde{\partial \Gamma} \rightarrow \PP(L) \subset \PTR
\end{equation*}
and $\xi^3$ is a continuous $\rho_Q$-equivariant curve
\begin{equation*}
  \xi^3: \widetilde{ \partial \Gamma} \rightarrow \PP( \RR^4/L) \subset
  \PTRd.
\end{equation*}

Let $\overline{\gamma}\in \overline{\Gamma}-\left<\tau\right>$ be an element
of zero translation with $(\tilde{t}_{+,\gamma}, \tilde{t}_{-,\gamma}) \in \widetilde{ \partial
  \Gamma}^{(2)}_{[0]}$ a pair of an attractive and a repulsive fixed point.
Then $\xi^1( \tilde{t}_{+,\gamma})$ is the eigenspace corresponding to the largest eigenvalue
(in modulus) of $\hol( \overline{\gamma})_{ |L}$. 
At least one $\rho_L( \overline{\gamma})$ is nontrivial since 
otherwise $\xi^1$ would be constant contradicting Proposition~\ref{prop_xi1notconstant}. 
For this particular element $\xi^1( \tilde{t}_{+,\gamma})$
is independent of the choice of $\tilde{t}_{+,\gamma}$ in its 
$\left<\tau\right>$-orbit. In other words $\xi^1( \tau \tilde{t}_{+,\gamma}) = \xi^1( \tilde{t}_{+,\gamma})$.
Thus $\xi^1$ is $\tau$-invariant, since the set
  \begin{equation*}
    \{ t \in \widetilde{ \partial \Gamma} \mid \xi^1( \tau t) = \xi^1( t) \}
  \end{equation*}
  is closed, $\overline{ \Gamma}$-invariant and non-empty. Consequently
  $\rho_L (\tau) = \textup{Id}_L$.  An analogous argument shows that $\rho_Q
  (\tau) = \textup{Id}_{\RR^4/L}$.
  
  In particular, the two representations $\rho_L$ and $\rho_Q$ factor as
  \begin{eqnarray*}
    \rho_L :  \Gamma &\longrightarrow& \PGL(L)\\
    \rho_Q : \Gamma &\longrightarrow& \PGL(  \RR^4/L).
  \end{eqnarray*}
  Since for any element $\overline{\g}$ of zero translation the eigenvalues of
  $\hol(\overline{\g})$ are of the same sign (by Lemma~\ref{lem_gammadiag}),
  these two representations have image in $\PSL_2( \RR)$. Moreover, they
  satisfy the hypothesis of Lemma~\ref{lem_curveteich} so $\rho_L$ and
  $\rho_Q$ are Fuchsian.
  
  Consider now the length functions of $\rho_L$ and $\rho_Q$
  (see Appendix~\ref{sec_appendix}).  

  Let $\gamma \in \Gamma-\{1\}$ and
  $\overline{\gamma}$ its lift in $\overline{ \Gamma}$ of zero translation
  with a pair of fixed points $(\tilde{t}_{+,\gamma},\tilde{t}_{-,\gamma}) \in \widetilde{ \partial
    \Gamma}^{(2)}_{[0]}$.
  By Lemma~\ref{lem_interxinotdev} the projective line $L$ does not meet
  $\dev( \widetilde{M})$. Therefore $L$ is tangent to $C_{ \tilde{t}_{+,\gamma}}$ at the point
  $x_+( \overline{ \gamma}) = \xi^1( \tilde{t}_{+,\gamma})$ and tangent to $C_{\tilde{t}_{-,\gamma}}$ at the
  point $x_+( \overline{ \gamma}^{-1}) = \xi^1( \tilde{t}_{-,\gamma})$. The eigenvalues of the
  restriction $ \hol( \overline{ \gamma})_{ | \xi^3( \tilde{t}_{+,\gamma})}$ are (we choose a
  lift of $\hol( \overline{ \gamma})$ with positive eigenvalues)
  \begin{equation*}
    \lambda_+( \overline{\gamma}) > \lambda_0( \overline{ \gamma}) >
    \lambda_-( \overline{ \gamma}) >0, 
  \end{equation*}
  where $\lambda_+$ and $\lambda_0$ are the eigenvalues corresponding to $L$
  which are distinct since $\rho_L$ is faithful. Similarly the eigenvalues of
  $ \hol( \overline{ \gamma}^{-1})_{ | \xi^3( \tilde{t}_{-,\gamma})}$ are:
  \begin{equation*}
    \lambda_+( \overline{\gamma}^{-1}) > \lambda_0( \overline{ \gamma}^{-1}) >
    \lambda_-( \overline{ \gamma}^{-1}) >0, 
  \end{equation*}
  with relations $\lambda_+( \overline{\gamma})^{-1} = \lambda_0( \overline{
    \gamma}^{-1})$ and $\lambda_0( \overline{\gamma})^{-1} = \lambda_+(
  \overline{ \gamma}^{-1})$. Therefore
  \begin{equation*}
    \ell_{ \rho_L}( \gamma) = \ln \Big( \frac{ \lambda_+(
      \overline{\gamma})}{ \lambda_0( \overline{ \gamma})} \Big)
    <  \ln \Big( \frac{ \lambda_-(
      \overline{\gamma}^{-1})^{-1}}{ \lambda_-( \overline{ \gamma})}
    \Big) = \ell_{ \rho_Q}( \gamma),
  \end{equation*}
  contradicting Fact~\ref{fact_lengthteich}. Therefore (Case 1)
  cannot occur.

\bigskip

If we are in {\bf(Case 2)}, that is for all 
  $t\in \widetilde{\partial\Gamma}$
  \begin{equation*}
     \xi^1(t) \oplus L = \xi^3(t),
  \end{equation*}
we show that $\rho_L$ and $\rho_Q$ have to be 
conjugate Fuchsian representations. Then we apply the analysis of 
 Section~\ref{sec_geomdesc} to get a contradiction.

  As above the
  $\rho_Q$-equivariant map $\xi^3: \widetilde{\partial\G} \to \PP(\RR^4/L)$ is
  $\tau$-invariant.  The representation $\rho_Q$ factors as $\rho_Q: \Gamma
  \to \PGL(\RR^4/L)$ satisfying the hypothesis of Lemma~\ref{lem_curveteich},
  so $\rho_Q$ is Fuchsian.
  
  Let $\overline{\gamma} \in \overline{ \Gamma}$ be an element of zero
  translation with fixed points $(\tilde{t}_{+,\gamma}, \tilde{t}_{-,\gamma})$ as before.  In $\xi^3( \tilde{t}_{+,\gamma})$ the
  line $L$ is spanned by $x_-( \gamma)$ and another invariant point $x_0(
  \gamma)$ distinct from $x_+( \gamma) = \xi^1( \tilde{t}_{+,\gamma})$. In particular $\rho_L(
  \overline{\gamma})$ is split over $\RR$ with distinct eigenvalues of the
  same sign. This implies that $\rho_L( \tau)$ is also split over $\RR$ in the
  same basis. By Lemma~\ref{lem_repgambarsl2} we get that
   $\rho_L( \tau) = \textup{Id}$. 

  
  Therefore $\rho_L$ factors through $\Gamma$ and again we have $\rho_L : \Gamma \rightarrow \PSL(L)$. Furthermore,
  for any nonzero $\gamma$, $\rho_L(\gamma)$ is split with two distinct real
  eigenvalues. Hence $\rho_L$ is faithful and discrete by
  Lemma~\ref{lem_fuchsian}, so $\rho_L$ is Fuchsian.
  
  We already observed that $\xi^1$ factors through a curve $\xi^1:\partial
  \Gamma \rightarrow \PTR$. For any $t\in \partial\Gamma$ the line $\xi^1(t)$
  is an eigenline for $\hol( \tau)$.  In particular $\hol( \tau)$ is split
  over $\RR$ and since $\hol( \tau)$ acts trivially on $L$ and on $\RR^4/L$,
  there are only two possibilities
  \begin{equation*}
    \hol( \tau) = \textup{Id} \textrm{ or } \hol( \tau) = \left(
      \begin{array}{cc}
        \textup{Id} & 0 \\ 0 & -\textup{Id}
      \end{array}\right).
  \end{equation*}
  In the first case the holonomy factors as $\hol: \Gamma \rightarrow \PSL_4(
  \RR)$, and in the second case we get $\hol: \widehat{ \Gamma} = \overline{
    \Gamma}/ \langle 2\tau \rangle \rightarrow \PSL_4( \RR)$. In both case
  what we already know about the eigenvalues implies that the holonomy lifts
  (non uniquely and up to taking a subgroup of index two) to $\widehat{
    \hol}: \Gamma \rightarrow \SL_4( \RR)$ or $\widehat{ \hol}:
  \widehat{\Gamma} \rightarrow \SL_4( \RR)$.
  
  Suppose that $\hol( \tau) = \textup{ Id}$. In a basis adapted to $L$, we
  have:
  \begin{eqnarray*}
    \hol : \Gamma & \longrightarrow & \SL_4( \RR) \\
    \gamma & \longmapsto & \hol( \gamma) = \left(
      \begin{array}{cc}
        \chi( \gamma)^{-1} \rho_Q( \gamma) & 0 \\
        * & \chi( \gamma) \rho_L( \gamma) 
      \end{array} \right),
  \end{eqnarray*}
  with $\chi: \Gamma \rightarrow \RR_{>0}$ a character 
and $\rho_L$, $\rho_Q: \Gamma \to \SL_2( \RR)$ Fuchsian. 
For all $\gamma$ in $\Gamma -\{1\}$
  \begin{eqnarray*}
    | \lambda_+( \gamma) | &=& \chi(\gamma)^{-1} \exp( \ell_{\rho_Q}(
    \gamma)/2) \\
    | \lambda_0( \gamma) | &=& | \lambda_-( \gamma^{-1})^{-1} |=
    \chi(\gamma) \exp( \ell_{\rho_L}( \gamma)/2) \\
    | \lambda_-( \gamma) | &=& | \lambda_0( \gamma^{-1})^{-1} |=
    \chi(\gamma) \exp( -\ell_{\rho_L}( \gamma)/2) \\
    | \lambda_+( \gamma^{-1})^{-1} | &=& \chi(\gamma)^{-1} \exp(
    -\ell_{\rho_Q}( \gamma)/2).
  \end{eqnarray*}
  From the inequalities:
  \begin{equation*}
    | \lambda_+( \gamma) | \geq | \lambda_0( \gamma) | > | \lambda_-(
    \gamma) | \geq | \lambda_+( \gamma^{-1})^{-1} |,
  \end{equation*}
  we obtain $\ell_{ \rho_L}( \gamma) \leq \ell_{ \rho_Q}( \gamma)$ hence
  $\rho_L$ is conjugate to $\rho_Q$ by Fact~\ref{fact_lengthteich}. So
  $\ell_{ \rho_L}( \gamma) = \ell_{ \rho_Q}( \gamma)$, this implies that
  $\chi^2( \gamma) \geq 1$ for all $\gamma$ and $\chi$ must be the trivial
  character.  In particular the holonomy homomorphism $\hol$ and its
  semisimplification $\hol_0$ satisfy the conditions of
  Lemma~\ref{lem_notpropconv}. Thus the convex sets $C_t$ must all be sectors,
  which is a contradiction.
  
  %

The other case $\hol(\tau) = \left( \begin{array}{cc} \textup{Id} & 0
    \\ 0 & -\textup{Id} \end{array}\right)$ cannot happen. Using the
    same argument as above, we can show that $\widehat{ \hol}$ has to
    be of the form
\begin{eqnarray*}
  \widehat{ \hol} : \Gamma \times \ZZ/2\ZZ & \longrightarrow & \SL_4(
  \RR) \\
  \gamma & \longmapsto & \left(
    \begin{array}{cc}
      \rho( \gamma) & 0 \\ 0 & \rho( \gamma)
    \end{array}\right)\\
  -1 & \longmapsto  & \left(
      \begin{array}{cc}
        \textup{Id} & 0 \\ 0 & -\textup{Id}
      \end{array}\right).
\end{eqnarray*}
The developing map $\dev$ factors through the quotient $\widehat{ M} = \widetilde{M}/
\langle 2 \tau \rangle$ which can be identified with the set of triples $(
\hat{t}_+, \hat{t}_0, \hat{t}_-)$ in $\widehat{ \partial \Gamma} = \widetilde{\partial \Gamma}/ \langle 2\tau \rangle $ of
pairwise distinct elements such that $( -\hat{t}_-, \hat{t}_+,
\hat{t}_0, \hat{t}_-)$ is oriented. Following the argument and the
notation of Section~\ref{sec_geomdesc}, one can write the developing
map as 
\begin{equation*}
  \dev(\hat{t}_+, \hat{t}_0, \hat{t}_-) = [ \eta_+( \hat{t}_+) + \vfi(
\hat{t}_+, \hat{t}_0, \hat{t}_-) \eta_-( \hat{t}_-)],
\end{equation*}
where $\vfi: \widehat{ M} \rightarrow \RR$ is continuous, non-zero and
satisfies $\vfi(-\hat{t}_+, -\hat{t}_0, -\hat{t}_-) = - \vfi(
\hat{t}_+, \hat{t}_0, \hat{t}_-)$. This is a contradiction.
\end{proof}

\begin{lem}\label{lem_int_xi3}\label{lem_xi1_xi3_fix}\label{lem_interxi3_fix}
Let $\gamma \in \overline{\Gamma} -\{ 1\}$ be an element of zero translation.
with fixed points $(\tilde{t}_{+,\gamma}, \tilde{t}_{-,\gamma}) \in 
\widetilde{ \partial \Gamma}^{(2)}_{[0]}$ 
Then
\begin{itemize}
\item $\xi^3(\tilde{t}_{+,\gamma}) \cap \xi^3(\tilde{t}_{-,\gamma})$ contains
  $x_-(\gamma)$ and $x_-(\gamma^{-1})$.
\item $\xi^1(\tilde{t}_{+,\gamma}) \oplus \xi^3(\tilde{t}_{-,\gamma}) = \RR^4
  =\xi^3(\tilde{t}_{+,\gamma}) \oplus \xi^1(\tilde{t}_{-,\gamma}).$
\item $\xi^3(\tilde{t}_{+,\gamma}) \cap \xi^3(\tilde{t}_{-,\gamma})  =
  \overline{x_0(\gamma) x_-(\gamma)}.$
\end{itemize}
Moreover $x_0(\gamma) = x_-(\gamma^{-1})$.
\end{lem}

\begin{proof}
Assume that $x_-(\gamma) \notin L:= \xi^3(\tilde{t}_{+,\gamma}) \cap
 \xi^3(\tilde{t}_{-,\gamma})$. Then $L = \overline{x_+(\gamma)
   x_0(\gamma)}$. Let $t \in \widetilde{\partial\Gamma}$
 be close to $\tilde{t}_{-,\gamma}$ and consider the line $L_t = \xi^3(t)\cap \xi^3(\tilde{t}_{+,\gamma})$. 
Then  $\lim_{n\to \infty}\hol(\gamma^{-n})L_t = L$. On the other hand the point $q=L_t \cap\overline{x_-(\gamma) x_+(\gamma)}$ satisfies 
$\lim_{n\to \infty}\hol(\gamma^{-n})q \in L$ if and only if $q =
x_+(\gamma)$, and similarly we get $p=L_t \cap \overline{x_-(\gamma)
  x_0(\gamma)} = x_0(\gamma)$. 
Thus $L_t = \overline{x_+(\gamma) x_0(\gamma)} = L$ for every $t$ close enough to $\tilde{t}_{-,\gamma}$. But then Lemma~\ref{lem_xilocline} implies 
 $\bigcap_{t\in \widetilde{\partial\Gamma}} \xi^3(t) = L$, which
contradicts Lemma~\ref{lem_interxi3notline}. This concludes the first claim.

Assume that $\xi^1(\tilde{t}_{+,\gamma}) \subset
\xi^3(\tilde{t}_{-,\gamma})$, then 
$\xi^3(\tilde{t}_{+,\gamma}) \cap \xi^3(\tilde{t}_{-,\gamma})  =
\overline{x_+(
\gamma) x_-(\gamma)}$. Let $t\in \widetilde{\partial\Gamma}$
 be close to $\tilde{t}_{-,\gamma}$ and consider the line $L_t =
 \xi^3(t)\cap \xi^3(\tilde{t}_{+,\gamma})$. Then  $\lim_{n\to
   \infty}\hol(\gamma^{-n})L_t = \overline{x_+(
\gamma) x_-(\gamma)}$. But the 
intersection $p= L_t \cap \overline{x_+(
\gamma) x_0(\gamma)}$ converges to $x_+(\gamma)$ only if 
$p =x_+(\gamma)$. Hence Lemma~\ref{lem_xilocline}, Lemma~\ref{lem_interxinotdev} and Lemma~\ref{lem_xilocpoint} imply 
that
$\bigcap_{t \in \widetilde{ \partial \Gamma}} \xi^3( t) = x_+(
\gamma) = \xi¹(\tilde{t}_{+,\gamma})$. Since $\overline{\Gamma}$ acts minimally on $\widetilde{ \partial
  \Gamma}$, $\xi^1$ is constant. This is a contradiction. 

From the above,
the equality 
\begin{equation*}
  \xi^3(\tilde{t}_{+,\gamma}) \cap \xi^3(\tilde{t}_{-,\gamma})  =
  \overline{x_0(\gamma) x_-(\gamma)} = \overline{x_0(\gamma^{-1}) x_-(\gamma^{-1})}
\end{equation*}
follows.
Since $\lambda_0(\gamma) >\lambda_-(\gamma)$ and $\lambda_0(\gamma^{-1})
 >\lambda_-(\gamma^{-1})$, we necessarily have $x_0(\gamma) = x_-(\gamma^{-1}).$
\end{proof}

\begin{lem}\label{lem_noint}\label{lem_alterxiconv}
For every $(t_+, t_-) \in \widetilde{\partial\Gamma}^{(2)}_{[0]}$ we have 
\bqn
\xi^3(t_+)\cap C_{t_-} = \emptyset.
\eqn
\end{lem}

\begin{proof}
  Let $\g\in \overline{\Gamma}$ be of zero translation, then
  $\xi^3(\tilde{t}_{-,\gamma}) \cap C_{\tilde{t}_{+,\gamma}}  =
  \emptyset$, since otherwise 
 $\xi^3(\tilde{t}_{-,\gamma}) \cap C_{\tilde{t}_{+,\gamma}} =
  \mathcal{D}(\tilde{t}_{+,\gamma}, \tilde{t}_{-,\gamma})$, which is spanned by
  $x_+(\gamma)$ and $x_-(\gamma)$, contradicting
  Lemma~\ref{lem_interxi3_fix}.


  Since the set $\{(t_+, t_-) \in \widetilde{\partial
    \Gamma}^{(2)}_{[0]} \mid \xi^3(t_+) \cap C_{ t_-} = \emptyset \}$
  is closed, we conclude with Lemma~\ref{lem_closedgeo}.
\end{proof}

\begin{remark}
  The order of $t_+$ and $t_-$ in
  Lemma~\ref{lem_noint} is of no
  importance.
\end{remark}

\begin{lem}
  \label{lem_interxi3notpoint}
  The intersection $\bigcap_{t \in \widetilde{ \partial \Gamma}} \xi^3( t)$ is
  not a point.
\end{lem}

\begin{proof}
Suppose that $x=\bigcap_{ t \in \widetilde{ \partial \Gamma}}
\xi^3(t)$ is a point. 
Then, for any $\gamma$, $x = x_-(\gamma)$ or  $x= x_-(\gamma^{-1})$.

In particular we have  
\bqn
x \in \mathcal{D}(\tilde{t}_{+,\gamma}, \tilde{t}_{-,\gamma}) \cup 
\mathcal{D}( \tilde{t}_{+, \gamma^{-1}}, \tilde{t}_{-, \gamma^{-1}}) 
= \mathcal{D}(\tilde{t}_{+,\gamma}, \tilde{t}_{-,\gamma})
 \cup \mathcal{D}(\tilde{t}_{-,\gamma}, \tau^{-1} \tilde{t}_{+,\gamma}).
\eqn

Thus the set 
\bqn
C= \{ (t_+, t_-) \in \widetilde{ \partial
  \Gamma}^{(2)}_{[0]} \mid x\in \mathcal{D}(t_{+}, t_{-})
 \cup \mathcal{D}(t_{-}, \tau^{-1} t_{+})\}
\eqn
is closed and contains
$(\tilde{t}_{+,\gamma}, \tilde{t}_{-,\gamma})$ for every $\gamma \in
 \overline{\Gamma}-\{1\}$ of zero translation, hence, by
 Lemma~\ref{lem_closedgeo}, $C = \widetilde{ \partial
  \Gamma}^{(2)}_{[0]}$.

Consider  $\gamma \in \overline{\Gamma}-\{1\}$ an element of zero
 translation with $(\tilde{t}_{+,\gamma},
 \tilde{t}_{-,\gamma})$. Then, by Lemma~\ref{lem_interxi3_fix} 
\bqn
\mathcal{D}(\tilde{t}_{+,\gamma}, \tilde{t}_{-,\gamma}) \cap 
\mathcal{D}(\tilde{t}_{+, \gamma^{-1}}, \tilde{t}_{-, \gamma^{-1}})  = \emptyset 
\eqn
and assume that $x \in \mathcal{D}(\tilde{t}_{+,\gamma}, \tilde{t}_{-,\gamma}) $, so $x \notin
 \mathcal{D}(\tilde{t}_{+,
 \gamma^{-1}}, \tilde{t}_{-, \gamma^{-1}})$. 
 Then by
 continuity of $\mathcal{D}$ there exists a neighborhood $U \subset
 \widetilde{\partial \Gamma}$ of $\tilde{t}_{+,\gamma^{-1}}$ such that $x
 \notin \mathcal{D}(t, \tilde{t}_{-,\gamma^{-1}}) = \mathcal{D}(t, \tau^{-1} \tilde{t}_{+, \gamma})$ for all
 $t\in U$.  
Thus for all $t\in U$ we have that $x \in \mathcal{D}(\tilde{t}_{+,\gamma}, t)$, hence  
$x = \bigcap_{t\in U} \mathcal{D}(\tilde{t}_{+,\gamma}, t) = \xi^1(\tilde{t}_{+,\gamma})$ by
 local injectivity of $\mathcal{D}$. This again implies that $\xi^1$ is constant and gives 
a contradiction.
%
%
\end{proof}

\begin{proof}[Proof of Proposition~\ref{prop_holtautriv}]
  The first statement has been proved by elimination. Now as we already
  observed, $\xi^3$ is $\tau$-invariant. This means that for all $t$,
  $\xi^3(t)$, as a line in the dual space of $\RR^4$, is invariant by $\hol(
  \tau)$. Now an element of $\PGL( V)$ having a continuous family of invariant
  lines generating the vector space $V$ is necessarily trivial, so $\hol( \tau)
  = \textup{Id}$.
\end{proof}

The main consequence of Proposition~\ref{prop_holtautriv} is the following
\begin{cor}\label{cor:hol_factors}
  Let $(\dev, \hol)$ be the developing pair defining a properly convex
  foliated projective structure on $M$, then the holonomy factors through a
  representation
\begin{equation*}
  \hol: \overline{\Gamma} \to  \Gamma \longrightarrow \PSL_4(\RR) \subset
  \PGL_4(\RR).
\end{equation*}
\end{cor}
We will consider now the holonomy as a homomorphism
\begin{equation*}
  \hol : \Gamma \longrightarrow \PSL_4( \RR)
\end{equation*}
and the developing map
\begin{equation*}
  \dev : \overline{M} = S \widetilde{ \Sigma} \longrightarrow \PTR.
\end{equation*}

Moreover, we get that the maps $\xi^3$, $\xi^1$ and $\mathcal{D}$ also factor
through the corresponding quotient, which is $\partial \Gamma$ and $\partial
\Gamma^{(2)}$ respectively.

\subsubsection{Some Other Consequences}

\begin{lem}
  {\bf Case (T)} never occurs, that is $|\lambda_+( \gamma)| > |\lambda_0(
  \gamma)|$ for all $\gamma$ in $\Gamma -\{1\}$
\end{lem}

\begin{proof}
  If $\lambda_+ = \lambda_0$ for some $\gamma\in \Gamma$ with fixed points
  $(t_+, t_-)$, then the line $L$ spanned by $x_+( \gamma)= \xi^1(t_+)$ and
  $x_0( \gamma)$ is pointwise fixed by $\hol( \gamma)$. Here, by
  Lemma~\ref{lem_int_xi3}, 
$x_0( \gamma)$ can be chosen to lie in the
  projective plane $\xi^3( t_-)$.
  
  For $t$ in $\partial \Gamma - \{t_+\}$, the sequence $( \gamma^{-n} t)_{ n
    \in \NN}$ converges to $t_-$, hence
  \begin{equation*}
    \lim_{n \to +\infty} \hol( \gamma)^{-n} \xi^3(t) = \xi^3(t_-).
  \end{equation*}
  In particular, applying negative powers of $\hol( \gamma)$ to the
  intersection $L \cap \xi^3(t)$, this point must converge to $x_0( \gamma)$.
  This implies that $x_0 =L \cap \xi^3(t)$ and in particular $x_0 \in \bigcap_{t \in
    \widetilde{ \partial \Gamma}} \xi^3( t)$ contradicting
  Proposition~\ref{prop_holtautriv}.
\end{proof}

\begin{lem}\label{lem_xi1_xi3}
Let $(t_+,t_-) \in \partial\Gamma^{(2)}$. Assume that $t_- = t_{-,\gamma}$ for some $\gamma \in \Gamma-\{1\}$. Then 
$\xi^1(t_+) \oplus \xi^3(t_-) = \RR^4 = \xi^1(t_-) \oplus \xi^3(t_+)$.
\end{lem}
\begin{proof}
Assume that $\xi^1(t_+) \subset \xi^3(t_{-,\gamma})$. Then 
\bqn
\xi^1(t_{+,\gamma})=  \lim_{n\to \infty} \hol(\g^n) \xi^1(t_+) \subset \xi^3(t_{-,\gamma}),
\eqn
which contradicts Lemma~\ref{lem_xi1_xi3_fix}. 
The other statement follows by a similar argument.
\end{proof}

\subsection{Convexity of 
$\xi^3$}\label{sec_convexityxi3}

In this section we will show that the curve 
\bqn
\xi^3: \partial \Gamma \longrightarrow \PP^3(\RR)^*
\eqn 
is convex. 
For this we will define for every $t \in \partial \Gamma$ an auxiliary convex set  $D_t \subset \xi^3(t)$ containing $C_t$.

Define for all $t\in \partial \Gamma$
\bqn
D_t:= \xi^3(t) - \bigcup_{t'\neq t} \xi^3(t) \cap \xi^3(t').
\eqn

Then, by Lemma~\ref{lem_convexset}, $D_t$ is a convex subset in the projective plane $\xi^3(t)$, which, 
by Proposition~\ref{prop_holtautriv}, is properly convex. By Lemma~\ref{lem_noint}, $D_t$ 
contains the properly convex set $C_t$. 
In particular $D_t$ has nonempty interior in $\xi^3(t)$
and by properties of convex sets we have 
\bqn
\overset{\circ}{D_t} \subset D_t \subset \overline{D_t}.
\eqn

\begin{remark}
Of course the two convex sets $C_t$ and $\overset{\circ}D_t$ should be
equal (see Section~\ref{sec_notedevmap}). But for the moment we work 
with $D_t$ because of the following semi-continuity property, which the reader might compare with Lemma~\ref{lem_semcontCt}.
\end{remark}

\begin{lem}\label{lem_semicontDt}
Suppose that $(t_n)_{n\in \NN}\subset \partial\Gamma$ is a sequence converging to $t\in \partial \Gamma$ such that the sequence of convex sets $(\overline{D_{t_n}})_{n\in \NN}$ converges in the Hausdorff topology for compact subsets of $\PP^3(\RR)$. Then 
\bqn
\lim_{n\to \infty} \overline{D_{t_n}} \subset \overline{D_t}.
\eqn
\end{lem}
\begin{proof}
Let $D_\infty := \lim_{n\to \infty} \overline{D_{t_n}}$. Then
$D_\infty$ is a convex set in $\xi^3(t)$ containing $C_t$
(Lemma~\ref{lem_semcontCt}). If $D_\infty  - \overline{D_t} \neq
\emptyset$ then, by properties of convex sets, 
$D_\infty  - \overline{D_t}$ contains an open set $U$. 

By definition of $D_t$ there exists then $t' \in \partial \Gamma$, $t'\neq t$ such that 
\bqn
\xi^3(t') \cap U \neq \emptyset.
\eqn
Fix some $t_0 \in \partial \Gamma$, $t_0 \neq t, t'$ such that for $n$ big enough 
$D_{t_n}$ is a convex set in the affine chart 
$\PP^3(\RR) -\xi^3(t_0)$. 
Choose coordinates $(x,y,z)$ in this affine chart 
such that $\xi^3(t) = \{(x,y,z) \mid x = 0\}$ and $\xi^3(t') = \{(x,y,z) \mid z = 0\}$.

Since $U$ intersects $\xi^3(t')$ nontrivially there are points $p,q,r \in U\subset D_\infty$ such that 
\bqn
p = (x_p, y_p, z_p), \,\, q=(x_q, y_q, z_q) \text{ and } r \notin
\overline{pq},
\eqn 
with $z_p >0$ and $z_q<0$. 
By hypothesis there exists points $p_n, q_n, r_n \in \overline{D_{t_n}}$ 
such that $\lim_{n\to \infty} p_n = p$, $\lim_{n\to \infty} q_n = q$ and  
$\lim_{n\to \infty} r_n = r$. Hence for $n$ big enough $p_n, q_n, r_n \in \PP^3(\RR) -\xi^3(t_0)$ and 
\bqn
p_n = (x_{p_n}, y_{p_n}, z_{p_n}) \, \text{ and } \, q_n=(x_{q_n}, y_{q_n}, z_{q_n})
\eqn 
with $ z_{p_n}> 0$ and $z_{q_n}< 0$. In particular, the open triangle with vertices $p_n, q_n, r_n$ intersects the hyperplane $\xi^3(t') = \{(x,y,z) \mid z = 0\}$ nontrivially. 
But this is a contradiction since by convexity this open triangle is contained in $D_{t_n}$ and 
$D_{t_n} \cap \xi^3(t') = \emptyset$ for all $t_n \neq t'$.
\end{proof}

Using Lemma~\ref{lem_semicontDt} we will now be able to show that $\xi^3(t) \cap \xi^3(t')$ is tangent to the convex set $D_t$, and this will suffice to show that the curve $\xi^3: \partial \Gamma \to \PP^3(\RR)^*$ is convex. 

\begin{lem}\label{lem_intinDt}
For all $(t_+, t_-) \in \partial \Gamma^{(2)}$, the point 
\bqn
\xi^3(t_+) \cap \xi^3(t_-) \cap \mathcal{D}(t_+, t_-) 
\eqn
is contained in $\overline{D_{t_+}}$.
\end{lem}
\begin{proof}
Since $\mathcal{D}(t_+, t_-)\cap C_{t_+} = \dev((t_+, t_-))$  
Lemma~\ref{lem_noint} implies that the intersection $ \xi^3(t_-) \cap \mathcal{D}(t_+, t_-) = \xi^3(t_+) \cap \xi^3(t_-) \cap \mathcal{D}(t_+, t_-) $ is a point. 

Note that it follows from Section~\ref{sec_actiononconv} and Lemma~\ref{lem_int_xi3} that 
the statement is true for the pairs of fixed points 
$(t_{+,\gamma}, t_{-,\gamma})$ for any $\gamma \in \Gamma-\{1\}$. 
So, by Lemma~\ref{lem_closedgeo} we only have to show that 
\bqn
C=\{(t_+, t_-) \mid \xi^3(t_+) \cap \xi^3(t_-) \cap \mathcal{D}(t_+, t_-) \in \overline{D_{t_+}}\}
\eqn 
is closed.

Let $\left((t_{+,n}, t_{-,n})\right)_{n\in \NN} \subset C$ and $(t_+, t_-) = \lim_{n\to \infty} (t_{+,n}, t_{-,n})$. 
Then Lemma~\ref{lem_semicontDt} implies that 
\bqn
\xi^3(t_+) \cap \xi^3(t_-) \cap \mathcal{D}(t_+, t_-) &=&\\ 
\lim_{n\to \infty} \xi^3(t_{+,n}) \cap \xi^3(t_{-,n}) \cap \mathcal{D}(t_{+,n}, t_{-,n}) &\in& \lim_{n\to \infty} \overline{D_{t_{+,n}}} 
\subset \overline{D_{t_+}}.
\eqn
Hence $(t_+, t_-) \in C$ and $C$ is closed.
\end{proof}

The intersection $\xi^3(t_+) \cap \xi^3(t_-) \cap \mathcal{D}(t_+, t_-)$ is one of the points of 
intersection 
\bqn
\mathcal{D}(t_+, t_-)\cap \partial {D_{t_+}} = \{e_+(t_+, t_-), e_-(t_+, t_-)\}, 
\eqn
where $e_-(t_+, t_-)$ is the point such that the open segment $]\xi^1(t_+),e_-(t_+, t_-)[$ 
in $\overline{D_{t_+}}$ intersects $C_{t_+}$.

For all $t_+ \in \partial \Gamma$ the map 
\bqn
e_- : \partial\Gamma-\{t_+\}& \to & \xi^3(t_+)\\
t_- &\mapsto& e_-(t_+, t_-)
\eqn
is continuous.

\begin{lem}\label{lem_inteminus}
For all $(t_+, t_-) \in \partial \Gamma^{(2)}$ we have 
\bqn
\xi^3(t_+) \cap \xi^3(t_-) \cap \mathcal{D}(t_+, t_-) = e_-(t_+, t_-). 
\eqn
\end{lem} 

\begin{proof}
We first prove that, for all $(t_+, t_-) \in \partial \Gamma^{(2)}$, we have 
\bqn
\xi^3(t_+) \cap \xi^3(t_-) \cap \mathcal{D}(t_+, t_-) \neq \xi^1(t_+).
\eqn

Assume the contrary. Then there exists $ (t_+, t_-) \in \partial \Gamma^{(2)}$ such that 
%
$\xi^1(t_+) \in \xi^3(t_+) \cap \xi^3(t_-)$ and hence $\xi^1(t_+) \notin D_{t_+}$. But since 
$\xi^1(t_+) \in \overline{D_{t_+}}$ (Lemma~\ref{lem_intinDt}) we have that 
$\xi^1(t_+) = e_+(t_+, t_-)$. 

By  Lemma~\ref{lem_xi1_xi3} we have that for all
$t_{-, \gamma}$, $\gamma \in \Gamma-\{1\}$ 
\bqn
\xi^3(t_+) \cap \xi^3(t_{-, \gamma}) \cap \mathcal{D}(t_+, t_{-, \gamma}) = 
e_-(t_+, t_{-,\gamma}) 
\neq \xi^1(t_+).
\eqn

Let $(t_{-,\gamma_n})_{n\in \NN} \subset \partial \Gamma$ be a sequence 
of fixed points of elements $\gamma_n \in \Gamma-\{1\}$ with 
$\lim_{n\to \infty} t_{-, \gamma_n} = t_-$.
Then, by continuity of $\xi^3$ and $\mathcal{D}$ the point of intersection 
\bqn
\xi^3(t_+) \cap \xi^3(t_{-, \gamma_n}) \cap \mathcal{D}(t_+, t_{-, \gamma_n})
\eqn 
converges to $\xi^1(t_+)$. But on the other hand
 by continuity of $e_-(t_+, \cdot)$ 
\bqn
\xi^3(t_+) \cap \xi^3(t_{-, \gamma_n}) \cap \mathcal{D}(t_+, t_{-, \gamma_n}) 
= e_-(t_+, t_{-,\gamma_n})
\eqn 
should converge to $e_-(t_+, t_-)$.
This is a contradiction. 

Note that the set 
\bqn
C=\{ (t_+, t_-) \mid \xi^3(t_+) \cap \xi^3(t_-) \cap \mathcal{D}(t_+, t_-) = e_-(t_+, t_-)\}
\eqn 
contains the set of pairs $(t_{+, \gamma}, t_{-,\gamma})$ of fixed points for any $\gamma \in \Gamma-\{1\}$.
So by Lemma~\ref{lem_closedgeo} we only have to show that $C$ is closed. 
Let $\left((t_{+,n}, t_{-,n})\right)_{n\in \NN} \subset C$ be a sequence converging to $(t_+, t_-)$. 
Then Lemma~\ref{lem_semicontDt} and the definition of $e_-$ imply that 
\begin{eqnarray*}
\xi^3(t_+) \cap \xi^3(t_-) \cap \mathcal{D}(t_+,t_-) &=& \lim_{n\to \infty} 
\xi^3(t_{+,n}) \cap \xi^3(t_{-,n}) \cap \mathcal{D}(t_{+,n},t_{-,n})\\
=\lim_{n\to \infty} e_-(t_{+,n}, t_{-,n}) &\in& [\xi^1(t_+),  e_-(t_+, t_-)] 
\subset \overline{D_{t_+}}.
\end{eqnarray*}
By the above $\xi^3(t_+) \cap \xi^3(t_-) \cap \mathcal{D}(t_+,t_-) \neq
 \xi^1(t_+)$, hence, since $\xi^3(t_+) \cap \xi^3(t_-) \cap
 \mathcal{D}(t_+,t_-)$ is not in $D_{t_+}$, we necessarily have 
$\xi^3(t_+) \cap \xi^3(t_-) \cap \mathcal{D}(t_+,t_-) = e_-(t_+, t_-)$
 and $C$ is closed.
\end{proof}

\begin{lem}\label{lem_geoinjective}
For all $t_+ \in \partial \Gamma$ the map 
\bqn
f_{t_+} : \partial\Gamma -\{ t_+\}  &\to& \PP(\xi^3(t_+)/\xi^1(t_+))\\
t_- &\to& \mathcal{D}(t_+, t_-)
\eqn
is injective. 
\end{lem}
\begin{proof}
By definition of $\xi^1$ we have $\xi^1 (t_+) \subset \mathcal{D}(t_+, t_-)$ for all $t_- \neq t_+$, so $f_{t_+} $ is well defined. 
Since $C_{t_+}$ is a convex set, $f_{t_+}$ is not surjective and we can think of $f_{t_+}$ as a map from $\RR$ to $\RR$. 

Assume that  $f_{t_+}$ is not injective. Then there exist a point $t_- \in \partial\Gamma$, 
a neighborhood $U_- \subset \partial\Gamma$ of $t_-$ and an interval 
\bqn
[l_-, l_0[ \subset \PP(\xi^3(t_+)/\xi^1(t_+))
\eqn
 such that 
\bqn
f_{t_+}(t_-) = l_- \, \, \text{ and } \,  f_{t_+} (U_-) \subset [l_-, l_0[. 
\eqn

This implies that there exist a point $(t_+, t_0, t_-)$  in the leaf 
$t_+ \subset \overline{M}$ 
and an open neighborhood $V$ of $(t_+, t_0, t_-)$ 
\bqn
\dev(V) \subset U = \bigcup_{l \in [l_-, l_0[} l \subset \xi^3(t_+), 
\eqn
but $U$ is not a neighborhood of $\dev(t_+, t_0, t_-) \in l_-$.
This contradicts the fact that $\dev $ is a local homeomorphism. 
\end{proof}

\begin{lem}\label{lem_noseginD}
For all $t_+ \in \partial \Gamma$ the map 
\bqn
e_- : \partial\Gamma-\{t_+\}& \longrightarrow & \xi^3(t_+)\\
t_- &\longmapsto& e_-(t_+, t_-)
\eqn
is injective. 
Moreover, for all $(t_+, t_-) \in \partial \Gamma^{(2)}$ there is no 
open segment in $\partial \overline{D_{t_+}}$ containing $e_-(t_+, t_-)$.
\end{lem}
\begin{proof}
Injectivity follows from the injectivity of $\mathcal{D}(t_+, \cdot)$ (Lemma~\ref{lem_geoinjective}).

Assume that $\partial \overline{D_{t_+}}$ contains a segment containing $e_-(t_+, t_-)$. 
Let $L$ be the projective line supporting this segment. 
Then for $t'$ close to $t_-$ the point 
$e_-(t_+, t') \in \partial \overline{D_{t_+}}$ will be an interior
point of this segment. Since 
$\xi^3(t_+)\cap \xi^3(t')$ is a line tangent to $D_{t_+}$ at $e_-(t_+,
t')$  this implies 
$\xi^3(t_+)\cap \xi^3(t') = L$. 
By Lemma~\ref{lem_xilocline} this implies $\bigcap_{t \in \partial \Gamma} \xi^3(t) = L$, which 
contradicts Lemma~\ref{lem_interxi3notline}.
\end{proof}
\begin{prop}\label{prop:curve_convex}
  The $\hol$-equivariant curve
\begin{equation*}
  \xi^3: \partial \Gamma \longrightarrow \PTRd
\end{equation*}
is convex, that is 
\begin{center}
  for all $(t_1, t_2, t_3,t_4)$ in $\partial \Gamma^4$ pairwise distinct
  $\bigcap_{ i=1}^{4} \xi^3( t_i) = \emptyset$.
\end{center}
\end{prop}
\begin{proof}
  Let us rename the four points $t_1, t_2, t_3,t_4 \in \partial \Gamma$ by $
  t_{-,1},t_{-,2},t_{-,3}, t_+$ such that $(t_+, t_{-,3},t_{-,2},t_{-,1})$ is
  positively oriented. By Lemma~\ref{lem_noseginD} the three lines
  $\xi^3(t_{-,i}) \cap\xi^3(t_+)$, $i=1,2,3$ are tangent to the convex $D_{ t_+}$ at the
  three distinct points $e_-(t_+, t_{-,1}) $, $e_-(t_+, t_{-,2})$ and $e_-(t_+, t_{-,3})$. By Lemma~\ref{lem_noseginD} they cannot intersect. In
  particular, the intersection
  \begin{equation*}
   \bigcap_{ i=1}^{4} \xi^3( t_i) =  \big( \xi^3( t_{-,1}) \cap
    \xi^3(t_+)\big) \cap \big( \xi^3( t_{-,2})
    \cap \xi^3(t_+)\big) \cap \big( \xi^3( t_{-,3}) \cap \xi^3(t_+)\big)
  \end{equation*}
  is empty.
\end{proof}
This shows that the holonomy representation 
$\rho=\hol$ of a properly convex foliated projective structure factors
through a convex representation and hence lies in the
Hitchin component in view of Theorem~\ref{thm_hiteqconvex}. 
In particular Theorem~\ref{thm_convexFrenet} states 
that there exists a unique $\hol$-equivariant curve 
\bqn
\xi= (\xi^1, \xi^2, \xi^3) : \partial \Gamma \longrightarrow \Flag(\RR^4).
\eqn
In the following section we want to describe a specific construction of the missing  
\bqn
\xi^2: \partial \Gamma \longrightarrow \Gr_{2}^{4}.
\eqn
\smallskip
The following figure shows the configuration of the images by the
developing map of two (weakly) stable leaves and serves as
motivation for the geometric construction of $\xi^2$ which is already pictured.

\begin{figure}[htbp]
    \centering \input{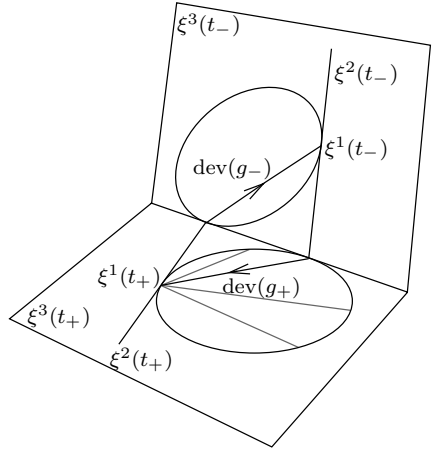}
    \caption{Two (weakly) stable leaves}
    \label{fig_twoleaves}
  \end{figure}

In Section~\ref{sec_notedevmap} we will compare the properly
convex foliated structure 
$(\dev, \hol)$ and the properly convex foliated structure associated
to $\rho$ in Section~\ref{sec_reptostruct}, this will complete the
proof of Theorem~\ref{thm:struct_hol}.

\subsection{Definition of $
\xi^2$}

The following is a consequence of Lemma~\ref{lem_inteminus}, which also
follows from the Frenet property of $\xi=(\xi¹, \xi², \xi³)$:
for all $(t_+, t_-) \in \partial\Gamma^{(2)}$ we have 
\bqn
\xi^1(t_+) \oplus \xi^3(t_-) = \RR^4 = \xi^1(t_-) \oplus \xi^3(t_+).
\eqn

This implies that the continuous $\hol$-equivariant map:
\begin{eqnarray*}
  \mathcal{P}: \partial \Gamma^{(2)} & \longrightarrow & \Gr_{2}^{4}
  (\RR)\\
  (t_+, t_-) & \longmapsto & \big( \mathcal{D}( t_+, t_-) \cap \xi^3(
  t_-) \big) \oplus \xi^1( t_-)
\end{eqnarray*}
is well defined.
\begin{lem}
  \label{lem_Ptang}
  For all $(t_+, t_-) \in \partial \Gamma^{(2)}$
  \begin{equation*}
    \mathcal{P}( t_+, t_-) \cap C_{t_-} = \emptyset.
  \end{equation*}
  For all $t_-$ in $\partial \Gamma$, the function
  \begin{eqnarray*}
    \mathcal{P}_{t_-}: \partial \Gamma - \{ t_-\} & \longrightarrow & \Gr_{2}^{4} (\RR)
    \\
    t_+ & \longmapsto & \mathcal{P}( t_+, t_-)
  \end{eqnarray*}
  is constant.

\end{lem}

\begin{proof}
  The subset
  \begin{equation*}
    \{ (t_+, t_-) \in \partial \Gamma^{(2)} \, |\, \mathcal{P}( t_+,
    t_-) \cap C_{t_-} =\emptyset\}
  \end{equation*}
  is closed. Since the pairs consisting of fixed points $(t_{+, \gamma}, t_{-,\gamma})$
  of a nontrivial element $\gamma \in \Gamma$ are contained in this set (see Figure~\ref{fig_twoleaves} and Lemma~\ref{lem_int_xi3}),
  it equals $\partial\Gamma^{(2)}$.

  The set
  \begin{equation*}
    \{t_- \in \partial \Gamma \, |\,\mathcal{P}_{t_-} \text{ is constant
    }\}
  \end{equation*}
  is $\Gamma$-invariant and closed by the continuity of $\mathcal{P}$.  This
  set does contain points $t_{-,\gamma}$ which  are the fix points of a
  non-trivial element $\gamma \in \Gamma$, since in that case the projective
  line $\mathcal{P}(t, t_-)$ is necessarily the unique tangent line to
  $C_{t_-}$ at $\xi^1( t_-)$. Hence, this set equals $\partial \Gamma$.
\end{proof}

\begin{defi}
  Define the map
 \begin{equation*}
   \xi^2 : \partial \Gamma \longrightarrow \Gr_{2}^{4}( \RR),
 \end{equation*}
 by setting $\xi^2(t_-)$ to be the value of the constant function
 $\mathcal{P}_{t_-}$.
\end{defi}

The map $\xi^2:\partial \Gamma \longrightarrow \Gr_{2}^{4}( \RR)$ is a 
continuous, $\hol$-equivariant injective map. 

\subsection{Proof of Theorem~\ref{thm:main}}
\label{sec_notedevmap}
In Section~\ref{sec_convexityxi3} we showed that the image of the holonomy map
\begin{equation*}
  \hol: \mathcal{P}_{pcf}(M) \longrightarrow \rep(\overline{\G}, \PGL_4(\RR))
\end{equation*}
is contained in the Hitchin component $\mathcal{T}^4(\Gamma)
\subset\rep(\overline{\G}, \PGL_4(\RR))$.  In Section~\ref{sec_reptostruct} we
defined a section
\begin{equation*}
  s : \mathcal{T}^4(\Gamma) \longrightarrow  \mathcal{P}_{pcf}(M)
\end{equation*}
by constructing a developing pair $(\dev_\xi, \hol)$ of a properly convex
foliated projective structure on $M$ starting with a
$\hol$-equivariant Frenet curve $\xi = (\xi^1, \xi^2, \xi^3)$.

To finish the proof of Theorem~\ref{thm:main} we have to show that the maps
$\hol$ and $s$ are inverse to each other. That $\hol\circ s = \id$ is
immediate. To show that $s \circ \hol = \id$ we have to show that the
developing pairs $(\dev, \hol)$ and $(\dev_\xi, \hol)$ are
equivalent (see Section~\ref{sec_folprostr}).

\begin{lem}\label{lem_devequaldev}
The developing map $\dev$ of a properly convex foliated projective structure is a homeomorphism onto its image
\bqn
\dev: \overline{M} \longrightarrow \bigcup_{t\in \partial\Gamma} C_t \subset \PP^3(\RR).
\eqn
Furthermore the image of $\dev$ equals the image of $\dev_\xi$. 
\end{lem}
\begin{proof}
The two convex sets $C_{t_+}$ and
$C_{s_+}$ intersect only when $t_+ = s_+$, and by
Lemma~\ref{lem_geoinjective} the two projective
lines $\mathcal{D}(t_+,t_-)$ and $\mathcal{D}(t_+,s_-)$ intersect in
the convex $C_{t_+}$ if and only if $t_-=s_-$. Therefore, 
if two points $(t_+,
t_0, t_-)$ and $(s_+,s_0, s_-)$ have the same image under $\dev$ then
necessarily $t_+ = s_+$ and $t_- = s_-$, but we already observed that the restriction
of $\dev$ to the geodesic $(t_+,t_-)$ is a homeomorphism onto its
image (see proof of Lemma~\ref{lem_notexchange}), so $(t_+, t_0, t_-) =
(s_+,s_0, s_-)$. 
Thus $\dev$ is a bijective local homeomorphism, hence a homeomorphism. 

The image $\dev(\overline{M})$ is a $\hol(\Gamma)$-invariant open
subset of $\dev_\xi(\overline{M}) = \Omega_\xi$. The group
$\hol(\Gamma)$ acts properly discontinuous with compact quotient on
both sets $\dev(\overline{M})$ and $\dev_\xi(\overline{M})$, so 
$\dev(\overline{M})/\hol(\Gamma)$ 
is a compact and open subset of $\dev_\xi(\overline{M})/
\hol(\Gamma)$, 
which is connected, hence  $\dev(\overline{M})/\hol(\Gamma) = \dev_\xi(\overline{M})/ \hol(\Gamma)$. 
\end{proof}

A direct corollary is that, for all $t \in \partial \Gamma$, we have 
\bqn
C_t = \overset{\circ}{D_t}.
\eqn

So $\dev^{-1} \circ
\dev_\xi: \overline{M} \rightarrow \overline{M}$ is a 
$\Gamma$-equivariant homeomorphism, hence it descends to a homeomorphism
$h:M \rightarrow M$. Since $\dev$ and $\dev_\xi$ send the geodesic
$(t_+, t_-)$ onto the same segment of $C_{t_+}$ (namely the segment
with endpoints $\xi^1( t_+)$ and $\xi^3(t_+) \cap \xi^2(t_-)= e_-(t_+,t_-)$, see Lemma~\ref{lem_inteminus}), the
map $h$ sends each geodesic leaf of $M$ into itself.
As noted in Remark~\ref{rem_folprostr} this is enough to
ensure that the two pairs are equivalent and define the same element in
$\mathcal{P}_{pcf}(M)$. This proves
Theorem~\ref{thm:main}.

\section{Projective Contact Structures}
\label{sec_selfdual}
\begin{quote}
  In this section we conclude the description of the Hitchin component for the
  symplectic group $\PSp_4(\RR)$ from Theorem~\ref{thm:main}.
\end{quote}

The $4$-dimensional irreducible representation
\begin{equation*}
  \rho_4: \PSL_2(\RR) \longrightarrow \PSL_4(\RR)
\end{equation*}
preserves a symplectic form $\omega$ on $\RR^4 = \textup{Sym}^3 \RR^2$, and so
takes values in $\PSp_4(\RR) = \PSp(\RR^4, \omega) \subset \PSL_4(\RR)$

The symplectic form $\omega$ on $\RR^4$ defines a canonical contact structure
on $\PTR$, that is a non integrable distribution.  Let $[l]$ be the line
spanned by a vector $l\in \RR^4$.  The tangent space $T_{[l]}\PTR$ is
naturally identified with the space of homomorphisms $\hom([l],\RR^4/[l])$.
Let $H_{l} = [l]^{\perp_\omega}$. Then the contact distribution is 
$$
\mathcal{H}_l = \hom([l], H_l/[l]) \subset \hom([l],\RR^4/[l]) \simeq
T_{[l]}\PTR. 
$$
We consider $\PTR$ with this contact structure. The maximal subgroup of
$\PSL_4(\RR)$ preserving this contact structure is $\PSp_4(\RR)$.


\begin{defi}
  A projective contact structure on $M$ is a projective structure $(U,\vfi_U)$
  on $M$, such that the coordinates changes $\vfi_V \circ \vfi_{U}^{-1}$ are locally
  in $\PSp_4(\RR)$.
\end{defi}
Two projective contact structures are equivalent if they are equivalent as
projective structure by a homeomorphism $h$ which preserves the contact
structure.

A projective contact structure on $M$ is given by a developing pair $(\dev,
\hol)$, where $\dev: \widetilde{M} \to \PTR$ is a local homeomorphism which is
equivariant with respect to the holonomy homomorphism $\hol: \overline{\G} \to
\PSp_4(\RR)$ with values in the symplectic group.  A projective contact
structure which is properly convex foliated as projective structure will be
called a properly convex foliated projective contact structure, and we denote
by $\mathcal{PC}_{pcf}(M)$ be the equivalence classes of properly convex
foliated projective contact structures.
\begin{thm}\label{thm_sympl}
The holonomy map is a homeomorphism between $\mathcal{PC}_{pcf}(M)$ and the Hitchin component 
${\mathcal T}( \Gamma, \PSp_{4}(\RR)) \subset \rep(\G,\PSp_{4}(\RR))$.
%
\end{thm}
This statement is indeed a direct Corollary of Theorem~\ref{thm:main} since
the Hitchin component ${\mathcal T}^4(\G, \PSp_4(\RR))$ is canonically
identified with the subset of representations in ${\mathcal T}^4(\G)$
preserving the symplectic form $\omega$ on $\RR^4$.

Using the existence of a $\hol$-equivariant
Frenet curve $\xi=(\xi^1,\xi^2, \xi^3): \partial \Gamma \to \Flag(\RR^4)$ for
a properly convex foliated projective structure $(\dev, \hol)$ on $M$ given in
Section~\ref{sec_proptorep} we get that properly convex foliated projective
contact structures carry indeed a much richer structure.

\begin{prop}
The contact distribution of a properly convex foliated projective
structure carries a contact vector field and a natural Hermitian structure. 
\end{prop} 
\begin{proof}
  Note that a Frenet curve equivariant with respect to a representation $\rho:
  \Gamma \to \PSp(4,\RR)$ is of the form $\xi=( \xi^1, \xi^2 =
  {\xi^2}^{\perp_\omega}, \xi^3 = {\xi^1}^{\perp_\omega})$. In particular for
  every $t\in \partial \Gamma$, the plane $\xi^2(t)$ is a Lagrangian subspace
  of $\RR^4$.  Moreover, for $(t_+, t_0, t_-)$ the three Lagrangians
  $\xi^2(t_+)$, $\xi^2(t_0)$, $ \xi^2(t_-)$ are pairwise transverse and
  define a complex structure $J_{(t_+, t_0, t_-)}$ on $\RR^4 =\xi^2(t_-)\oplus
  \xi^2(t_+)$ such that the symmetric bilinear form $\omega(\cdot, J_{(t_+,
    t_0, t_+)} \cdot)$ is positive definite (this is true in greater 
generality for maximal representations \cite{Burger_Iozzi_Labourie_Wienhard},
 here this statement can be proved by deformation starting from Fuchsian
 representations).

Let us recall this construction for the reader's convenience: 
If $(L_+, L_0, L_-)$ is a triple of pairwise transverse Lagrangians,
  then $L_0$ can be realized as the graph of $F_+ \in \hom(L_+, L_-)$
  and as graph of $F_- \in \hom(L_-, L_+)$. Then $F_+ \circ F_- =
  \textup{Id}_{ L_-}$ and $F_- \circ F_+ = \textup{Id}_{ L_+}$, and
  the assignment 
$$l_+ + l_- \in L_+ \oplus L_- \longmapsto -F_-(l_-) +
  F_+(l_+)$$ 
defines a complex structure $J_{L+, L_0, L_-}$. We write
  $J_{(t_+,t_-,t_0)}$ for $J_{ \xi^2(t_+),
    \xi^2(t_0),\xi^2(t_-)}$. 
We can use the complex structure $J$ to define a 
$\hol$-equivariant vector field: $V_J(t_+, t_0, t_-) \in T_{\dev(t_+,
  t_0, t_-)} \PTR$ is the linear map $ l \mapsto J_{ (t_+,
  t_0,t_-)} \cdot l $. The definiteness of the quadratic
  form implies that $V_J$ is
orthogonal to the distribution $\mathcal{H}_l$ 
and descends to a contact vector field on $M$. Moreover $J_{(t_+, t_0,
 t_-)}$ induces a Hermitian structure on $\mathcal{H}_l$.
\end{proof}

\appendix

\section{Some Useful Facts}
\label{sec_appendix}

Let $\G=\pi_1(\Sigma)$ be the fundamental group of a closed oriented
surface of genus $g\geq 2$. We recall some facts
about representations of $\Gamma$ and its central extension 
$\overline{ \Gamma}$ into $\PSL_2(\RR)$ and $\PSL_2(\CC)$.

\subsection{Equivariant curves}

Fuchsian representations are the representations coming from a uniformization
of the surface 
$\Sigma$. They are precisely the faithful and discrete representations:
\begin{lem}\label{lem_fuchsian}
  If $\iota : \Gamma \rightarrow \PSL_2( \RR)$ is faithful and discrete, then
  it is Fuchsian. This will be for example the case if for all $\gamma$ in $\Gamma-\{1\}$,  
$\rho(\gamma)$ is a nontrivial hyperbolic element of $\PSL_2(\RR)$. 
\end{lem}

\begin{proof}
  The hypothesis imply that the action of $\Gamma$ on $\HH^2$ is
  proper and free. The quotient surface $\iota( \Gamma) \backslash
  \HH^2$ is a surface whose fundamental group is $\Gamma = \pi_1(
  \Sigma)$ hence is diffeomorphic to $\Sigma$. So we get a
  uniformization $\Sigma \simeq \iota( \Gamma) \backslash \HH^2$ and
  its holonomy is $\iota$.

  For the second statement, the hypothesis already implies that $\rho$
  is faithful. To show that it is discret we need to show that the neutral
  component $\overline{ \rho( \Gamma)}^\circ$ is trivial. This subgroup is
  normalized by $\rho( \Gamma)$ and by $\PSL_2( \RR)$ since $\rho$ is
  Zariski dense (the proper Lie subgroups of $\PSL_2(\RR)$ are
  virually solvable, so a representation having value in one of this
  subgroup has a nontrivial kernel). Since the elliptic elements in
  $\PSL_2(\RR)$ form an open subset, $\rho( \Gamma)$ cannot be dense.
  This implies that $\overline{ \rho( \Gamma)}^\circ$ is trivial.
\end{proof}

\begin{lem}
  \label{lem_curveteich}
  Let $\rho : \Gamma \rightarrow \PSL_2( \RR)$ be a
  representation. Suppose that
  there exists a continuous non-constant, $\rho$-equivariant curve $\xi:
  \partial \Gamma \rightarrow \POR$ then $\rho$ is a discrete
  and faithful representation, hence Fuchsian.
\end{lem}

\begin{proof}
  Since $\Gamma$ is torsion-free both faithfulness and discreteness will 
  follow from the property:
  \begin{quote}
    there is no sequence $( \gamma_n)_{n \in \NN}$ in $\Gamma$ such that
    $\lim \gamma_n = \infty$ and $\lim \rho(\gamma_n) = \textup{Id}$.
  \end{quote}
  Suppose that such a sequence exists. Applying the Theorem of 
  Abels, Margulis and So{\u\i}fer (see next Lemma) there exists an element $f$ 
  in $\Gamma$ such that, up to extracting a subsequence, the sequence
  $\delta_n = \gamma_n f$ satisfies that 
  \begin{center}
    $\lim t_{+, \delta_n} = t_+$ and $\lim t_{-, \delta_n} = t_-$ exist with
    $t_+ \neq t_-$.
  \end{center}
  The sequence $(\delta_n)_{n \in \NN}$ has the following property:
  \begin{center}
    $\lim \delta_n = \infty$ and $\lim \rho( \delta_n) = \rho(f)$.
  \end{center}
  For any $t \neq t_+$ the limit of $( \delta_{n}^{-1} \cdot t)$ equals $t_-$
  so that
  \begin{equation*}
    \rho( f)^{-1} \xi(t) = \lim \rho( \delta_n)^{-1} \xi(t) = \xi(t_-).
  \end{equation*}
  Hence $\xi$ is constant on $\partial \Gamma - \{ t_+ \}$ and therefore
  constant by continuity, contradicting the hypothesis.
\end{proof}

To be complete we explain now how the statement we used follows easily from 
Abels, Margulis and So{\u\i}fer's Theorem.

\begin{lem}
  There exists a finite subset $F$ in $\Gamma$ such that for any sequence
  $(\gamma_n)_{n \in \NN}$ in $\Gamma$ converging to infinity there exists
  some $f$ in $F$ such that (up to extracting a subsequence) the sequence $\delta_n =
  \gamma_n f$ satisfies:
  \begin{center}
    $t_+ = \lim t_{+, \delta_n} \neq t_- = \lim t_{-, \delta_n}$, 
  \end{center}
 with $t_{\pm, \delta}$ being the fixed points of $\delta$ 
  in $\partial \Gamma$.

  Also for $(\delta_n)_{n \in \NN}$ one has the following property:
  \begin{center}
  for all $t \neq t_+$, $\lim \delta_{n}^{-1} \cdot t = t_-$.
  \end{center}
\end{lem}

\begin{proof}
  First we realize $\Gamma$ as a cocompact lattice in $\SL_2( \RR)$ so 
  that the boundary at infinity is identified with $\PP^1( \RR)$.

  Recall that an $\RR$-split element $A$ of $\SL_2( \RR)$ is called
  $(r, \varepsilon)$-proximal if, with $t_{\pm, A} \in \PP^1( \RR)$ being its
  eigenlines, we have
  \begin{itemize}
  \item	$\textup{d}_{ \PP^1( \RR)}( t_{+,A}, t_{-,A} ) \geq 2r$,
  \item	if $x_1, x_2$ in $\PP^1( \RR)$ are two points satisfying 
    $\textup{d}_{ \PP^1( \RR)}( x_i, t_{-,A} ) \geq \varepsilon$ then
    \begin{eqnarray*}
      \textup{d}_{ \PP^1( \RR)}( A \cdot x_i, t_{+,A} ) &\leq& \varepsilon \\
      \textup{d}_{ \PP^1( \RR)}( A \cdot x_1, A \cdot x_2 ) &\leq& 
      \varepsilon \textup{d}_{ \PP^1( \RR)}( x_1, x_2 ).
    \end{eqnarray*}
  \end{itemize}
  Here $\textup{d}_{ \PP^1( \RR)}$ is a distance on $\PP^1( \RR)$ coming 
  from a norm on $\RR^2$.
  
  The Theorem of Abels, Margulis and So{\u\i}fer (see \cite{AbelsMargulisSoifer}) 
  states that there exist 
  $r > \varepsilon > 0$ and a finite subset $F \subset \Gamma$ such that 
  \begin{quote}
    for any $\gamma$ in $\Gamma$ there exists $f$ in $F$ such that $\gamma f$ 
    is $(r, \varepsilon)$-proximal.
  \end{quote}

  In particular for $\gamma_n$ we find $f_n \in F$ such that 
  $\delta_n = \gamma_n f_n$ is $(r, \varepsilon)$-proximal hence 
  $\textup{d}_{ \PP^1( \RR)}( t_{+, \delta_n}, t_{-, \delta_n} ) \geq 2r$.
  By compacity we can extract  a subsequence such that the limits
  $t_+ = \lim t_{+, \delta_n}$ and $t_- = \lim t_{-, \delta_n}$ exist
  and such that $(f_n)$ is constant equal to $f \in F$.
  Then $\textup{d}_{ \PP^1( \RR)}( t_{+}, t_{-} ) \geq 2r$ 
  and $t_+ \neq t_-$.

  Denote by $\lambda_n > 1$ the maximum absolute value of the eigenvalues 
  of $\delta_n$. Then $\lim \lambda_n = + \infty$, because if not 
  there exists a subsequence such that  $\lim \lambda_n = \lambda$ and $(\delta_n)$ 
  converges to an element of $\SL_2( \RR)$ with eigenlines $t_\pm$ and 
  eigenvalues $ \lambda^{\pm 1}$ (or $ -\lambda^{\pm 1}$), contradicting $\lim \delta_n = \infty$.

  Also for any $t \neq t_+$, if $n$ is big enough, then the absolute 
 value of the crossratio
  \begin{equation*}
    | [t, \delta_{n}^{-1} \cdot t, t_{+, \delta_n}, t_{-, \delta_n} ] | =
    \frac{ 1}{ \lambda_n},
  \end{equation*}
  if $t \neq t_{-, \delta_n}$. Hence for any accumulation point $s$ of 
  the sequence $(\delta_{n}^{-1} \cdot t)_{ n \in \NN}$ the equality 
  \begin{equation*}
    [t,s,t_+, t_-] = 0
  \end{equation*}
  holds, which implies that  $s=t_-$ and  $\lim \delta_{n}^{-1} \cdot t = t_-$.
\end{proof}



\subsection{The Length Function}
\label{sec_lengthfunc}
Let $\rho : \Gamma \rightarrow \PSL_2( \RR)$ be a Fuchsian representation. 
The length function for $\rho$ is defined by:
\begin{center}
  for any $\gamma$ in $\Gamma - \{1\}$, $\ell_\rho (\gamma) = t$
  $\Leftrightarrow$ $\rho(\gamma)$ is conjugate to $\displaystyle 
  \left( \begin{array}{cc}
      e^{t/2} & 0 \\ 0 & e^{-t/2}
    \end{array} \right)$.
\end{center}
We recall the well known
\begin{fact}
  \label{fact_lengthteich}
  Let $\rho, \rho': \Gamma \rightarrow \PSL_2( \RR)$ be two discrete and
  faithful representation, such that for all $\g \in \Gamma$ the inequality
  $\ell_\rho( \gamma) \leq \ell_{ \rho'}( \gamma)$ holds, then $\rho$ and
  $\rho'$ are conjugate by an element of $\SL^\pm_2(\RR)$.  In particular,
  $\ell_\rho = \ell_{\rho'}$.
\end{fact}

\begin{proof}
  Denote by $K(\rho, \rho')$ the supremum
  \begin{equation*}
    K(\rho, \rho') = \sup_\gamma \log \frac{\ell_{ \rho}( \gamma)}{\ell_{ \rho'}( \gamma)}.
  \end{equation*}
  Then \cite[Theorem~3.1]{Thurston} states that, if $\rho$ and $\rho'$ are not conjugate, then $K(\rho, \rho') >0$.
\end{proof}

\subsection{Conjugate Representations}
\label{sec_conjrep}

\begin{lem}
  \label{lem_adhzarofpair}
  Let $\rho_1$, $\rho_2: \Gamma \rightarrow \PSL_2( \CC)$ be two
  representations with Zariski dense image (this will be automatically the
  case when $\rho_i$ is faithful). We have the following alternative:
  \begin{itemize}
  \item either $\rho_1$ and $\rho_2$ are conjugate (by an element of
    $\SL_2(\CC)$),
  \item or the representation $\rho=( \rho_1, \rho_2)$ of $\Gamma$ into
    $\PSL_2( \CC) \times \PSL_2( \CC)$ has Zariski-dense image.
  \end{itemize}
\end{lem}

\begin{proof}
  Let $G$ denote the Zariski closure of $\rho( \Gamma)$ and $p_1, p_2 :\PSL_2(
  \CC) \times \PSL_2( \CC) \rightarrow \PSL_2( \CC)$ the two projections. By
  hypothesis the restrictions $p_{i|G}$ are onto, in particular the two
  kernels $\ker p_{i|G}$ have the same dimension.
  
  Suppose that this dimension is not zero then $p_2( \ker p_1 \cap G)$ has
  positive dimension and is normalized by $\rho_2(\Gamma)$, hence $p_2( \ker
  p_1 \cap G) = \PSL_2(\CC)$. Since this holds also for the first projection,
  we have $G = \PSL_2( \CC) \times \PSL_2( \CC)$.
  
  If the above dimension is zero, the group $p_2( \ker p_1 \cap G)$ is finite
  and normal hence trivial. This means that $G \simeq \PSL_2(\CC)$ and
  $\rho_1, \rho_2$ are conjugate.
\end{proof}

\subsection{Representations of $\overline{\Gamma}$}
\label{sec_repgammabar}

\begin{lem}
\label{lem_repgambarsl2}
  Let $\rho: \overline{ \Gamma} \to \PSL_2(\RR)$ a representation, then 
  $\rho(\tau)$ cannot be nontrivial $\RR$-split.
\end{lem}

\begin{proof}
  Up to conjugation, we can suppose $\rho( \tau) = 
  \left( \begin{array}{cc} \lambda & 0 \\ 
  0 & \lambda^{-1} \end{array} \right)$ with $\lambda >1$ then 
  $\rho( \overline{ \Gamma})$ is contained in the centralizer of this 
  element hence in the subgroup of diagonal matrices. Hence $\rho$ is trivial 
  on the commutator subgroup $[ \overline{ \Gamma}, \overline{ \Gamma}]$, 
  but $\tau^{2g}$ is in this group with $\rho( \tau)^{2g} \neq \textup{Id}$.
\end{proof}

\subsection{Convex sets}
\label{sec:convsetinp2}

\begin{lem}\label{lem_convexset}
Let $I$ be an interval and $\Ll: I \to  \PP^2(\RR)^*$ a continuous map. 
Then the subset 
\bqn
D = \PP^2(\RR) - \bigcup_{t\in I}\Ll(t) 
\eqn
is a convex set in $\PP^2(\RR)$.
\end{lem}
\begin{proof}
Suppose that $D$ is not empty. 
Assume that $D$ is not convex. Then there exists $p,q \in D$, $p\neq q$ such that 
the intersection of the line $L$ spanned by $p$ and $q$ in $\PP^2(\RR)$ with $D$ is not connected. 
Therefore there are $t_1 <t_2 \in I$ such that $\Ll(t_1)$ and $\Ll(t_2)$ intersect the two different connected components of $L -\{p,q\}$.

We choose coordinates $(x,y)$ in the affine chart $\PP^2(\RR) - L$ such that 
$\Ll(t_1)$ is the $y$-axis $\{(x,y)\mid x=0\}$ and  $\Ll(t_2)$ is the $x$-axis $\{ (x,y) \mid y=0\}$.
Since $p,q \in D$ we have for all $t\in I$ that $\Ll(t) \neq L$, and there exist continuous functions 
$a,b,c : I \to \RR$ such that 
\bqn
\Ll(t) = \{(x,y) \mid a(t) x + b(t) y + c(t) = 0\}.  
\eqn
In particular $a(t_1) \neq 0$, $b(t_1) = c(t_1) = 0$ and $a(t_2) = c(t_2) = 0$, $b(t_2) \neq 0$.

The points $p,q$ in $ L = \PP(\{(x,y)\})$ have homogeneous coordinates $[x_p, y_p]$ and $[x_q, y_q]$, and since $p$ and $q$ do not belong to $\Ll(t_1)$ or $\Ll(t_2)$ the 
products $x_py_p \neq 0$ and $x_qy_q \neq 0$ are nonzero. Moreover since  $\Ll(t_1)$ and $\Ll(t_2)$ 
intersect $L -\{p,q\}$ in different connected components we have 
\bq\label{eq_different_sign}
\mathrm{sign}(x_py_p) \neq \mathrm{sign}(x_qy_q).
\eq
Furthermore, since $\Ll(t)$ does contain neither $p$ nor $q$, 
\bq\label{eq_nonzero}
a(t) x_p + b(t) y_p \neq  0 \,\text{ and }\,
a(t) x_q + b(t) y_q \neq 0.
\eq
In particular the quantities in Equation~\eqref{eq_nonzero} do not change their sign for all 
$t\in I$. 
Therefore 
\bqn
\mathrm{sign}\left(a(t_1) x_p\right)= \mathrm{sign}\left(b(t_2) y_p\right)\quad \text{and} \quad \mathrm{sign}\left(a(t_1) x_q\right) =\mathrm{sign}\left(b(t_2) y_q\right), 
\eqn
which implies 
\bqn
\mathrm{sign}(x_py_p) = \mathrm{sign}\left(a(t_1) b(t_2)\right) = \mathrm{sign}\left(x_qy_q\right), 
\eqn
contradicting Equation~\eqref{eq_different_sign}.
\end{proof}

\providecommand{\bysame}{\leavevmode\hbox to3em{\hrulefill}\thinspace}
\providecommand{\MR}{\relax\ifhmode\unskip\space\fi MR }
\providecommand{\MRhref}[2]{%
  \href{http://www.ams.org/mathscinet-getitem?mr=#1}{#2}
}
\providecommand{\href}[2]{#2}

\end{document}